\newcommand\SEKImasterusepackages{\usepackage[T1]{fontenc}}
\newcommand\SEKIusersusepackages
{
\usepackage{tikz}
\usepackage{amssymb}
\usepackage{named}
\usepackage{makeidx}
\bibliographystyle{named}
\makeindex
\input headerhot
\mathcommand\ident[1]{\mathsf{#1}}
\newcommand\plussymbol  {\ident{+}}
\newcommand\minussymbol {\ident{-}}
\newcommand\dividesymbol{\ident{/}}
\newcommand\timessymbol {\ident{*}}

\newcommand\set     {\ident{set}}

\newcommand\naturalssymbol{\ident{naturals}}
\newcommand\gensymsymbol{\ident{gensym}}
\mathcommand\mbpsymbol{\ident{m\hspace{-0.055em}b\hspace{-0.045em}p}}

\newcommand\csymbol     {\ident c}
\newcommand\esymbol     {\ident e}
\newcommand\fsymbol     {\ident f}
\newcommand\gsymbol     {\ident g}
\newcommand\hsymbol     {\ident h}
\newcommand\ksymbol     {\ident k}
\newcommand\psymbol     {\ident p}
\newcommand\ssymbol     {\ident s}
\newcommand\Everysymbol {\ident{Every}}
\newcommand\Permsymbol {\ident{Perm}}
\newcommand\RExistssymbol{\ident{Rexists}}
\newcommand\invertsymbol{\ident{invert}}
\newcommand\invsymbol{\ident{inv}}
\newcommand\abssymbol   {\ident{abs}}
\newcommand\cnssymbol   {\ident{cons}}
\mathcommand\cnsindexsymbol[1]{\ident{cons}_{#1}}
\newcommand\carsymbol   {\ident{car}}

\newcommand\cdrsymbol   {\ident{cdr}}
\newcommand\lengthsymbol{\ident{length}}
\newcommand\sizesymbol{\ident{size}}
\newcommand\dlsymbol    {\ident{dl}}
\newcommand\dloncesymbol{\ident{delfirst}}
\newcommand\rcsymbol    {\ident{rc}}
\newcommand\brsymbol    {\ident{br}}
\newcommand\revtailsymbol{\ident{revtail}}
\newcommand\revsymbol{\ident{rev}}
\newcommand\appendsymbol {\ident{append}}
\newcommand\zeropredicatesymbol{\ident{zerop}}
\newcommand\eqsymbol        {\ident{eq}}
\newcommand\ifthensymbol    {\mbox{\ident{If{}Then}}}
\newcommand\ifthenelsesymbol{\mbox{\ident{If{}ThenElse}}}
\mathcommand\eqindexsymbol        [1]{\eqsymbol        _{#1}}
\mathcommand\ifthenindexsymbol    [1]{\ifthensymbol    _{#1}}
\mathcommand\ifthenelseindexsymbol[1]{\ifthenelsesymbol_{#1}}
\newcommand\orsymbol    {\ident{or}}
\newcommand\andsymbol   {\ident{and}}
\newcommand\leqsymbol   {\ident{leq}}
\newcommand\lessymbol   {\ident{less}}
\newcommand\lexlessymbol{\ident{lexless}}
\newcommand\lexlimlessymbol{\ident{lexlimless}}
\newcommand\lexsymbol   {\ident{lex}}
\newcommand\acksymbol   {\ident{ack}}
\newcommand\switchsymbol{\ident{switch}}
\newcommand\swatchsymbol{\ident{swatch}}
\newcommand\diveinssymbol{\ident{div1}}
\newcommand\divzweisymbol{\ident{div2}}
\newcommand\divrestsymbol{\ident{divrest}}
\newcommand\diveinstailsymbol{\ident{div1tail}}
\newcommand\divzweitailsymbol{\ident{div2tail}}
\newcommand\remsymbol{\ident{rem}}
\newcommand\divsymbol{\ident{div}}

\newcommand\turingmachinesymbol{\ident T}
\newcommand\terminatespsymbol  {\ident{terminatesp}}
\newcommand\statesymbol        {\ident{state}}
\newcommand\cmdsymbol          {\ident{cmd}}
\newcommand\nthsymbol          {\ident{nth}}
\newcommand\doublesymbol       {\ident{double}}

\newcommand\ppsymbol           {\ident{p}}
\newcommand\qpsymbol           {\ident{q}}
\newcommand\Epsymbol           {\ident{E}}
\newcommand\Ppsymbol           {\ident{P}}
\newcommand\Qpsymbol           {\ident{Q}}
\newcommand\Marriessymbol      {\ident{Marries}}
\newcommand\Lovessymbol        {\ident{Loves}}
\newcommand\StolenBysymbol     {\ident{StolenBy}}
\newcommand\Humansymbol        {\ident{Human}}
\newcommand\Evensymbol         {\ident{Even}}
\newcommand\Oddsymbol          {\ident{Odd}}
\newcommand\Primesymbol        {\ident{Prime}}
\newcommand\EveryPairsymbol   {\ident{EveryPair}}
\newcommand\Givesymbol         {\ident{Give}}
\newcommand\Fathersymbol       {\ident{Father}}
\newcommand\Elephantpsymbol    {\ident{Elephant}}
\newcommand\Flowerpsymbol    {\ident{Flower}}
\newcommand\Germanpsymbol      {\ident{German}}
\newcommand\Bicyclepsymbol     {\ident{Bicycle}}
\newcommand\Hugepsymbol        {\ident{Huge}}
\newcommand\Animalpsymbol      {\ident{Animal}}
\newcommand\Malepsymbol        {\ident{Male}}
\newcommand\Boypsymbol         {\ident{Boy}}
\newcommand\Girlpsymbol        {\ident{Girl}}
\newcommand\Femalepsymbol      {\ident{Female}}
\newcommand\Roundpsymbol       {\ident{Round}}
\newcommand\Quadrangularpsymbol{\ident{Quadrangular}}
\newcommand\Metpsymbol         {\ident{Met}}
\newcommand\Kissedpsymbol      {\ident{Kissed}}
\newcommand\Bishopsymbol       {\ident{Bishop}}
\newcommand\mindexsymbol[1]{\existsvari w{#1}}

\newcommand\nonnegpsymbol      {\ident{nonnegp}}
\newcommand\wellsymbol         {\ident{well}}
\newcommand\welltailsymbol     {\ident{welltail}}
\newcommand\varsymbol          {\ident{var}}
\newcommand\aritysymbol        {\ident{arity}}

\newcommand\whilesymbol        {\ident{while}}

\newcommand\nullsymbol         {\ident{null}}
\newcommand\hdsymbol           {\ident{hd}}
\newcommand\tlsymbol           {\ident{tl}}
\newcommand\insymbol           {\ident{in}}
\newcommand\applysymbol        {\ident{app}}
\newcommand\termsymbol         {\ident{term}}
\newcommand\russellsymbol      {\ident{russell}}
\newcommand\sqrtindordsymbol[1]{\ident{sqrtio#1}}
\mathcommand\tightim{\longrightarrow}
\mathcommand\im{\ \tightim\ }
\mathcommand\rs{\:\rulesugar\:\:}
\mathcommand\rulesugar{\longleftarrow}

\mathcommand\doublepp[1]      {\doublesymbol   \beginargs{#1}\allargs}
\mathcommand\aritypp[1]      {\aritysymbol   \beginargs{#1}\allargs}
\mathcommand\lengthpp[1]      {\lengthsymbol   \beginargs{#1}\allargs}
\mathcommand\sizepp[1]      {\sizesymbol   \beginargs{#1}\allargs}
\mathcommand\wellpp[1]      {\wellsymbol   \beginargs{#1}\allargs}
\mathcommand\welltailpp[1]      {\welltailsymbol   \beginargs{#1}\allargs}
\mathcommand\varpp[1]      {\varsymbol   \beginargs{#1}\allargs}
\mathcommand\rempp[2]    {\remsymbol\beginargs{#1}\separgs{#2}\allargs}
\mathcommand\divpp[2]    {\divsymbol\beginargs{#1}\separgs{#2}\allargs}
\mathcommand\divrestpp[2]    {\divrestsymbol\beginargs{#1}\separgs{#2}\allargs}
\mathcommand\diveinspp[2]    {\diveinssymbol\beginargs{#1}\separgs{#2}\allargs}
\mathcommand\divzweipp[3]    {\divzweisymbol\beginargs{#1}\separgs{#2}
\separgs{#3}\allargs}
\mathcommand\diveinstailpp[4]    {\diveinstailsymbol\beginargs{#1}\separgs{#2}
\separgs{#3}\separgs{#4}\allargs}
\mathcommand\divzweitailpp[6]    {\divzweitailsymbol\beginargs{#1}\separgs{#2}
\separgs{#3}\separgs{#4}\separgs{#5}\separgs{#6}\allargs}
\mathcommand\mbppp[2]         {\mbpsymbol   \beginargs{#1}\separgs{#2}\allargs}
\mathcommand\revpp[1]     
{\revsymbol\beginargs{#1}\allargs}
\mathcommand\revppi[2]     
{\revsymbol^{#1}\beginargs{#2}\allargs}
\mathcommand\revtailpp[2]     
{\revtailsymbol\beginargs{#1}\separgs{#2}\allargs}
\mathcommand\revtailppi[3]
{\revtailsymbol^{#1}\beginargs{#2}\separgs{#3}\allargs}
\mathcommand\Permpp[2]     
{\Permsymbol\beginargs{#1}\separgs{#2}\allargs}
\mathcommand\Permppi[3]
{\Permsymbol^{#1}\beginargs{#2}\separgs{#3}\allargs}
\mathcommand\appendpp[2]      
{\appendsymbol \beginargs{#1}\separgs{#2}\allargs}
\mathcommand\appendppi[3]      
{\appendsymbol^{#1}\beginargs{#2}\separgs{#3}\allargs}
\mathcommand\Everypp[2]      
{\Everysymbol \beginargs{#1}\separgs{#2}\allargs}
\mathcommand\RExistspp[1]      
{\RExistssymbol \beginargs{#1}\allargs}
\mathcommand\appendlongpp[2]      
{\appendsymbol\left(\begin{array}{@{}l@{}}{#1}\separgs\\{#2}\end{array}\right)}
\mathcommand\cnspp[2]         {\cnssymbol   \beginargs{#1}\separgs{#2}\allargs}
\mathcommand\cnsppi[3]       {\cnssymbol^{#1}\beginargs{#2}\separgs{#3}\allargs}
\mathcommand\cnsindexpp[3]
{\cnsindexsymbol{#1}\beginargs{#2}\separgs{#3}\allargs}
\mathcommand\dlpp[2]          {\dlsymbol    \beginargs{#1}\separgs{#2}\allargs}
\mathcommand\dloncepp[2]      {\dloncesymbol\beginargs{#1}\separgs{#2}\allargs}
\mathcommand\dlonceppi[3]{\dloncesymbol^{#1}\beginargs{#2}\separgs{#3}\allargs}
\mathcommand\rcpp[2]          {\rcsymbol    \beginargs{#1}\separgs{#2}\allargs}
\mathcommand\brpp[2]          {\brsymbol    \beginargs{#1}\separgs{#2}\allargs}
\mathcommand\orpp[2]          {\orsymbol    \beginargs{#1}\separgs{#2}\allargs}
\mathcommand\andpp[2]         {\andsymbol   \beginargs{#1}\separgs{#2}\allargs}
\mathcommand\shortcnspp[2]    {\csymbol     \beginargs{#1}\separgs{#2}\allargs}
\mathcommand\tightshortcnspp[2]
{\csymbol\beginargs{#1}\tightsepargs{#2}\allargs}
\mathcommand\spp[1]           {\ssymbol     \beginargs{#1}\allargs}
\mathcommand\sppiterated[2]   {\ssymbol^{#1}\beginargs{#2}\allargs}
\mathcommand\sqrtindordpp[3]
                       {\sqrtindordsymbol{#1}\beginargs{#2}\separgs{#3}\allargs}
\mathcommand\ppp[1]           {\psymbol     \beginargs{#1}\allargs}
\mathcommand\pppiterated[2]   {\psymbol^{#1}\beginargs{#2}\allargs}
\mathcommand\zeropp           {\ident 0}
\mathcommand\Julietpp         {\ident{Juliet}}
\mathcommand\Romeopp          {\ident{Romeo}}
\mathcommand\Ipp              {\ident I}
\mathcommand\onepp            {\ident1}
\mathcommand\twopp            {\ident2}
\mathcommand\threepp          {\ident3}
\mathcommand\invertpp[1]      {\invertsymbol\beginargs{#1}\allargs}
\mathcommand\invpp[1]         {\invsymbol\beginargs{#1}\allargs}
\mathcommand\abspp[1]         {\abssymbol\beginargs{#1}\allargs}
\mathcommand\naturalspp[1]    {\naturalssymbol\beginargs{#1}\allargs}
\mathcommand\gensympp[1]      {\gensymsymbol\beginargs{#1}\allargs}
\mathcommand\nilpp            {\ident{nil}}
\mathcommand\falsepp          {\ident{false}}
\mathcommand\truepp           {\ident{true}}
\mathcommand\FALSEpp          {\ident{FALSE}}
\mathcommand\TRUEpp           {\ident{TRUE}}
\mathcommand\UNDEFpp          {\ident{UNDEF}}
\mathcommand\weirdppp         {\ident{weirdp}}
\mathcommand\ambigppp         {\ident{ambigp}}
\mathcommand\zeropredicatepp[1]{\zeropredicatesymbol\beginargs{#1}\allargs}
\mathcommand\cppeins       [1]{\csymbol     \beginargs{#1}\allargs}
\mathcommand\cppzwei       [2]{\csymbol\beginargs{#1}\separgs{#2}\allargs}
\mathcommand\eppeins       [1]{\esymbol     \beginargs{#1}\allargs}
\mathcommand\fppeins       [1]{\fsymbol     \beginargs{#1}\allargs}
\mathcommand\fppeinsindex  [2]{\fsymbol_{#1}\beginargs{#2}\allargs}
\mathcommand\fppeinsiterated[2]{\fsymbol^{#1}\beginargs{#2}\allargs}
\mathcommand\gppeins       [1]{\gsymbol     \beginargs{#1}\allargs}
\mathcommand\gppzwei       [2]{\gsymbol     \beginargs{#1}\separgs{#2}\allargs}
\mathcommand\hppeins       [1]{\hsymbol     \beginargs{#1}\allargs}
\mathcommand\kppeins       [1]{\ksymbol     \beginargs{#1}\allargs}
\mathcommand\appzero          {\ident a}
\mathcommand\bppzero          {\ident b}
\mathcommand\cppzero          {\ident c}
\mathcommand\dppzero          {\ident d}
\mathcommand\eppzero          {\ident e}
\mathcommand\eqindexpp[3]{\eqindexsymbol{#1}\beginargs{#2}\separgs{#3}\allargs}
\mathcommand\ifthenindexpp
[3]{\ifthenindexsymbol{#1}\beginargs{#2}\separgs{#3}\allargs}
\mathcommand\ifthenelseindexpp
[4]{\ifthenelseindexsymbol{#1}\beginargs{#2}\separgs{#3}\separgs{#4}\allargs}
\mathcommand\eqpp[2]{\eqsymbol\beginargs{#1}\separgs{#2}\allargs}
\mathcommand\leqpp[2]{\leqsymbol\beginargs{#1}\separgs{#2}\allargs}
\mathcommand\lespp[2]{\lessymbol\beginargs{#1}\separgs{#2}\allargs}
\mathcommand\lexlespp[2]{\lexlessymbol\beginargs{#1}\separgs{#2}\allargs}
\mathcommand\lexlimlespp[3]
               {\lexlimlessymbol\beginargs{#1}\separgs{#2}\separgs{#3}\allargs}
\mathcommand\lexpp[3]{\lexsymbol\beginargs{#1}\separgs{#2}\separgs{#3}\allargs}
\mathcommand\ackpp[2]{\acksymbol\beginargs{#1}\separgs{#2}\allargs}
\mathcommand\switchpp[1]{\switchsymbol\beginargs{#1}\allargs}
\mathcommand\swatchpp[1]{\swatchsymbol\beginargs{#1}\allargs}
\mathcommand\whilepp[2]{\whilesymbol\beginargs{#1}\separgs{#2}\allargs}
\mathcommand\nullpp[1]{\nullsymbol\beginargs{#1}\allargs}
\mathcommand\nullppiterated[2]{\nullsymbol^{#1}\beginargs{#2}\allargs}
\mathcommand\hdpp[1]{\hdsymbol\beginargs{#1}\allargs}
\mathcommand\hdppiterated[2]{\hdsymbol^{#1}\beginargs{#2}\allargs}
\mathcommand\carpp[1]{\carsymbol\beginargs{#1}\allargs}
\mathcommand\cdrpp[1]{\cdrsymbol\beginargs{#1}\allargs}
\mathcommand\tlpp[1]{\tlsymbol\beginargs{#1}\allargs}
\mathcommand\tlppiterated[2]{\tlsymbol^{#1}\beginargs{#2}\allargs}
\mathcommand\inpp[2]{\insymbol\beginargs{#1}\separgs{#2}\allargs}
\mathcommand\inppiterated[3]{\insymbol^{#1}\beginargs{#2}\separgs{#3}\allargs}
\mathcommand\applypp[2]{\applysymbol\beginargs{#1}\separgs{#2}\allargs}
\mathcommand\applyppiterated
[3]{\applysymbol^{#1}\beginargs{#2}\separgs{#3}\allargs}
\mathcommand\termpp[2]{\termsymbol\beginargs{#1}\separgs{#2}\allargs}
\mathcommand\setpp[1]{\set\beginargs{#1}\allargs}
\mathcommand\russellpp[1]{\russellsymbol\beginargs{#1}\allargs}

\mathcommand\Tpp[6]{\turingmachinesymbol\beginargs{#1}\separgs{#2}\separgs
{#3}\separgs{#4}\separgs{#5}\separgs{#6}\allargs}
\mathcommand\Tppseven[7]{\turingmachinesymbol\beginargs{#1}\separgs{#2}\separgs
{#3}\separgs{#4}\separgs{#5}\separgs{#6}\separgs{#7}\allargs}
\mathcommand\foreverppp[6]{\ident{foreverp}\beginargs{#1}\separgs{#2}\separgs
{#3}\separgs{#4}\separgs{#5}\separgs{#6}\allargs}
\mathcommand\terminatesppp[6]{\terminatespsymbol\beginargs{#1}\separgs
{#2}\separgs{#3}\separgs{#4}\separgs{#5}\separgs{#6}\allargs}
\mathcommand\terminatespppone[1]{\terminatespsymbol \beginargs{#1}\allargs}
\mathcommand\statepp
[3]{\statesymbol\beginargs{#1}\separgs{#2}\separgs{#3}\allargs}
\mathcommand\tightstatepp
[3]{\statesymbol\beginargs{#1}\tightsepargs{#2}\tightsepargs{#3}\allargs}
\mathcommand\cmdpp
[3]{\cmdsymbol  \beginargs{#1}\separgs{#2}\separgs{#3}\allargs}
\mathcommand\tightcmdpp
[3]{\cmdsymbol  \beginargs{#1}\tightsepargs{#2}\tightsepargs{#3}\allargs}
\mathcommand\stoppp           {\ident{stop}}
\mathcommand\leftpp           {\ident{left}}
\mathcommand\rightpp          {\ident{right}}
\mathcommand\nthpp         [2]{\nthsymbol  \beginargs{#1}\separgs{#2}\allargs}
\mathcommand\pppp          [1]{\ppsymbol\beginargs{#1}            \allargs}
\mathcommand\qppp          [2]{\qpsymbol\beginargs{#1}\separgs{#2}\allargs}
\mathcommand\Eppp          [1]{\Epsymbol\beginargs{#1}            \allargs}
\mathcommand\Epppzwei      [2]{\Epsymbol\beginargs{#1}\separgs{#2}\allargs}
\mathcommand\Pppp          [1]{\Ppsymbol\beginargs{#1}            \allargs}
\mathcommand\Ppppeinsindex [2]{\Ppsymbol_{#1}\beginargs{#2}\allargs}
\mathcommand\Ppppdrei      
[3]{\Ppsymbol\beginargs{#1}\separgs{#2}\separgs{#3}\allargs}
\mathcommand\Ppppvier
[4]{\Ppsymbol\beginargs{#1}\separgs{#2}\separgs{#3}\separgs{#4}\allargs}
\mathcommand\Qppp          [2]{\Qpsymbol\beginargs{#1}\separgs{#2}\allargs}
\mathcommand\Qpppeins      [1]{\Qpsymbol\beginargs{#1}\allargs}
\mathcommand\Qpppeinsindex [2]{\Qpsymbol_{#1}\beginargs{#2}\allargs}
\mathcommand\Qpppdrei      
[3]{\Qpsymbol\beginargs{#1}\separgs{#2}\separgs{#3}\allargs}
\mathcommand\Fatherpp      [2]{\Fathersymbol\beginargs{#1}\separgs{#2}\allargs}
\mathcommand\Marriespp     [2]{\Marriessymbol\beginargs{#1}\separgs{#2}\allargs}
\mathcommand\Lovespp       [2]{\Lovessymbol\beginargs{#1}\separgs{#2}\allargs}
\mathcommand\StolenBypp    [2]
{\StolenBysymbol\beginargs{#1}\separgs{#2}\allargs}
\mathcommand\Humanpp       [1]{\Humansymbol\beginargs{#1}\allargs}
\mathcommand\Evenpp        [1]{\Evensymbol\beginargs{#1}\allargs}
\mathcommand\Evenppi       [2]{\Evensymbol^{#1}\beginargs{#2}\allargs}
\mathcommand\Oddpp         [1]{\Oddsymbol\beginargs{#1}\allargs}
\mathcommand\Primepp       [1]{\Primesymbol\beginargs{#1}\allargs}
\mathcommand\EveryPairpp  [2]{\EveryPairsymbol\beginargs{#1}\separgs
{#2}\allargs}
\mathcommand\mindexppeins  [2]{\mindexsymbol{#1}\beginargs{#2}\allargs}
\mathcommand\Givepp        [3]{\Givesymbol
\beginargs{#1}\separgs{#2}\separgs{#3}\allargs}
\mathcommand\mindexppzwei  [3]{\mindexsymbol
{#1}\beginargs{#2}\separgs{#3}\allargs}
\mathcommand\mindexppdrei  [4]{\mindexsymbol
{#1}\beginargs{#2}\separgs{#3}\separgs{#4}\allargs}

\mathcommand\nonnegppp     [1]{\nonnegpsymbol\beginargs{#1}\allargs}

\mathcommand\anonymouscsymbol{c}
\mathcommand\anonymouscindexsymbol[1]{\anonymouscsymbol_{#1}}
\mathcommand\anonymousfsymbol{f}
\mathcommand\anonymouscpp
[2]{\anonymouscsymbol\beginargs{#1}\separgs\ldots\separgs{#2}\allargs}
\mathcommand\anonymouscindexpp
[3]{\anonymouscindexsymbol{#1}\beginargs{#2}\separgs\ldots\separgs{#3}\allargs}
\mathcommand\anonymousfpp
[2]{\anonymousfsymbol\beginargs{#1}\separgs\ldots\separgs{#2}\allargs}
\mathcommand\coerceindexpp[3]{[#3]_{#1}^{#2}}

\mathcommand\Elephantppp    [1]{\Elephantpsymbol\beginargs{#1}\allargs}
\mathcommand\Flowerppp      [1]{\Flowerpsymbol  \beginargs{#1}\allargs}
\mathcommand\Bicycleppp     [1]{\Bicyclepsymbol \beginargs{#1}\allargs}
\mathcommand\Germanppp      [1]{\Germanpsymbol  \beginargs{#1}\allargs}
\mathcommand\Hugeppp        [1]{\Hugepsymbol    \beginargs{#1}\allargs}
\mathcommand\Animalppp      [1]{\Animalpsymbol  \beginargs{#1}\allargs}
\mathcommand\Maleppp        [1]{\Malepsymbol    \beginargs{#1}\allargs}
\mathcommand\Boyppp         [1]{\Boypsymbol     \beginargs{#1}\allargs}
\mathcommand\Girlppp        [1]{\Girlpsymbol    \beginargs{#1}\allargs}
\mathcommand\Femaleppp      [1]{\Femalepsymbol  \beginargs{#1}\allargs}
\mathcommand\Roundppp       [1]{\Roundpsymbol   \beginargs{#1}\allargs}
\mathcommand\Bishoppp       [1]{\Bishopsymbol   \beginargs{#1}\allargs}
\mathcommand\Quadrangularppp[1]{\Quadrangularpsymbol  \beginargs{#1}\allargs}
\mathcommand\Kissedppp[2]{\Kissedpsymbol\beginargs{#1}\separgs{#2}\allargs}
\mathcommand\Metppp[2]   {\Metpsymbol   \beginargs{#1}\separgs{#2}\allargs}

\newcommand\bound     {{\rm bound}}
\newcommand\free      {{\rm free}}

\mathcommand\Vtripleindex[3]{\V\!_{{#1},\,{#2},\,{#3}}}
\mathcommand\Vdoubleindex[2]{\V\!_{{#1},\,{#2}}}
\mathcommand\Vsingleindex[1]{\V\!_{{#1}}}

\mathcommand\Erel[1]{\Gammaoffont\!_{#1}}
\mathcommand\Urel[1]{\Deltaoffont_{#1}}

\mathcommand\theRprimefromstrongtoweak{
  \inparenthesesinlinetight{
     \domres\id{\Vwall\cup\Vsome\setminus\RAN\varsigma}
     \nottight{\nottight\uplus}
     \reverserelation\varsigma
  }
  \nottight{\circ}
  \ranres
    {\transclosureinline R}
    {\Vwall\cup\Vsome\setminus\RAN\varsigma}
  \nottight{\nottight{\nottight{\uplus}}}
  \Vsome\tighttimes\Vsall
}

\mathcommand\deltaminus{\delta^-}
\mathcommand\deltaplus{\delta^+}
\mathcommand\deltaplusplus{\delta^{+^+}}
\mathcommand\deltastar{\delta^*}
\mathcommand\deltastarstar{\delta^{*^*}}

\mathcommand\Vall     {\Vsingleindex\indexdelta         }
\mathcommand\Vwall    {\Vsingleindex\indexdeltaminu     }
\mathcommand\Vsall    {\Vsingleindex\indexdeltaplus     }
\mathcommand\Vgsome   {\Vsingleindex\indexgammaplus     }
\mathcommand\Vsome    {\Vsingleindex\indexgamma         }
\mathcommand\Vfree    {\Vsingleindex\indexfree          }
\mathcommand\Vbound   {\Vsingleindex\indexbound         }
\mathcommand\Vsomesall{\Vsingleindex\indexgammadeltaplus}

\mathapplycommand\VARall      {\VARsingleindex\indexdelta         }
\mathapplycommand\VARwall     {\VARsingleindex\indexdeltaminu     }
\mathapplycommand\VARsall     {\VARsingleindex\indexdeltaplus     }
\mathapplycommand\VARgsome    {\VARsingleindex\indexgammaplus     }
\mathapplycommand\VARsome     {\VARsingleindex\indexgamma         }
\mathapplycommand\VARfree     {\VARsingleindex\indexfree          }
\mathapplycommand\VARbound    {\VARsingleindex\indexbound         }
\mathapplycommand\VARsomesall {\VARsingleindex\indexgammadeltaplus}
\mathcommand\displayVARsall[1]{\VARsingleindex\indexdeltaplus
\!\!\!\:\left(\begin{array}{@{}c@{}}#1\end{array}\right)}

\mathcommand\rigidvari     [2]{#1_{#2}^\indexgammadeltaplus}
\mathcommand\existsvari    [2]{#1_{#2}^\indexgamma    }
\mathcommand\forallvari    [2]{#1_{#2}^\indexdelta    }
\mathcommand\freevari      [2]{#1_{#2}^\indexfree     }
\mathcommand\wforallvari   [2]{#1_{#2}^\indexdeltaminu}
\mathcommand\sforallvari   [2]{#1_{#2}^\indexdeltaplus}
\mathcommand\gexistsvari   [2]{#1_{#2}^\indexgammaplus}
\mathcommand\boundvari     [2]{#1_{#2}}
\mathcommand\vari          [2]{#1_{#2}}
\mathcommand\wforallvarilow[2]{#1_{#2}^
{\raisebox{-.82ex}{\math\indexdeltaminu}}}

\newcommand\indexhelper[1]{{\scriptscriptstyle#1\:\!\!}}
\newcommand\indexdeltaplus
{\indexhelper{\delta^{\raisebox{-.17ex}{\fvesf\hskip-0.14em +}}}}
\newcommand\indexdeltaminu
{\indexhelper{\delta^{\mbox{\fvesf\hskip-0.14em\rule[.2ex]{.7em}{.15ex}}}}}
\newcommand\indexgammaplus
{\indexhelper{\gamma^{\mbox{\fvesf\hskip-0.14em +}}}}
\newcommand\indexgammadeltaplus
{\indexhelper{\gamma\delta^{\raisebox{-.17ex}{\fvesf\hskip-0.14em +}}}}

\newcommand\indexdelta{\indexhelper\delta}
\newcommand\indexgamma{\indexhelper\gamma}
\newcommand\indexfree
{{\scriptscriptstyle\free}}
\newcommand\indexbound
{{\scriptscriptstyle\bound}}

\newcommand\Wellfsymb{\ident{Wellf}}
\mathapplycommand\Wellfpp{\Wellfsymb}

\mathcommand\beginargs{(}
\mathcommand\allargs  {)}
\mathcommand\separgs  {,\,}
\mathcommand\tightsepargs{,}

\mathcommand\minusppnoparentheses  [2]{{#1}\,\minussymbol\,{#2}}
\mathcommand\tightminusppnoparentheses  [2]{{#1}\minussymbol{#2}}
\mathcommand\divideppnoparentheses [2]{{#1}\,\dividesymbol\,{#2}}
\mathcommand\plusppnoparentheses   [2]{{#1}\,\plussymbol \,{#2}}
\mathcommand\plusppnoparenthesesi  [3]{{#2}\,\plussymbol^{#1}\,{#3}}
\mathcommand\tightplusppnoparentheses   [2]{{#1}\plussymbol{#2}}
\mathcommand\timesppnoparentheses  [2]{{#1}\,\timessymbol\,{#2}}
\mathcommand\undppnoparentheses    [2]{{#1}\und            {#2}}
\mathcommand\oderppnoparentheses   [2]{{#1}\oder           {#2}}
\mathcommand\impliesppnoparentheses[2]{{#1}\implies        {#2}}
\mathcommand\leqinfixppnoparentheses[2]{{#1}\,\tight\leq\,{#2}}
\mathcommand\geqinfixppnoparentheses[2]{{#1}\,\tight\geq\,{#2}}
\mathcommand\dividepp [2]{(\divideppnoparentheses {#1}{#2})}
\mathcommand\minuspp  [2]{(\minusppnoparentheses  {#1}{#2})}
\mathcommand\pluspp   [2]{(\plusppnoparentheses   {#1}{#2})}
\mathcommand\timespp  [2]{(\timesppnoparentheses  {#1}{#2})}
\mathcommand\undpp    [2]{(\undppnoparentheses    {#1}{#2})}
\mathcommand\oderpp   [2]{(\oderppnoparentheses   {#1}{#2})}
\mathcommand\impliespp[2]{(\impliesppnoparentheses{#1}{#2})}

\input headertheorem
\def\citep{\cite}
\def\citet#1{\citeauthor{#1} \shortcite{#1}}
\newcommand\startcite{{\raise.2ex\hbox{[}}}
\newcommand\stopcite {\raise.2ex\hbox{]}}
\newcommand\citehelper[1]{\startcite #1\stopcite}
\newcommand\makeaciteoftwo[2]
{\citehelper{\citeauthor{#1}, \citeyear{#1}; \citeyear{#2}}}

\newcommand\makeaciteofthree[3]
{\citehelper{\citeauthor{#1}, \citeyear{#1}; \citeyear{#2}; \citeyear{#3}}}

}

\newcommand\SEKIREVIEWSTATUS{2}

\newcommand\SEKIREVIEWER
{{\sc Erica Melis} 
\\FR Informatik, Universit\"at des Saarlandes, D--66123 Saarbr\"ucken, Germany
\\E-mail: {\tt melis@activemath.org}
\\WWW: \url{http://www.activemath.org/~melis}
}

\newcommand\SEKIREFEREEPROCESS
{The previous version \cite{wirth-heijenoort} passed the review process of 
 \logicauniversalisname.}

\newcommand\SEKIDOCUMENTCLASS{0}

\newcommand\SEKIYEAR{2014}
\newcommand\SEKINUMBEROFYEAR{01}
\newcommand\SEKITYPE{0}
\newcommand\SEKIPUBLISHEDTYPE{3}
\newcommand\SEKIPUBLISHEDAS
{\cite{wirth-heijenoort}.}

\title{\herbrandsfundamentaltheorem:
\\The Historical Facts and
\\their Streamlining}
\author{\wirthname\\\Institute\\\emailcp}
\date{\small\May\,27, 2014\\Thorough Revision \Aug\,14, 2014}

\newcommand\frenchfont{\sf}
\newcommand\germanfont{}
\newcommand\germanfontfootnote{}
\newcommand\repair{correction}
\newcommand\repairverb{correct}
\newcommand\repaired{corrected}
\newcommand\heijenoortsrepairindex{\index{Heijenoort's correction}}
\newcommand\heijenoortsrepair{\heijenoortsrepairindex\heijenoort's \repair}
\newcommand\bernaysrepairindex{\index{Bernays' correction}}
\newcommand\bernaysrepair{\bernaysrepairindex\bernays' \repair}
\newcommand\goedelsrepairindex
  {\index{Goedel's and Dreben's correction@G{\"o}del's and Dreben's correction}}
\newcommand\goedelsrepair{\goedelsrepairindex\goedel's \repair}
\newcommand\goedelsanddrebensrepair
                             {\goedelsrepairindex\goedel's and \dreben's \repair}
\newcommand\termeins{\app{\forallvari m{}}{\forallvari v{},\forallvari w{}}}%
\newcommand\termzwei{\app{\forallvari m{}}{\forallvari u{},\termeins}}%
\newcommand\boxeins{
{\termeins}}%
\newcommand\boxu{
{\forallvari u{}}}%
\newcommand\boxv{
{\forallvari v{}}}%
\newcommand\boxw{
{\forallvari w{}}}%
\newcommand\termsofdepth[2]{\app{{\mathcal T}_{#1}}{#2}}
\newcommand\propertyCstar
{Property\nolinebreak\hskip.2em\nolinebreak C\/\math{^*}} 
\newcommand\PropertyCstar
{Property\nolinebreak\hskip.2em\nolinebreak C\/\math{^*}} 
 
\newcommand\propertyC{Property\nolinebreak\hskip.2em\nolinebreak C} 
\newcommand\PropertyC{Property\nolinebreak\hskip.2em\nolinebreak C} 
\newcommand\qed
{\relax\ifmmode\hskip2em \Box\else\unskip\nobreak\hskip1em $\Box$\fi}

\input seki-deckblatt-5

\usepackage{amssymb}

\usepackage[T1]{fontenc}

\newcommand\actualpagebreaknospace{{\rm\textbar}}
\newcommand\actualpagebreak{\actualpagebreaknospace\ }

\mathchardef\Gammaoffont="7000
\mathchardef\Gamma="0100
\mathchardef\Deltaoffont="7001
\mathchardef\Delta="0101
\mathchardef\Thetaoffont="7002
\mathchardef\Theta="0102
\mathchardef\Lambdaoffont="7003
\mathchardef\Lambda="0103
\mathchardef\Xioffont="7004
\mathchardef\Xi="0104
\mathchardef\Pioffont="7005
\mathchardef\Pi="0105
\mathchardef\Sigmaoffont="7006
\mathchardef\Sigma="0106
\mathchardef\Upsilonoffont="7007
\mathchardef\Upsilon="0107
\mathchardef\Phioffont="7008
\mathchardef\Phi="0108
\mathchardef\Psioffont="7009
\mathchardef\Psi="0109
\mathchardef\Omegaoffont="700A
\mathchardef\Omega="010A
\mathchardef\itype="017B

\catcode`\@=11

\gdef\allowhyphens{\penalty\@M \hskip\z@skip}

\gdef\set@low@box#1{\setbox\tw@\hbox{,}\setbox\z@\hbox{#1}\dimen\z@\ht\z@
     \advance\dimen\z@ -\ht\tw@
     \setbox\z@\hbox{\lower\dimen\z@ \box\z@}\ht\z@\ht\tw@ \dp\z@\dp\tw@ }
\gdef\set@low@boxsingle#1{\setbox\tw@\hbox{\rm,}\setbox\z@\hbox{#1}\dimen\z@\ht\z@
     \advance\dimen\z@ -\ht\tw@
     \setbox\z@\hbox{\lower\dimen\z@ \box\z@}\ht\z@\ht\tw@ \dp\z@\dp\tw@ }

\gdef\@glqq{%
\ifhmode\edef\@SF{\spacefactor\the\spacefactor}%
\else\let\@SF\empty
\fi
\CheckFamily\font\fraknomath\ifSameFamily ``\relax
\else\CheckFamily\font\swab\ifSameFamily ``\relax
\else\leavevmode\set@low@box{''}\box\z@\kern-.04em\allowhyphens\@SF\relax
\fi\fi}
\gdef\glqq{\protect\@glqq\kern+.07em}
\gdef\@grqq{%
\ifhmode\edef\@SF{\spacefactor\the\spacefactor}%
\else\let\@SF\empty 
\fi 
\CheckFamily\font\fraknomath\ifSameFamily ''\relax
\else\CheckFamily\font\swab\ifSameFamily ''\relax
\else\kern+.07em``\kern.07em\@SF\relax
\fi\fi}
\gdef\grqq{\protect\@grqq}
\gdef\@glq{{\ifhmode \edef\@SF{\spacefactor\the\spacefactor}\else
     \let\@SF\empty \fi \leavevmode
     \set@low@boxsingle{'\/}\box\z@\kern-.04em\allowhyphens\@SF\relax}}
\gdef\glq{\protect\@glq\kern+.07em}
\gdef\@grq{\ifhmode \edef\@SF{\spacefactor\the\spacefactor}\else
     \let\@SF\empty \fi \kern-.0125em`\kern.07em\@SF\relax}
\gdef\grq{\protect\@grq}

\catcode`\@=12

\newcommand\closequotecommanospace{''\nolinebreak\hskip-0.23em,}
\newcommand\closequotecomma      {\closequotecommanospace\         \,}

\newcommand\closequotefullstop   {\closequotefullstopnospace\      \,}
\newcommand\closequotefullstopextraspace   
                                 {\closequotefullstopnospace\    \ \,}

\newcommand\closequotefullstopnospace
                                 {''\nolinebreak\hskip-0.20em\@.}

\newcommand\closesinglequotefullstopnospace{'\nolinebreak\hskip-0.20em\@.}
\newcommand\closesinglequotefullstop{\closesinglequotefullstopnospace\ \,}

\makeatletter
\if \@ptsize 0
   \newfont{\scriptscriptscriptgoth}{ygoth scaled 760}
   \newfont{\scriptscriptgoth}{ygoth scaled 833}
   \newfont{\scriptgoth}{ygoth scaled 912}
   \newfont{\gothnomath}{ygoth}
   \newfont{\Goth}{ygoth scaled \magstephalf}
   \newfont{\GOth}{ygoth scaled \magstep1}
   \newfont{\GOTh}{ygoth scaled \magstep2}
   \newfont{\GOTH}{ygoth scaled \magstep3}

   \newfont{\scriptscriptscriptswab}{yswab scaled 760}
   \newfont{\scriptscriptswab}{yswab scaled 833}
   \newfont{\scriptswab}{yswab scaled 912}
   \newfont{\swab}{yswab}
   \newfont{\Swab}{yswab scaled \magstephalf}
   \newfont{\SWab}{yswab scaled \magstep1}
   \newfont{\SWAb}{yswab scaled \magstep2}
   \newfont{\SWAB}{yswab scaled \magstep3}

   \newfont{\scriptscriptscriptfrak}{yfrak scaled 760}
   \newfont{\scriptscriptfrak}{yfrak scaled 833}
   \newfont{\scriptfrak}{yfrak scaled 912}
   \newfont{\fraknomath}{yfrak}
   \newfont{\Frak}{yfrak scaled \magstephalf}
   \newfont{\FRak}{yfrak scaled \magstep1}
   \newfont{\FRAk}{yfrak scaled \magstep2}
   \newfont{\FRAK}{yfrak scaled \magstep3}

   \newfont{\init}{yinit}
   \newfont{\Init}{yinit scaled \magstephalf}
   \newfont{\INit}{yinit scaled \magstep1}
   \newfont{\INIt}{yinit scaled \magstep2}
   \newfont{\INIT}{yinit scaled \magstep3}
\fi
\if \@ptsize 1
   \newfont{\scriptscriptscriptgoth}{ygoth scaled 833}
   \newfont{\scriptscriptgoth}{ygoth scaled 912}
   \newfont{\scriptgoth}{ygoth}
   \newfont{\gothnomath}{ygoth scaled \magstephalf}
   \newfont{\Goth}{ygoth scaled \magstep1}
   \newfont{\GOth}{ygoth scaled \magstep2}
   \newfont{\GOTh}{ygoth scaled \magstep3}
   \newfont{\GOTH}{ygoth scaled \magstep4}

   \newfont{\scriptscriptscriptswab}{yswab scaled 833}
   \newfont{\scriptscriptswab}{yswab scaled 912}
   \newfont{\scriptswab}{yswab}
   \newfont{\swab}{yswab scaled \magstephalf}
   \newfont{\Swab}{yswab scaled \magstep1}
   \newfont{\SWab}{yswab scaled \magstep2}
   \newfont{\SWAb}{yswab scaled \magstep3}
   \newfont{\SWAB}{yswab scaled \magstep4}

   \newfont{\scriptscriptscriptfrak}{yfrak scaled 833}
   \newfont{\scriptscriptfrak}{yfrak scaled 912}
   \newfont{\scriptfrak}{yfrak}
   \newfont{\fraknomath}{yfrak scaled \magstephalf}
   \newfont{\Frak}{yfrak scaled \magstep1}
   \newfont{\FRak}{yfrak scaled \magstep2}
   \newfont{\FRAk}{yfrak scaled \magstep3}
   \newfont{\FRAK}{yfrak scaled \magstep4}

   \newfont{\init}{yinit scaled \magstephalf}
   \newfont{\Init}{yinit scaled \magstep1}
   \newfont{\INit}{yinit scaled \magstep2}
   \newfont{\INIt}{yinit scaled \magstep3}
   \newfont{\INIT}{yinit scaled \magstep4}
\fi
\if \@ptsize 2
   \newfont{\scriptscriptscriptgoth}{ygoth scaled 912}
   \newfont{\scriptscriptgoth}{ygoth}
   \newfont{\scriptgoth}{ygoth scaled \magstephalf}
   \newfont{\gothnomath}{ygoth scaled \magstep1}
   \newfont{\Goth}{ygoth scaled \magstep2}
   \newfont{\GOth}{ygoth scaled \magstep3}
   \newfont{\GOTh}{ygoth scaled \magstep4}
   \newfont{\GOTH}{ygoth scaled \magstep5}

   \newfont{\scriptscriptscriptswab}{yswab scaled 912}
   \newfont{\scriptscriptswab}{yswab}
   \newfont{\scriptswab}{yswab scaled \magstephalf}
   \newfont{\swab}{yswab scaled \magstep1}
   \newfont{\Swab}{yswab scaled \magstep2}
   \newfont{\SWab}{yswab scaled \magstep3}
   \newfont{\SWAb}{yswab scaled \magstep4}
   \newfont{\SWAB}{yswab scaled \magstep5}

   \newfont{\scriptscriptscriptfrak}{yfrak scaled 833}
   \newfont{\scriptscriptfrak}{yfrak}
   \newfont{\scriptfrak}{yfrak scaled \magstephalf}
   \newfont{\fraknomath}{yfrak scaled \magstep1}
   \newfont{\Frak}{yfrak scaled \magstep2}
   \newfont{\FRak}{yfrak scaled \magstep3}
   \newfont{\FRAk}{yfrak scaled \magstep4}
   \newfont{\FRAK}{yfrak scaled \magstep5}

   \newfont{\init}{yinit scaled \magstep1}
   \newfont{\Init}{yinit scaled \magstep2}
   \newfont{\INit}{yinit scaled \magstep3}
   \newfont{\INIt}{yinit scaled \magstep4}
   \newfont{\INIT}{yinit scaled \magstep5}
\fi

\newcommand{\mscriptscriptfrak}      [1]{\mbox{\scriptscriptscriptfrak#1}}
\newcommand{\mscriptfrak}            [1]{\mbox{\scriptscriptscriptfrak#1}}

\newcommand{\mfrak}[1]{\mbox{\fraknomath#1}}

\makeatother

\newif\ifSameFamily
\def\CheckFamily#1#2{\GetFamilyName{#1}\ArgOne
        \GetFamilyName{#2}\ArgTwo
        \ifx\ArgOne\ArgTwo\SameFamilytrue\else\SameFamilyfalse\fi}
\def\GetFamilyName#1{\edef\Tempa{#1}\def\Tempb{#1}\ifx\Tempa\Tempb
        \edef\Tempa{\fontname#1}\fi
        \edef\Tempa{\Tempa\space}%
        \expandafter\iGetFamilyName\Tempa\\}
\def\iGetFamilyName#1 #2\\#3{\def#3{#1}}
\def\DefFontName#1#2{{\escapechar-1\expandafter\expandafter\expandafter
        \iDefFontName\expandafter{\csname#2\endcsname}%
        \xdef#1{\expandafter\string\Tempa}}}
\def\iDefFontName{\def\Tempa}

\newcommand\unprotectedae
{\CheckFamily\font\fraknomath\ifSameFamily *a\else
 \CheckFamily\font\swab\ifSameFamily\char'212\else\"a\fi\fi}
\newcommand\unprotectedoe
{\CheckFamily\font\fraknomath\ifSameFamily 
*o\else\CheckFamily\font\swab\ifSameFamily\char'232\else\"o\fi\fi}
\newcommand\unprotectedue
{\CheckFamily\font\fraknomath\ifSameFamily 
*u\else\CheckFamily\font\swab\ifSameFamily\char'237\else\"u\fi\fi}

\DefFontName\eccclarge{eccc1200}
\DefFontName\eccc{eccc1000}
\DefFontName\ecccsmall{eccc0900}
\DefFontName\ecccfootnotesize{eccc0800}

\newcommand\unprotectedes
{\CheckFamily\font\fraknomath\ifSameFamily\char'215\else
\CheckFamily\font\swab\ifSameFamily\char'215\else  
s\fi\fi}

\newcommand\unprotectedesi
{\CheckFamily\font\fraknomath\ifSameFamily\char'215\else
\CheckFamily\font\swab\ifSameFamily\char'215\else  
\mbox{s\hskip.04em}\fi\fi
\-\ignorespaces}

\newcommand\unprotectedmyparagraphsymbol
{\CheckFamily\font\fraknomath\ifSameFamily 
\char'244\else\CheckFamily\font\swab\ifSameFamily
\char'244\else\S\fi\fi}

\renewcommand\ae{\protect\unprotectedae}
\renewcommand\oe{\protect\unprotectedoe}
\newcommand\ue  {\protect\unprotectedue}

\newcommand\es  {\protect\unprotectedes}

\newcommand\esi {\protect\unprotectedesi}  

\newcommand\myparagraphsymbol{\protect\unprotectedmyparagraphsymbol}

\def\mathfrak#1{%
\mathchoice
{{\mfrak{#1}}}
{{\mfrak{#1}}}
{{\mscriptfrak{#1}}}
{{\mscriptscriptfrak{#1}}}
}

\newcommand\namefont{}

\usepackage{url}

\newcommand\majorfootroom{\raisebox{-1.9ex}{\rule{0ex}{.5ex}}}

\newcommand\footroom{\raisebox{-1.5ex}{\rule{0ex}{.5ex}}}

\newcommand\smallfootroom{\raisebox{-1.0ex}{\rule{0ex}{3ex}}}

\newcommand\majorheadroom{\rule{0ex}{3.2ex}}
\newcommand\headroom{\rule{0ex}{2.8ex}}
\newcommand\mediumheadroom{\rule{0ex}{2.4ex}}

\newcommand\claus    {Clau\es}
\newcommand\david    {{\namefont David}}

\newcommand\gerhard  {Ger\-hard}

\newcommand\irving   {Irving}

\newcommand\jacques  {{\namefont Jacque\es}}

\newcommand\jean     {Jean}
\newcommand\john     {John}

\newcommand\kurt     {Kurt}

\newcommand\leopold  {Leo\-pold}

\newcommand\paul     {{\namefont Paul}}
\newcommand\peter    {Peter}

\newcommand\raymond  {Raymond}

\newcommand\thoralf  {{\namefont Thoralf}}

\newcommand\warren   {Warren}

\newcommand\aanderaa        {{\namefont Aanderaa}}
\newcommand\aanderaaname    {{\namefont St\aa l \aanderaa}}

\newcommand\anellisindex    {\index{Anellis, Irving H. (1946--2013)}}

\newcommand\anellis         {{\namefont Anelli\es}}
\newcommand\anellisname     {{\namefont\anellisindex\irving\ H. \anellis}}

\newcommand\anellislifetime {(1946--2013)}

\newcommand\andrewsindex    {\index{Andrews, Peter B. (*1937)}}

\newcommand\andrews         {\mbox{\namefont Andrew\es}}
\newcommand\andrewsname     {\andrewsindex{\namefont\peter\ B. \andrews}}
\newcommand\andrewsnamewithoutmiddleinitial
                            {\andrewsindex{\namefont\peter\ \andrews}}
\newcommand\andrewslifetime {(*1937)}

\newcommand\autexierindex   {\index{Autexier, Serge (*1971)}}

\newcommand\autexier        {\mbox{\namefont Autexier}}
\newcommand\autexiername    {\autexierindex{\namefont Serge \autexier}}

\newcommand\bernaysindex    {\index{Bernays, Paul (1888--1977)}}

\newcommand\bernays         {\mbox{\namefont Ber\-nay\es}}
\newcommand\bernaysname     {\bernaysindex{\namefont \paul\       \bernays}}

\newcommand\bernayslifetime {(1888--1977)}

\newcommand\churchindex     {\index{Church, Alonzo (1903--1995)}}

\newcommand\church          {{\namefont Church}}
\newcommand\churchname      {\churchindex{\namefont Alonzo \church}}
\newcommand\churchlifetime  {(1903--1995)}

\newcommand\denton          {{\namefont Denton}}
\newcommand\dentonname      {{\namefont\john\ \denton}}

\newcommand\drebenindex                 {\index{Dreben, Burton (1927--1999)}}

\newcommand\dreben          {{\namefont Dreben}}
\newcommand\drebenname      {\drebenindex{\namefont Burton \dreben}}
\newcommand\drebenlifetime  {(1927--1999)}

\newcommand\engels          {\mbox{\namefont Engels}}

\newcommand\fefermanindex   {\index{Feferman, Sol(omon) (*1928)}}

\newcommand\feferman        {{\namefont Fefer\-man}}
\newcommand\fefermanname    {\fefermanindex{\namefont Sol \feferman}}
\newcommand\fefermanlifetime{(*1928)}

\newcommand\fermat          {\mbox{\namefont Fermat}}

\newcommand\fermatbirthyear 
{160?}

\newcommand\fermatslasttheorem{\fermat's Last Theorem}

\newcommand\fourcolortheorem{Four-Color Theorem}

\newcommand\frege           {{\namefont Frege}}

\newcommand\Begriffsschrift {Begriff\esi schrift}

\newcommand\michaelgabbayindex
{}

\newcommand\gentzenindex    {\index{Gentzen, Gerhard (1909--1945)}}
\newcommand\gentzen         {{\namefont Gentzen}}
\newcommand\gentzenname     {\gentzenindex{\namefont\gerhard\ \gentzen}}
\newcommand\gentzenlifetime {(1909--1945)}
\newcommand\gentzensHauptsatz{{\namefont\gentzen}'\es\ {\Hauptsatz}}

\newcommand\Hauptsatz{Haupt\-satz}

\newcommand\goedelindex     {\index{Goedel, Kurt@G\"odel, Kurt (1906--1978)}}

\newcommand\goedel          {{\namefont G\oe del}}
\newcommand\goedelname      {\goedelindex{\namefont \kurt\ \goedel}}
\newcommand\goedellifetime  {(1906--1978)}

\newcommand\secondincompletenesstheorem
{second \incompletenesstheorem}
\newcommand\secondIncompletenessTheorem
{Second \IncompletenessTheorem}

\newcommand\firstincompletenesstheorem
{first \incompletenesstheorem}
\newcommand\firstIncompletenessTheorem
{First \IncompletenessTheorem}
\newcommand\incompletenesstheorem{incompleteness theorem}
\newcommand\IncompletenessTheorem{Incompleteness Theorem}

\newcommand\goldfarbindex   {\index{Goldfarb, Warren}}

\newcommand\goldfarb        {{\namefont Gold\-farb}}
\newcommand\goldfarbname    {\goldfarbindex{\namefont\warren\ 
                              \goldfarb}}

\newcommand\heijenoortindex {\index{Heijenoort, Jean van (1912--1986)}}

\newcommand\heijenoort      {{\namefont Heijen\-oort}}
\newcommand\vanheijenoort   {{\namefont van Heijen\-oort}}
\newcommand\heijenoortname  {\heijenoortindex{\namefont\jean\ 
                             \vanheijenoort}}
\newcommand\heijenoortdeathyear{1986}
\newcommand\heijenoortlifetime{(1912--\heijenoortdeathyear)}

\newcommand\henkin          {{\namefont Henkin}}

\newcommand\herbrandindex     {\index{Herbrand, Jacques (1908--1938)}}

\newcommand\herbrand        {{\namefont Herbrand}}
\newcommand\herbrandname    {\herbrandindex{\namefont \jacques\ \herbrand}}
\newcommand\herbranddeathyear{1931}
\newcommand\herbrandlifetime{(1908--\herbranddeathyear)}

\newcommand\herbrandsfundamentaltheoremindex
                                       {\index{Herbrand's Fundamental Theorem}}

\newcommand\herbrandsfundamentaltheoremnoindex
                                             {\herbrand'\es\ \fundamentaltheorem}
\newcommand\herbrandsfundamentaltheorem
           {\herbrandsfundamentaltheoremindex\herbrandsfundamentaltheoremnoindex}
\newcommand\herbrandsfalselemmaindex{\index{Herbrand's ``False Lemma''}}
\newcommand\herbrandsfalselemma
                           {\herbrandsfalselemmaindex\herbrand's ``False Lemma''}
\newcommand\herbrandsfalselemmacomma
             {\herbrandsfalselemmaindex\herbrand's ``False Lemma\closequotecomma}
\newcommand\herbrandsfalselemmafullstop
          {\herbrandsfalselemmaindex\herbrand's ``False Lemma\closequotefullstop}
\newcommand\herbrandsfalselemmawithinquotedtext
                           {\herbrandsfalselemmaindex\herbrand's `False Lemma'}

\newcommand\herbrandsfalselemmawithoutherbrandcomma
                         {\herbrandsfalselemmaindex``False Lemma\closequotecomma}
\newcommand\fundamentaltheorem{Fundamental Theorem}

\newcommand\hilbertindex                 {\index{Hilbert, David (1862--1943)}}

\newcommand\hilbert         {\mbox{\namefont Hilbert}}
\newcommand\hilbertname     {\hilbertindex{\namefont \david\ \hilbert}}
\newcommand\hilbertlifetime {(1862--1943)}

\newcommand\hilbertsepsilonlongindex{\index{Hilbert's epsilon}}

\newcommand\hilbertsepsilon
                     {\hilbertsepsilonlongindex\hilbert'\es\ \nlbmath\varepsilon}

\newcommand\hintikka        {{\namefont Hintikka}}

\newcommand\loewenheimindex {\index
{Loewenheim@L\"owenheim, Leopold (1878--1957)}}

\newcommand\loewenheim      {{\namefont L\oe wen\-heim}}
\newcommand\loewenheimname  {\loewenheimindex{\namefont\leopold\ \loewenheim}}
\newcommand\loewenheimlifetime{(1878--1957)}
\newcommand\loewenheimskolem{\loewenheim--\skolem}
\newcommand\loewenheimskolemtheorem{\index
                    {Loewenheim-Skolem Theorem@L{\oe}wenheim--Skolem Theorem}%
                                                     \loewenheimskolem\ Theorem}

\newcommand\marx            {{\namefont Marx}}

\newcommand\peirce          {{\namefont Peirce}}

\newcommand\russell         {{\namefont Russell}}

\newcommand\schroeder       {{\namefont Schr\oe der}}

\newcommand\schuetteindex   {\index{Schuette@Sch\"utte, Kurt (1909--1998)}}

\newcommand\schuette        {{\namefont Sch\ue tte}}
\newcommand\schuettename    {\schuetteindex{\namefont\kurt\ \schuette}}
\newcommand\schuettelifetime{(1909--1998)}

\newcommand\skolemindex     {\index{Skolem, Thoralf (1887--1963)}}

\newcommand\skolem          {{\namefont Skolem}}
\newcommand\skolemname      {\skolemindex\thoralf\ 
                             \skolem}
\newcommand\skolemlifetime  {(1887--1963)}

\newcommand\skolemization   {\skolem\-ization}
\newcommand\skolemizedform  {\skolem\-ized form}   
\newcommand\skolemizedForm  {\skolem\-ized Form}   
\newcommand\skolemnormalform
{\index{Skolem normal form}\skolem\ normal form}

\newcommand\smullyan        {{\namefont Smullyan}}
\newcommand\smullyanname    {{\namefont\raymond\ M. \smullyan}}
\newcommand\smullyanlifetime{(*1919)}

\newcommand\PM              {Principia Mathematica}

\newcommand\wirthindex                 {\index{Wirth, Claus-Peter (*1963)}}

\newcommand\wirth           {{\namefont Wirth}}
\newcommand\wirthnamenoindex{{\namefont\claus-\peter\ \wirth}}
\newcommand\wirthname       {\wirthindex\wirthnamenoindex}

\newcommand\Aug  {Aug.}

\newcommand\FB   {FB}

\newcommand\FBautinfveryshort{\FB\ AI}

\newcommand\f    {\mbox{}{f.}}   

\newcommand\freely{free}
\newcommand\defi {Definition} 

\newcommand\eMAIL{e-mail}
\newcommand\ifandonlyif{if \nolinebreak and \onlyif}

\newcommand\onlyif{only \nolinebreak if}

\newcommand\getittotheright[1]  
{\hfill\mbox{}\penalty 100\mbox{\ \,}\nolinebreak
\hfill\hfill\hfill\nolinebreak\mbox{#1}\ignorespaces}

\newcommand\role{r\^ole}

\newcommand\theo {Theorem}

\newcommand\Vol  {Vol.}
\newcommand\WWW  {WWW}

\newcommand\aswell{as \nolinebreak well}
\newcommand\aswellas{\aswell\ \nolinebreak as}

\newcommand\Cf   {Cf.}
\newcommand\cf   {cf.}

\newcommand\Cfnlb{\Cf\nolinebreak}
\newcommand\cfnlb{\cf\nolinebreak}

\newcommand\CS   {Computer \Sci}

\newcommand\eg   {e.g.}
\newcommand\Eg   {E.g.}

\newcommand\esp  {esp.}

\newcommand\sententialtautology{sentential tautology}
\newcommand\SententialTautology{Sentential Tautology}
\newcommand\sententialtautologies{sentential tautologies}

\newcommand\firstorder{first-order}

\newcommand\ie   {i.e.}

\newcommand\udiff{\ if\ }

\newcommand\May  {May}

\def\note{Note}

\newcommand\p    {p.}
\newcommand\pp   {pp.}
\newcommand\PP[2]{\pp\,\ignorespaces#1--\ignorespaces#2}

\newcommand\PhD  {PhD}
\newcommand\PhDthesis{\PhD\ thesis}
\newcommand\PhDThesis{\PhD\ Thesis}

\newcommand\sect {\myparagraphsymbol} 
\newcommand\sects{\myparagraphsymbol\myparagraphsymbol}

\newcommand\Sci  {Sci.}

\newcommand\singulary{singulary}

\newcommand\wrt  {w.r.t.}

\newcommand\thewordand{and}

\newcommand\litnoteref[1]{\note\,#1}

\newcommand\litfiguref[1]{Figure\,#1}

\newcommand\littheoref[1]{\theo\,#1}
\newcommand\litsectref[1]{\sect\,#1} 

\newcommand\litchapref[1]{Chapter\,#1} 

\newcommand\litsectfromtoref[2]{\sects\,#1--#2}

\newcommand\litdefiref[1]{\defi\,#1}
\newcommand\Examplename{Ex\-am\-ple}
\newcommand\litexamref[1]{\Examplename\,#1}
\newcommand\litcoexref[1]{Counter\-ex\-am\-ple\,#1}
\newcommand\litlemmref[1]{Lem\-ma\,#1}

\newcommand\litcororef[1]{Co\-rollary\,#1}

\newcommand\litlemmrefs[2]
{Lem\-mas        \nolinebreak #1 \thewordand\ \nolinebreak #2}

\newcommand\litcororefs[2]
{Co\-rol\-laries \nolinebreak #1 \thewordand\ \nolinebreak #2}

\newcommand\litsectrefs[2]
{\sects\         \nolinebreak #1 \thewordand\ \nolinebreak #2}

\newcommand\cororef[1]{\litcororef{\ref{#1}}}
\newcommand\lemmref[1]{\litlemmref{\ref{#1}}}

\newcommand\examref[1]{\litexamref{\ref{#1}}}
\newcommand\coexref[1]{\litcoexref{\ref{#1}}}
\newcommand\defiref[1]{\litdefiref{\ref{#1}}}
\newcommand\theoref[1]{\littheoref{\ref{#1}}}

\newcommand\figuref[1]{\litfiguref{\ref{#1}}}
\newcommand\noteref[1]{\litnoteref{\ref{#1}}}
\newcommand\sectref[1]{\litsectref{\ref{#1}}}
\newcommand\nlbsectref[1]{\nolinebreak\sectref{#1}}

\newcommand\lemmrefs[2]{\litlemmrefs{\ref{#1}}{\ref{#2}}}

\newcommand\cororefs[2]{\litcororefs{\ref{#1}}{\ref{#2}}}

\newcommand\sectrefs[2]{\litsectrefs{\ref{#1}}{\ref{#2}}}

\newcommand\sectfromtoref[2]{\litsectfromtoref{\ref{#1}}{\ref{#2}}}

\newcommand\nthpositioner[2]
{#1\raisebox{0.52ex}{\scriptsize\hspace{0.07em}#2}}
\newcommand\nth[1]{\nthtinypositioner{#1}{\nthstring{#1}}}
\newcommand\nthtinypositioner[2]{#1\raisebox{0.52ex}{\tiny\hspace{0.07em}#2}}

\newcommand\mthpositioner[2]
{\math{#1}\raisebox{0.52ex}{\scriptsize\hspace{0.07em}#2}}
\newcommand\modulointocountzero[2]
{\count1=#1
\count2=#2
\count0=\count1
\divide  \count0 by \count2
\multiply\count0 by-\count2
\advance \count0 by \count1}
\newcommand\absolutevalueintocountzero[1]
{\count0=#1
\ifnum\count0<0\multiply\count0 by -1\fi}
\newcommand\nthstring[1]
{\def\myargone{#1}\ifcat a\myargone th\else\nthstringnochar{#1}\fi}
\newcommand\nthstringnochar[1]
{\absolutevalueintocountzero{#1}%
\modulointocountzero{\count0}{100} 
\ifnum\count0>9\ifnum\count0<20 th\else\stupidnthstring\fi
                                  \else\stupidnthstring\fi}
\newcommand\stupidnthstring
{\modulointocountzero{\count0}{10}
\ifnum\count0=1 \hskip-0.2em st\else
\ifnum\count0=2 nd\else
\ifnum\count0=3 rd\else 
                th\fi\fi\fi}

\newcommand\writeascents
[1]{\count4=#1
\ifnum\count4<0 
-\multiply\count4 by -1\fi
\modulointocountzero{\count4}{10}
\divide\count4 by 10
\count3=\the\count0
\modulointocountzero{\count4}{10}
\divide\count4 by 10
\the\count4
.\the\count0
\the\count3
}

\newcommand\frenchnthstring[1]
{\def\myargone{#1}\ifcat a\myargone th\else\frenchnthstringnochar{#1}\fi}
\newcommand\frenchnthstringnochar[1]
{\absolutevalueintocountzero{#1}%
\modulointocountzero{\count0}{100} 
\ifnum\count0>9\ifnum\count0<20 th\else\frenchstupidnthstring\fi
                                  \else\frenchstupidnthstring\fi}
\newcommand\frenchstupidnthstring
{\modulointocountzero{\count0}{10}
\ifnum\count0=1 \hskip-0.2em re\else
\ifnum\count0=2 me\else
\ifnum\count0=3 rd\else 
                th\fi\fi\fi}

\newcommand\CLAM      {{\rm CL\kern-.36em\raise.39ex\hbox{\sc a}\kern-.15emM}}

\newcommand\TEXMACS   {{\sc T\kern-.1667em\lower.5ex\hbox{E}\kern-.125emX\kern-.1em\lower.5ex\hbox{\textsc{m\kern-.05ema\kern-.125emc\kern-.05ems}}}}

\newcommand\Wernigerode    {Werni\-gerode}

\newcommand\plzwernigerode{\mbox{38855}}

\def       \emailcp      {{\tt wirth@logic.at}}

\newcommand\Institutedept
{\FBautinfveryshort}
\newcommand\Instituteinst
{\mbox{Hochschule Harz}}
\newcommand\Instituteplac
{\plzwernigerode\,\Wernigerode}
\newcommand\Institutecoun{Germany}

\newcommand\Institute
{\Institutedept, \Instituteinst, \Instituteplac, \Institutecoun}

\newcommand\academicpress{Academic Press (\elsevier)}

\newcommand\elsevier{Elsevier}

\newcommand\newspaperreference[5]
{\def\nameofjournalpress{#2}#1, #4 #5, #3\if?\nameofjournalpress
\else, #2\fi}

\newcommand\dateinjournal[1]{}

\newcommand\journalreference[6]
{\def\nameofjournalpress{#2}#1\nolinebreak\hskip.2em%
\dateinjournal{(#3) }{\mbox{\bf #4}}, \PP{#5}{#6}\if?\nameofjournalpress
\else, #2\fi}

\newcommand\journalreferenceprintyear[6]
{\def\nameofjournalpress{#2}#1 
#4:#5--#6, #3%
\if?\nameofjournalpress
\else, #2\fi}

\newcommand\journalreferenceprintyearaspartofnumber[6]
{\def\nameofjournalpress{#2}#1 
{#4/#3}, \PP{#5}{#6}\if?\nameofjournalpress
\else, #2\fi}

\newcommand\jscname
{J. Symbolic Computation}

\newcommand\jscprintyear
{\journalreferenceprintyear{\jscname}\academicpress}

\newcommand\logicauniversalisname{Logica Universalis}

\newcommand\tcsname{Theoretical \CS}
\newcommand\tcsjournal
{\journalreference\tcsname\elsevier}
\newcommand\tcsjournalprintyear
{\journalreferenceprintyear\tcsname\elsevier}

\urldef\wirthkuhnurl\url
{http://wirth.bplaced.net/SEKI/welcome.html#SWP-2007-01}

\urldef\urlsrninetythreedashzerofive\url
{http://wirth.bplaced.net/SEKI/welcome.html#SR-93-05}

\urldef\urlsreightyeighttwelve\url
{http://wirth.bplaced.net/SEKI/welcome.html#SR-88-12}

\renewcommand\namefont{\sc}
\mathcommand\formulaexampleC
{\exists a\stopq\inparentheses
{\inpit{\neg\exists b\stopq
{\app p b}}\oder{\app p a}}}%
\mathcommand\formulaexampleCouterskolemizedform
{\exists a\stopq\inparentheses{
{\neg
{\app p{\app{\forallvari b{}}a}}}\oder{\app p a}}}
\newcommand\outerskolemizedformwithoutform
                       {\mbox{\deltaminus\hskip-.1em-\nolinebreak\skolem}\-ized}
\newcommand\outerskolemizedform          {\outerskolemizedformwithoutform\ form}
\newcommand\deltaplusplusskolemizedform{\math{\deltaplusplus\!}-\skolemizedform}
\begin{document}
\makecover
\maketitle
\begin{abstract}%
\sloppy
Using \heijenoort's unpublished 
{\em generalized}\/ rules of quantification,
we discuss the proof of \herbrandsfundamentaltheorem\ in the form
of \heijenoortsrepair\ of \mbox{\herbrandsfalselemma}
and present a didactic example.
Although we are mainly concerned with the inner structure of
\herbrandsfundamentaltheorem\ and the questions of its quality and its 
depth,
we also discuss the outer questions of its historical context and
why \bernays\ called it 
``the central theorem of predicate logic'' 
and considered the form of its expression to
be \mbox{``concise and felicitous''.}
\yestop\yestop\yestop
\end{abstract}
\tableofcontents\cleardoublepage
\section{Motivation}\label
{section Motivation}
\begin{sloppypar}
When 
working on our handbook article~\cite{herbrand-handbook}
on \herbrandindex\herbrand's work in logic,
we were fascinated by the following three points:
\begin{description}\item[\herbrandindex\herbrand's Unification Algorithm]%
\mbox{}\par\noindent
In addition to \heijenoortindex\heijenoort's 
reprint of \herbrandindex\herbrand's \PhDthesis\
in \cite{herbrand-ecrits-logiques}, \hskip.2em
on \nolinebreak which also the standard translation of the 1960s
by \dreben\ and \heijenoortindex\heijenoort\ in 
\cite{heijenoort-source-book}, \cite{herbrand-logical-writings}
is based, \hskip.1em
we accessed the original print~\cite{herbrand-PhD}. \hskip.3em
For some unknown reason,
\herbrandindex\herbrand's unification algorithm in the reprint~\cite{herbrand-ecrits-logiques}
has the confusing wording ``{\frenchfont\'egalit\'es associ\'ees}''\/
instead of the original 
``{\frenchfont\'egalit\'es normales}\closequotefullstopextraspace 
Based on the correct original text,
our English translation of 
\herbrandindex\herbrand's unification algorithm
resulted in a modern equational formulation,
which is more lucid than the one of \dreben\ and \heijenoortindex\heijenoort\
and strongly supports the thesis of~\cite{abeles-herbrand-unification}
that \herbrandindex\herbrand's unification algorithm is very close to 
the equational unification algorithm of~\cite{martelli-montanari}.\footnote{%
 Note that that thesis of~\cite{abeles-herbrand-unification}
 seems to have been developed
 without access to the uncontorted original~\cite{herbrand-PhD}. \hskip.4em
 While our translation may well be influenced by our knowledge of 
 the development of unification algorithms in the 1970s
 (see \cite{martelli-montanari}, \cite{PatersonWegman78} for references),
 \hskip.2em we noticed the article~\cite{abeles-herbrand-unification}
 only after our translation had already been completed.%
}\item
[\heijenoortsrepair\ of \herbrandsfalselemma]\label{item repair}\mbox{}%
\par\noindent We were surprised that \heijenoortsrepairindex\heijenoort's 
lucid and modern \repair\ of 
\herbrandsfalselemma\ in \cite{heijenoort-herbrand} was still 
unpublished.\item
[\herbrandindex\herbrand's Construction of Linear First-Order Derivations]\label
{item modus ponens elimination}\mbox{}%
\par\noindent
We were fascinated by that part of the proof of 
\herbrandsfundamentaltheorem\ which shows how to construct
a {\em linear}\/ derivation of a first-order formula \nlbmath A,
given as input (besides \nlbmath A) \hskip.1em
only the natural number \nlbmath n \hskip.1em
such that \math A has the {\em purely sentential}\/ 
\PropertyC\ of order \nlbmaths n.

Here, ``linear'' means 
that the derivation does not contain any inference step
with more than one 
premise 
(such as an application of \mbox{{\em modus ponens}\/),} \hskip.2em
such that it has no branching when it is viewed as a tree.
Moreover, all function and predicate symbols 
within this proof occur already in \nlbmath A,
and all formulas in the proof are similar to \nlbmath A \hskip.1em
in \nolinebreak the sense that they have the 
so-called {\em``sub-formula'' property}\/ \wrt\ \nlbmaths A.

As a consequence of this linearity,
\herbrandindex\herbrand's {\it modus ponens}\/-free calculus and 
\gentzen's classical sequent calculus 
represent first steps toward human-orientation
in the sense that a human being
with a semantical understanding of an informal proof has a good chance to 
construct a formal proof in these calculi.\footnote{
 \majorheadroom
 Calculi that are both human- and machine-oriented are
 in demand for the synergetic combination of
 mathematicians and computing machines in interactive theorem
 proving systems with strong automation support~\cite{wirth-manifesto-ljigpl}.}

\pagebreak

We \nolinebreak think that \herbrandindex\herbrand's construction of 
{\it modus ponens}\/-free proofs
should be part of the standard education of every logician,
just as well as the famous construction of Cut-free
proofs according to \gentzensHauptsatz\ \cite{gentzen},
wherein \gentzen\
---~some years later~---
picked up \herbrandindex\herbrand's idea of formal proofs without detours
in a different form.
\end{description}
The last two items are strongly connected 
because the calculus of \heijenoortsrepair\
makes \herbrandindex\herbrand's construction of linear first-order derivations
particularly elegant and noteworthy. 

Therefore, we will present in this paper a novel elaboration 
of the last two 
points.%
\end{sloppypar}%
\section{Introduction}\label
{section Introduction}
\begin{sloppypar}
Before we start, let us briefly mention the complementarity
of the lives of \herbrandname\ \herbrandlifetime\ 
and \heijenoortname\ \heijenoortlifetime:
\herbrand\ had created
a significant part of the world heritage of modern logic
when he died in \herbranddeathyear. \
\heijenoort\ became the leading conservator 
and splendid commentator and translator 
of the
world heritage of modern logic in the 1960s,
and realized it as a field for historical studies.

The scope of work on history
ranges from the gathering and documenting of facts
over interpretation, translation, commentation, 
and the speculation on the gaps in the documented facts,
up to the manipulation of knowledge on history
for politics or for efficient didactics, 
such as in the notes on history in natural-science textbooks 
\cite[\litsectref{XIII}]{kuhn}. \

\heijenoortindex\heijenoort's interest in history
partly resulted
from his involvement into politics and his realization of the contortion 
of knowledge for the manipulation of 
society~\makeaciteoftwo{heijenoort-life}{heijenoort-work}. \
His work on the mathematics on 
\marx\ and \engels\ \nolinebreak\cite{heijenoort-engels}
obviously had the intention to free himself and the
society from this kind of manipulation.


Our main subject in this paper
is how \heijenoortindex\heijenoort\ himself 
has slightly manipulated the historical
facts on \herbrandsfundamentaltheorem,
mostly in his papers~\makeaciteoftwo
{heijenoort-herbrand}{heijenoort-modern-logic},
which he did not properly publish himself,
but also in his published paper \cite{heijenoort-work-herbrand}. \hskip.3em
In these papers,
\heijenoortindex\heijenoort\ changed the inference rules of 
\herbrandindex\herbrand's 
{\it modus ponens}\/-free calculus 
(which we will present in \nolinebreak
 \sectref{section herbrand's Original Rules}), \hskip.2em
probably for didactical reasons. \hskip.1em
We \nolinebreak think that \heijenoortindex\heijenoort\ 
wanted to put more emphasis on 
the 
elegance
of \herbrandsfundamentaltheorem\
and to \repairverb\ \herbrandsfalselemma\ in a more reasonable way
than \bernays, \goedel, and \dreben\nolinebreak\ did.

Regarding \heijenoortindex\heijenoort's merits in publishing, 
translating, 
and commenting \herbrandindex\herbrand's original work,\footnote{%
See \makeaciteoftwo{herbrand-ecrits-logiques}
{herbrand-logical-writings} and \cite{heijenoort-source-book}
for \heijenoort's publishing, translating, and commenting
of \herbrand's work as a logician. 
Moreover,
see
\makeaciteofthree
{heijenoort-oeuvre-herbrand}
{heijenoort-work-herbrand}
{heijenoort-modern-logic} for \heijenoort's secondary 
discussion of \herbrand's work 
as a logician.%
} \hskip.1em
we \nolinebreak have absolutely no \nolinebreak \mbox{reason} 
to blame \heijenoortindex\heijenoort\ for this
slight manipulation. \hskip.2em 
Our motivation here is a different one:
We will discuss the inference system resulting from
\heijenoortindex\heijenoort's changes 
to \herbrandindex\herbrand's {\it modus ponens}\/-free calculus
in detail in \sectfromtoref
{section calculus}
{section The Generalized Rules of Quantification in the Literature}
(and \nolinebreak also \nolinebreak present 
 an up-to-date free-variable version of \nolinebreak it), \hskip.2em
\mbox{because} it helps us to put \herbrandsfundamentaltheorem\
into a new light,
\mbox{emphasizing} its elegance and practical and theoretical relevance.%
\pagebreak

We \nolinebreak think that 
\herbrandindex\herbrand's work in logic still deserves to be
even more widely known than it is today.
We \nolinebreak also think that our partly didactical intentions in this paper
go well together
with \heijenoortindex\heijenoort's original ones 
regarding his \nolinebreak work \nolinebreak on 
\nolinebreak\herbrandindex\herbrand.

The still contentious assessment of \herbrandsfundamentaltheorem\ is also
expressed at the end of 
the following text, written
%
%
%
by \fefermanname\ \fefermanlifetime\ in 1993,\footnotemark\
which also serves as a further introduction into our subject:
\begin{quote}\sloppy
``\herbrandindex\herbrand's fundamental theorem from his 1930 dissertation gives
a kind of reduction, 
in the limit,
of quantificational logic to propositional logic.
This has led to considerable further work in proof theory,
among which the most noteworthy is that of \gentzenname\
(beginning in 1934); \hskip.2em
it has also been used in recent years as the basis for an approach to automated
deduction.
However \herbrandindex\herbrand's own proofs were flawed;
already in 1939 \bernays\ remarked that 
``\herbrandindex\herbrand's proof is hard to follow\closequotecomma
and in 1943, \goedel\ uncovered an essential gap,
though his notes on this were never published.
Counter-examples to two important lemmas in \herbrandindex\herbrand\
were produced by \drebenname\ in collaboration with his students
and colleagues \andrewsnamewithoutmiddleinitial\ and \aanderaaname\ 
in an article in 1963,
in which they outlined how the arguments could be repaired.
A detailed proof of the crucial lemma was later given by \dreben\
and \dentonname\ in 1966.
\par
Despite the known flaws,
\vanheijenoort\ regarded \herbrandindex\herbrand's theorem as one of the deepest results
of logic.
(In \nolinebreak the \nolinebreak minds of most logicians 
that assessment is debatable,
but the importance of the theorem is undeniable.)''
\par\getittotheright{\cite[\p\,383, \litnoteref{11} omitted]{heijenoort-work}}%
\end{quote}
We will refer to that ``crucial lemma'' as 
{\em``\/\herbrandsfalselemmawithinquotedtext\,''}\/ \hskip.1em
and discuss it in 
\nlbsectref{section herbrand's ``False Lemma'' and its Corrections}.%
\footnotetext{%
 Today, 
 however,
 there is better information available on the work of \gentzenname\
 \gentzenlifetime\ in~\makeaciteoftwo
 {menzler-gentzen-german}{menzler-gentzen-english}, and \goedel's
 notes on \herbrandsfalselemma\ have been published in 
 \cite{goldfarb-herbrand-goedel}.}%

After having made all that material most obvious,
we will not be able to avoid the debate on the depth and centrality of
\herbrandsfundamentaltheorem\ in \sectref{section Why is Central}, \hskip.2em
supporting some of the 
well-known assessments of \heijenoortindex\heijenoort, \bernays, and \feferman.

We will conclude our paper with \sectref{section Conclusion}.
\end{sloppypar}
\vfill\pagebreak
\section
{\herbrand's {\it Modus Ponens}\/\hskip.07em-Free Calculus 
\\in the Eyes of \protect\citet{heijenoort-herbrand}}\label
{section calculus}
\begin{sloppypar}
When \heijenoortindex\heijenoort\ refers to 
\herbrandindex\herbrand's calculus in \makeaciteofthree
{heijenoort-oeuvre-herbrand}{heijenoort-work-herbrand}{heijenoort-modern-logic},
\hskip.3em
he actually describes his own calculus 
of \cite{heijenoort-herbrand}, \hskip.2em
which lacks \herbrandindex\herbrand's 
rules of passage and compensates this
with generalized versions of the rules 
of quantification.\footnotemark\end{sloppypar}

Partly as a further introduction to our subject,
we will now quote from \heijenoortindex\heijenoort's unpublished 
draft paper \cite{heijenoort-herbrand},
but 
---~as a few of \heijenoortindex\heijenoort's technical terms are awkward and 
vary from paper to paper~---
we \nolinebreak will replace them in our quotations 
consistently with modern terminology,\newcommand\existentialoid
{\math\gamma}\newcommand\existentialoidquantifier
{\mbox{\math\gamma-quantifier}}\newcommand\Anexistentialoidquantifier
{A \nolinebreak\existentialoidquantifier}\newcommand\anexistentialoidquantifier
{a \existentialoidquantifier}\newcommand
\Ruleofexistentialoidquantificationfirstpart
{Generalized rule of}\newcommand\Ruleofexistentialoidquantificationsecondpart
{\mbox{\math\gamma-quantification}}\newcommand
\Ruleofexistentialoidquantification
{\Ruleofexistentialoidquantificationfirstpart\
 \Ruleofexistentialoidquantificationsecondpart}\newcommand
\ruleofexistentialoidquantification
{generalized rule of \mbox{\math\gamma-quantification}}\newcommand
\Ruleofsimplificationfirstpart
{Generalized rule of}\newcommand\Ruleofsimplificationsecondpart
{simplification}\newcommand\Ruleofsimplification
{\Ruleofsimplificationfirstpart\ \Ruleofsimplificationsecondpart}\newcommand
\ruleofuniversaloidquantification
{generalized rule of 
{\deltaminus\nolinebreak\hskip-.16em\nolinebreak-\nolinebreak 
quanti\-fi\-ca\-tion}}%
\newcommand\Ruleofuniversaloidquantificationfirstpart
{Generalized rule of}\newcommand
\Ruleofuniversaloidquantificationsecondpart
{\deltaminus\nolinebreak\hskip-.16em\nolinebreak-quantification}\newcommand
\Ruleofuniversaloidquantification
{\Ruleofuniversaloidquantificationfirstpart\ 
 \Ruleofuniversaloidquantificationsecondpart}\newcommand\universaloid
{\math\delta}\newcommand\universaloidquantifier
{\mbox{\math\delta-quantifier}}\newcommand\auniversaloidquantifier
{a \universaloidquantifier}\footnotemark\ \hskip.2em 
based on the classification of reductive inference rules by 
\smullyanname\ \smullyanlifetime\ into 
\math\alpha~(sentential+non-branching),
\math\beta~(sentential+\underline branching),
\math\gamma, 
and \nlbmath\delta, \hskip.2em
\cfnlb\ \cite{smullyan}.\begin{quote}%
``In 1929 \herbrandname\ stated a theorem that, 
as \bernays\ \shortcite{bernays-herbrand} writes,
is `the central theorem of predicate logic\closesinglequotefullstop
The theorem has many applications, 
one of them being a new approach to the problem of ascertaining
the validity of a formula.
\par
\herbrand's theorem, stated in a general form, is this: \
Given a formula \nlbmath F of classical quantification theory,
we can effectively generate an infinite sequence of quantifier-free
formulas \nlbmaths{F_1,F_2,\ldots},
such that \nlbmath F is provable in (any standard system of)
quantification theory \ifandonlyif\ there is a \nlbmath k such that
\math{F_k} \nolinebreak is (sententially) valid;
moreover,
\math F \nolinebreak can be retrieved from \nlbmath{F_k} 
through applications of
some rules that have remarkable properties.
There are several ways
of generating an adequate sequence of quantifier-free formulas;
and,
moreover,
many details can vary; \hskip.2em
here we shall consider two kinds of sequences:
expansions and disjunctions.\par
We consider a system \nlbmath Q adequate for classical quantification theory,
without identity, \hskip.1em
in which the sentential connectives are negation, disjunction, 
and \mbox{conjunction.} \hskip.2em
A formula of \nlbmath Q is {\em rectified}\/
\ifandonlyif\ it contains no vacuous quantifier,
no variable has free and bound occurrences in it,
and no two quantifiers bind occurrences of the same variable.
A quantifier in a rectified formula is {\em\existentialoid} \hskip.1em
\ifandonlyif\ it is existential and in the scope of an even number of
negation symbols,
or universal and in the scope
of an odd number of negation symbols; \hskip.2em
otherwise the quantifier is {\em\universaloid}. \hskip.4em
(\Anexistentialoidquantifier\ turns up as existential in a prenex form of
 the formula, 
 and \auniversaloidquantifier\ as universal.) \hskip.3em
A variable in a rectified formula is {\em\existentialoid}\/ \hskip.1em
\ifandonlyif\ it is bound by \anexistentialoidquantifier; \hskip.2em
it is \actualpagebreak
\universaloid\ \hskip.1em
\ifandonlyif\ it is free or bound by 
\auniversaloidquantifier. \hskip.2em
(\herbrandindex\herbrand\ uses `{\frenchfont restreint}' 
 for `\existentialoid' and 
 `{\frenchfont general}'
 for `\universaloid'; \hskip.1em
 in work on \herbrandindex\herbrand's theorem in English
 `restricted' and `general' have been used.)\par
A quantifier is {\em accessible}\/ in a formula \nlbmath F \ifandonlyif\ it is
in the scope of no quantifier of \nlbmaths F. \hskip.3em
A variable is {\em accessible}\/ in a formula \ifandonlyif\ it is either
free or bound by an accessible quantifier.''
\par\getittotheright{\cite[\PP 1 2\@. \ See our \noteref{note rewriting}]
{heijenoort-herbrand}}\pagebreak\end{quote}

\addtocounter{footnote}{-1}\footnotetext{%
 \label{note on heijenoort-work-herbrand}%
 Indeed, \herbrandindex\herbrand\ did not have any {\em generalized}\/ rules
 of quantification; 
 so \herbrandindex\herbrand\ additionally needed 
 his rules of passage to achieve completeness,
 resulting in a calculus of at least four rules.
 Therefore, \heijenoortindex\heijenoort\ definitely refers to his own three-rule 
 calculus when he 
 writes:\begin{quote}
 ``A second result of \herbrandindex\herbrand, connected with the first,
 is that the formula \nlbmath F can be recovered from the formula
 \nlbmath{F_k} by means of three rules of a very definite character.
 These rules are the generalized rule of \mbox{\math\gamma-quantification},
 the \ruleofuniversaloidquantification\
 and the generalized rule of simplification.''
 \getittotheright
 {\cite[\p\,99\@. \ See our \noteref{note rewriting}]
 {heijenoort-work-herbrand}}\end{quote}
 And, one page later, \heijenoortindex\heijenoort\ explicitly declares 
 his own three-rule calculus
 to be \herbrandindex\herbrand's achievement in history:\begin{quote}``%
  The system based on \herbrandindex\herbrand's three rules is, \hskip.1em
  historically, \hskip.1em
  the first example of what we \nolinebreak call \nolinebreak today 
  a cut-free system; \hskip.3em
  it also has the so-called subformula property.''
  \getittotheright
  {\cite[\p\,100]{heijenoort-work-herbrand}}\end{quote}}%
\addtocounter{footnote}{1}\footnotetext{\label{note rewriting}%
 \majorheadroom
 Besides adding a dot after \heijenoort's quantifiers, \hskip.1em
 this means that the following technical terms will be tacitly rewritten when 
 quoting \cite{heijenoort-herbrand,heijenoort-work-herbrand}:\nopagebreak
 \yestop\par\noindent\LINEnomath{\begin{tabular}{@{}l | l | l | l@{}}
  \citet{heijenoort-herbrand} 
 &\citet{heijenoort-work-herbrand} 
 &\citet{heijenoort-modern-logic} 
 &Our wording here
\\\hline
  existentialoid
 &existentialoid
 &restricted
 &\existentialoid
\\existentialoid quantifier
 &
 &
 &\existentialoidquantifier
\\Rule of existentialoid 
 &Generalized rule of
 &Rule of 
 &\Ruleofexistentialoidquantificationfirstpart
\\[-.5ex]\mbox{}\hfill quantification
 &\mbox{}\hfill existentialization
 &\mbox{}\hfill existentialization
 &\mbox{}\hfill\Ruleofexistentialoidquantificationsecondpart
\\Rule of simplification
 &\Ruleofsimplificationfirstpart
 &Rule of
 &\Ruleofsimplificationfirstpart
\\[-.5ex]
 &\mbox{}\hfill\Ruleofsimplificationsecondpart 
 &\mbox{}\hfill simplification
 &\mbox{}\hfill\Ruleofsimplificationsecondpart
\\Rule of universaloid 
 &Generalized rule of 
 &Rule of 
 &\Ruleofuniversaloidquantificationfirstpart
\\[-.5ex]\mbox{}\hfill quantification
 &\mbox{}\hfill universalization
 &\mbox{}\hfill universalization
 &\mbox{}\hfill\Ruleofuniversaloidquantificationsecondpart
\\universaloid 
 &universaloid
 &nonrestricted
 &\universaloid
\\universaloid quantifier
 &
 &
 &\universaloidquantifier
\\\end{tabular}}}%
Note that both \herbrandindex\herbrand\ and \heijenoort\ 
consider equality of formulas only up to
  renaming of bound variables and often implicitly assume that 
  a formula is rectified.

As \heijenoort's 
three-rule inference system 
is difficult to understand
and not complete in the given form, 
let us present each of the rules twice in the following sub-sections: 
first in quotation of \cite{heijenoort-herbrand}, \hskip.2em
second in a hopefully more readable transcription.

Note that 
the three rules can be understood as \mbox{operating} 
on rectified formulas only, \hskip.2em
but this is not necessary for their soundness.
\subsection{\math\gamma-Quantification}\label
{section gamma-Quantification}
\heijenoort\ writes:
\begin{quote}
``
To pass from \maths{F\inparentheses{Q x\stopq H\nottight/H\inpit{x/t}}}, \
where \math{Q x.} \nolinebreak 
is an accessible \existentialoidquantifier\ of \nlbmaths F, \
\math H \nolinebreak is the scope of \nlbmaths{Q x.}, \
and \math t \nolinebreak is a term of the system, \
to \nlbmaths F.\,''
\getittotheright{\cite[\p\,4\@. \ See our \noteref{note rewriting}]
{heijenoort-herbrand}}\end{quote}
Herein ``\math{H\inpit{x/t}\,}''
denotes a capture-avoiding replacement 
of all occurrences of \nlbmath x in \nlbmath H 
with \nlbmath t,
where those bound variables
in \nlbmath H are renamed whose quantifiers would otherwise 
capture free variables in \nlbmath t.

Let \math{A[\ldots]} be the context such that 
\heijenoortindex\heijenoort's \nlbmath F is
\math{A[Q x\stopq H]}. \
Moreover, let us 
write \nolinebreak``\math{H\{x\mapsto t\}\,}'' for the 
result of the replacement 
of all occurrences of \nlbmath x in \nlbmath H 
with \nlbmath t (without renaming). \
Then
we get the following formulation.\pagebreak\begin{quote}
{\bf\Ruleofexistentialoidquantification:} \ 
\LINEmaths{\begin{array}[c]{l}{A[H\{x\mapsto t\}]}
\\\hline
  {A[Q x\stopq H]}
\\\end{array}}{}
where\begin{enumerate}\noitem\item
  \math{Q x.} is an accessible\footnote{%
The restriction of accessibility of the introduced quantifier 
in the \ruleofexistentialoidquantification\ is not necessary for soundness 
(contrary to the \ruleofuniversaloidquantification).
This restriction is introduced, however, 
because already the restricted rule suffices 
for completeness and
for the constructions in the proof of
\herbrandsfundamentaltheorem.
Moreover, the restriction of accessibility guarantees
the equivalence of the generalized rules of quantification 
with their non-generalized versions 
via the rules of passage; \hskip.2em
\cfnlb\ \cite[\p\,6]{heijenoort-tree-herbrand}.}
\existentialoidquantifier\ of 
\nlbmath{A[Q x\stopq H]}, \ and\item the free variables of the term \nlbmath t 
\hskip.1em must\footnote{%
  \majorheadroom
  If a variable \nlbmath z \hskip.1em
  is bound by quantifiers in \nlbmath H
  (as in \examref
   {example Application of the rule of existentialoid quantification}) \hskip.1em
  and occurs free in the term \nlbmath t 
  (in \nolinebreak violation of Side-Condition\,2), \hskip.2em
  an implicit renaming of the bound occurrences of \nlbmath z in 
  \nlbmath H is admitted to enable a backward \mbox{(\ie\,{\em reductive}\/)} 
  application of the generalized rule of 
  \mbox{\math\gamma-quantification}. \hskip.3em
  Note that this is never required for reductive 
  proof search
  if \nolinebreak the formula \nlbmath{A[Q x\stopq H]} is rectified, \hskip.1em
  because then, \hskip.06em 
  in case of violation of Side-Condition\,2, \hskip.1em
  the variable \nlbmath z cannot occur in \nlbmaths{A[\ldots]}, \hskip.2em
  and thus choosing another name for \nlbmath z in \nlbmath t
  can circumvent the renaming of bound occurrences of \nlbmath z in \nlbmath H
  ---~without destroying 
  the property of being a (completed) formal proof later on,
  provided that this other name is chosen instead of \nlbmath z consistently.%
}
not be bound by quantifiers in \nlbmaths H.%
\majorfootroom\end{enumerate}\end{quote}

\yestop
\begin{example}[Application of the \ruleofexistentialoidquantification]\label
{example Application of the rule of existentialoid quantification}%
\mbox{}
\\\noindent If the variable \math z does not occur free in the term \nlbmaths t,
we get the following two inference steps with identical premises
by application of the \ruleofexistentialoidquantification\
at two different positions:\begin{itemize}\noitem\item
\math{\begin{array}[c]{r c r}\inpit{t\tightprec t}
 &\oder
 &\neg\forall z\stopq\inpit{t\tightprec z}
\\\hline
  \inpit{t\tightprec t}
 &\oder
 &\exists x\stopq\neg\forall z\stopq\inpit{x\tightprec z}
\\\end{array}}{} \
via the meta-level substitution
\getittotheright{\maths{
\{\ \ \ \ A[\ldots]\ 
\mapsto\ \inpit{t\tightprec t}
\oder [\ldots]\comma\ \ \
H\ \mapsto\ \neg\forall z\stopq\inpit{x\tightprec z}
\comma\ \ \
Q\ \mapsto\ \exists
\ \ \ \ \}\footroom\headroom};\mbox{~~~~~\,}}\noitem\item
\maths{\begin{array}[c]{r c r}\inpit{t\tightprec t}
 &\oder
 &\neg\forall z\stopq\inpit{t\tightprec z}
\\\hline
  \inpit{t\tightprec t}
 &\oder
 &\neg\forall x\stopq\forall z\stopq\inpit{x\tightprec z}
\\\end{array}}{} \
via the meta-level substitution
\getittotheright{\maths{
\{\ \ \ \ A[\ldots]\ \mapsto\ 
\inpit{t\tightprec t}
\oder\neg[\ldots]
\comma\ \ \
H\ \mapsto\ \forall z\stopq\inpit{x\tightprec z}
\comma\ \ \
Q\ \mapsto\ \forall
\ \ \ \ \}\footroom\headroom}.}\qed\end{itemize}
\end{example}

\halftop\noindent
Today's calculi for proof search have an additional kind of 
free variables, called ``free variables'' 
in \cite{fitting}, \cite{SR--2011--01}, \hskip.2em
or ``(free) \math\gamma-variables'' in \makeaciteoftwo{wirthcardinal}{wirth-jal}.
\hskip.4em
We \nolinebreak will mark such a variable with the upper index \nlbmath\gamma,
such as in 
``\existsvari x{}\nolinebreak\hskip.1em\nolinebreak\closequotefullstopextraspace
By taking the meta-variable \nlbmath t of the generalized rule of 
\mbox{\math\gamma-quantification} to be such a 
free \mbox{\math\gamma-variable}, these {\em free-variable calculi}\/
admit us to delay the choice of the witnessing term \nlbmath t
in a backward application of the \ruleofexistentialoidquantification\ in 
reductive proof search. \hskip.1em
Later, 
when the state of the reductive proof attempt provides sufficient information
which witness is promising, \hskip.1em
the free \mbox{\math\gamma-variable} can be 
globally instantiated in the whole proof tree.

\begin{sloppypar}
Although free-variable calculi can still use the 
generalized rule of \mbox{\math\gamma-quantification} 
in the previous formulation,
we will refer several times to the following sub-rule:
\end{sloppypar}

\begin{quote}
{\bf Generalized rule of restricted \math\gamma-quantification:} \ 
\LINEmaths{\begin{array}[c]{l}{A[H\{x\mapsto\existsvari x{}\}]}
\\\hline
  {A[Q x\stopq H]}
\\\end{array}}{}
where \math{Q x.} is an accessible \existentialoidquantifier\ of 
\nlbmath{A[Q x\stopq H]}.\\\mbox{}\majorfootroom\end{quote}\pagebreak
\subsection{\math\delta-Quantification}\label
{section delta-Quantification}
\heijenoortindex\heijenoort\ writes:\notop\halftop
\begin{quote}\sloppy
``To pass from \bigmaths{F\inparentheses{Q y\stopq H\nottight/H}}, \
where \math{Q y.} \nolinebreak 
is an accessible \universaloidquantifier\ of \nlbmaths F, \
\math H \nolinebreak is the scope of \nlbmaths{Q y.}, \
and \math y does not occur free in \nlbmaths F, \
to \nlbmaths F.\,''
\getittotheright{\cite[\p\,4\@. \ See our \noteref{note rewriting}]
{heijenoort-herbrand}}\end{quote}
Let \math{A[\ldots]} be the context such that 
\heijenoortindex\heijenoort's \math F is
\nlbmaths{A[Q y\stopq H]}. \
Then
we \nolinebreak get \nolinebreak the following formulation.\begin{quote}
{\bf\Ruleofuniversaloidquantification:} \LINEmaths{\begin{array}[c]{l}{A[H]}
  \\\hline
    {A[Q y\stopq H]}
  \\\end{array}}{} \
where\begin{enumerate}\item 
\math{Q y.} is an accessible \universaloidquantifier\
of \nlbmaths{A[Q y\stopq H]}, \ and\item
the variable \nlbmath y must\footnotemark\
not occur 
free 
in the context \nlbmath{A[\ldots]}.\majorfootroom\end{enumerate}\end{quote}
\footnotetext{%
  If \math y occurs in the context \nlbmath{A[\ldots]}, an implicit
  renaming of the bound occurrences of \nlbmath y in 
  \bigmaths{Q y\stopq H}{} is admitted to enable backward 
  (\ie\ {\em reductive}\/)
  application of the inference rule. \
  Note that this is never required for reductive 
  proof search
  if the formula \nlbmath{A[Q y\stopq H]} is rectified,
  because then \nlbmath y cannot occur in \nlbmaths{A[\ldots]}.
}%
For soundness reasons, free-variable calculi 
(\ie\ calculi with free \math\gamma-variables) \hskip.1em
need a different rule
of \mbox{\math\delta-quantification}. \hskip.3em
We \nolinebreak take here a generalized version 
of the so-called \mbox{\deltaplusplus-rule} \mbox{\cite{deltaplusplus}}. \
This \mbox{\deltaplusplus-rule} is a {\em liberalized}\/ \math\delta-rule
in the sense that, compared to the simple \math\delta-rule
(also called \deltaminus-rule), \hskip.1em
it admits additional proofs that are more easily found
by man and machine. (See
\sectref{section Formal Proof in Our Free-Variable Calculus}, 
and also
\cite[\litsectref{2.1.5}]{wirthcardinal},
\cite{nonpermut}.)

\begin{quote}
{\bf Generalized rule of \deltaplusplus-quantification:} \\
 \LINEmaths{\begin{array}[c]{l}
 {A[H\{y\mapsto\inpit{Q y\stopq H,\existsvari x 1,\ldots,\existsvari x m}
 ^\delta(\existsvari x 1,\ldots,\existsvari x m)\}]}
  \\\hline
    {A[Q y\stopq H]}
  \\\end{array}}{} \
where\begin{enumerate}
\item
\math{Q y.} is an accessible \universaloidquantifier\
of \nlbmaths{A[Q y\stopq H]},
\item\sloppy\math{\existsvari x 1,\ldots,\existsvari x m} are exactly the 
\math m~free \math\gamma-variables occurring in \nlbmaths{Q y\stopq H}, \
and
\item
the \math\delta\ in 
\math{\inpit{Q y\stopq H,\existsvari x 1,\ldots,\existsvari x m}^\delta} 
denotes the application of the function \nlbmath\delta\ \hskip.1em to
\nlbmaths{\inpit{Q y\stopq H,\existsvari x 1,\ldots,\existsvari x m}}, 
\hskip.2em
resulting in a function symbol of arity \nlbmath m \hskip.1em
that is not part of the original signature, \hskip.2em
such that, \hskip.2em
in \nolinebreak case of 
\bigmaths{\inpit{B,\existsvari x 1,\ldots,\existsvari x m}^\delta=
\inpit{C,\existsvari z 1,\ldots,\existsvari z{m'}}^\delta}, 
either \inpit{C,\existsvari z 1,\ldots,\existsvari z {m'}} 
\nolinebreak is 
\nlbmaths{\inpit{B\existsvari x 1,\ldots,\existsvari x m}}, \ 
or \math C results from \nlbmath B \hskip.1em
by renaming of bound variables
and by a bijective renaming of the free \math\gamma-variables 
\math{\existsvari x 1,\ldots,\existsvari x m} via 
\bigmaths{\{
\existsvari x 1\tight\mapsto\existsvari z 1\comma\ldots\comma
\existsvari x m\tight\mapsto\existsvari z{m'}
\}}.\end{enumerate}\end{quote}
In case of \bigmaths{m\tightequal 0},
\ie\ in the absence of free \math\gamma-variables,
there is no essential difference between the generalized rules
of \deltaminus- and \math{\deltaplusplus\!}-quantification,
except that the \deltaplusplus-rule may violate Side-Condition\,2 of the
\ruleofuniversaloidquantification\
(when reductively \mbox{applied} a second time to the same formula).
\subsection{Simplification}\label{section generalized rule of simplification}
\heijenoortindex\heijenoort\ writes:
\begin{quote}
``To pass from \bigmaths{F\inparentheses{H\nottight/\inpit{H\oder H'}}}, \
where \math H \nolinebreak is a subformula of \nlbmath F and
\math{H'} \nolinebreak is a variant of \nlbmaths H, \
to \nlbmaths F.\,''
\getittotheright{\cite[\p\,4]{heijenoort-herbrand}}\end{quote}
To satisfy the needs of the proof construction in 
\herbrandsfundamentaltheorem,
we have to correct \heijenoortindex\heijenoort's version a bit:
\mbox{On the} one hand,
\mbox{we remove} \heijenoortindex\heijenoort's admission of  
\bigmaths{H\oder H'}{} at {\em negative positions}\/
(\ie\nolinebreak\ in \nolinebreak 
 the scope of an odd number of negation symbols),
and, on the other hand,
we admit also \bigmaths{H\und H'}{} at negative positions.\footnote{%
 See the formal proof in \sectref{section example},
 where 
 we need the operator ``\tightund'' twice at negative positions.%
}

\yestop\noindent Let \math{A[\ldots]} be the context such that 
\heijenoortindex\heijenoort's \math F is \math{A[H]}. \
Then we get the following formulation.\begin{quote}
{\bf\Ruleofsimplification:} \bigmaths{\begin{array}[c]{l}{A[H\circ H']}
\\\hline
  {A[H]}
\\\end{array}} \ where\begin{enumerate}\item
``\math{\circ}'' stands for ``\tightoder'' if \math{[\ldots]} \nolinebreak
denotes a positive position in \nlbmath{A[\ldots]}, \ 
and for ``\tightund'' if this position is negative,\footnote{%
 \majorheadroom
 In the terminology of \smullyan's classification,
 this means that \math\circ\ is always an \mbox{\math\alpha-operator}.%
} \ and\item\math{H'} is a variant of the sub-formula \nlbmath H \
(\ie, \math{H'} is \math H or 
 can be obtained from \nlbmath H by the renaming of variables
bound in \nlbmath H).\\\mbox{}\footroom\end{enumerate}\end{quote}
Moreover, the 
{\em generalized rule of \math\gamma-simplification} 
  is the sub-rule for the case that
  \math H \nolinebreak is of the form \nlbmath{Q y\stopq C}
  and \math{Q y.} is \anexistentialoidquantifier\
  of \nlbmaths{A[Q y\stopq C]}.

\begin{sloppypar}\yestop\noindent
As the only essential function of the generalized rule of simplification
is to increase \mbox{\math\gamma-multiplicity,} \hskip.1em
the generalized rule of \math\gamma-simplification is sufficient for 
completeness and for the proof of \herbrandsfundamentaltheorem.%
\end{sloppypar}
\vfill\pagebreak
\section{\protect\PropertyC}\label
{section property C}

\yestop\subsection{Prerequisites for \protect\PropertyC}\label
{section Prerequisites for Property C}%
\yestop\halftop\begin{definition}[Height of a Term, Champ Fini 
\nlbmath{\termsofdepth n F}]\label
{definition champs finis}\mbox{}
\par\noindent
We \nolinebreak use \nlbmath{\CARD t}
to denote the {\em height}\/ of a term \nlbmath t,
which is given by
\par\noindent\LINEmaths{\CARD{\anonymousfpp{t_1}{t_{m}}}
  \nottight{\nottight{\nottight=}}1+\max\{0,\CARD{t_1},\ldots,\CARD{t_m}\}}.
\par\noindent For a positive natural number \nlbmath n
and a formula \nlbmath F,
as a finite substitute for a typically infinite, 
full term universe,
\herbrandindex\herbrand\ uses what he calls a 
{\em{\frenchfont\em champ fini} of order \nlbmaths{n},} \hskip.1em
which we \nolinebreak will denote with 
\nlbmath{\termsofdepth n F}.  \ 
The terms of \nlbmath{\termsofdepth n F} 
are constructed from the symbols that occur free
in \nlbmath F: \
the function symbols, the constant symbols
(which we will tacitly subsume 
 under the function symbols in what follows), \hskip.1em
and the free variable symbols
(which can be seen as constant symbols here). \hskip.3em
Such a 
{\frenchfont champ fini}\/ differs from a full term universe in
containing only the terms \nlbmath t with \bigmaths{\CARD t\prec n}{\,.} \ 
\par\noindent 
So we have \bigmaths{\termsofdepth 1 F\tightequal\emptyset}. \hskip.2em
\par\noindent To guarantee \bigmaths
{\termsofdepth n F\tightnotequal\emptyset}{} \hskip.1em
for \bigmaths{n\succ 1},
in case that neither constants nor free variable symbols occur in \nlbmaths F,
\hskip.1em we will assume that a fresh constant symbol 
\nolinebreak ``\nlbmath\bullet\nolinebreak\hskip.02em\nolinebreak''
(which does not occur elsewhere)
\hskip .2em is \nolinebreak
included in the term construction in addition to the free symbols of \nlbmath F.
\getittotheright\qed\end{definition}

\begin{sloppypar}
\yestop\yestop\halftop\noindent\herbrand's definition of an expansion 
follows the traditional
idea that ---~for a {\em finite}\/ domain~--- universal (existential) 
quantification can be seen as a finite conjunction \mbox{(disjunction)} over the 
elements of the domain. \hskip.2em
To 
reduce the size of formal 
proofs according to \herbrandsfundamentaltheorem,
we will consider sub-expansions in addition.\end{sloppypar}

\yestop\halftop\begin{definition}[\opt{Sub-} Expansion]\label
{definition sub-expansion}\sloppy\mbox{}\\Let \nlbmath{\mathcal T} be a 
finite set of terms. \
To simplify substitution,
let \nlbmath A be a rectified formula whose bound variables do not occur in 
\nlbmaths{\mathcal T}. \ 
A \nolinebreak 
formula \nlbmath B is a {\em sub-expansion of \nlbmath A \wrt\ \nlbmath
{\mathcal T}} \hskip.3em
if \hskip.2em\mbox{\maths A, \math B, \math{\mathcal T}} 
satisfy the recursive definition of the following table,
where ``s-e'' abbreviates ``is a sub-expansion'':
\par\noindent
\getittotheright{\begin{tabular}{l|l|l}
  form of \nlbmath A
 &form of \nlbmath B
 &required properties
\\\hline
  quantifier free
 &\math A
 &---
\\\math{\neg A'}
 &\math{\neg B'}
 &\math{B'} s-e of \math{A'} \wrt\ \nlbmath{\mathcal T}
\\\math{A_1'\oder A_2'}
 &\math{B_1'\oder B_2'}
 &\math{B_i'} s-e of \math{A_i'} \wrt\ \nlbmaths{\mathcal T\!}, \hskip.2em 
  for \math{i\tightin\{1,2\}}
\\\math{A_1'\und A_2'}
 &\math{B_1'\und B_2'}
 &\math{B_i'} s-e of \math{A_i'} \wrt\ \nlbmaths{\mathcal T\!}, \hskip.2em 
  for \math{i\tightin\{1,2\}}
\\\math{\exists x\stopq A'}
 &\math{\bigvee_{t\in\mathcal T'}B'\{x\tight\mapsto t\}}
 &\math{B'} s-e of \math{A'} \wrt\ \nlbmaths{\mathcal T\!}, \hskip.2em 
  and \nlbmath{\emptyset\tightnotequal\mathcal T'\tightsubseteq\mathcal T}
\\\math{\forall x\stopq A'}
 &\math{\bigwedge_{t\in\mathcal T'}B'\{x\tight\mapsto t\}}
 &\math{B'} s-e of \math{A'} \wrt\ \nlbmaths{\mathcal T\!}, \hskip.2em 
  and \nlbmath{\emptyset\tightnotequal\mathcal T'\tightsubseteq\mathcal T}
\\\end{tabular}}
\par\halftop\halftop\noindent
If we restrict \math{\mathcal T'} in this table to be equal to 
\nlbmaths{\mathcal T\!}, \hskip.2em
then we get an {\em expansion \nlbmath B of \nlbmath A}
(instead of a sub-expansion \nlbmath B of \nlbmath A) 
\wrt\ \nlbmath{\mathcal T\!}, \hskip.2em
and write \nlbmath{A^{\mathcal T}} for \nlbmaths B.
\getittotheright\qed\end{definition}
\vfill\pagebreak

\begin{definition}[\skolemizedForm s]\label
{definition skolemized forms}\mbox{}\par\noindent Let \math B be a rectified
formula.
\par\noindent The {\em\outerskolemizedform\ of a formula \nlbmath B}\/
(also called: ``outer \skolemizedform'') \hskip.1em
results from \nlbmath B \hskip.1em 
by removing every \universaloidquantifier\ and replacing its bound
variable \nlbmath y with \nlbmaths{\app{\forallvari y{}}{x_1,\ldots,x_m}}, \ 
where \math{\forallvari y{}} \nolinebreak is a fresh (``\skolem'') symbol
and \maths{x_1,\ldots,x_m}, \hskip.1em
in this order, \hskip.1em
are exactly the \mbox{\math m~variables} 
of the \mbox{\math m~\math\gamma-quantifiers} in
whose scope the \mbox{\math\delta-quantifier} occurs.

\par\noindent The {\em\deltaplusplusskolemizedform\ of \nlbmath B} 
results from \nlbmath B \hskip.1em
by repeating the following procedure
until all \mbox{\math\delta-quantifiers} have been removed: \ 
Let \math{Q y.} be 
an outermost \universaloidquantifier\ 
\mbox{in \math{A[Q y\stopq H]}} \math= \maths B. \hskip.7em
Replace \nlbmath B \hskip.1em
with \bigmathnlb
{A[H\{y\mapsto\inpit{Q y\stopq H, x_1,\ldots,x_m}^\delta(x_1,\ldots,x_m)\}]}, 
\hskip.1em
where
\maths{x_1,\ldots,x_m}{}
are exactly the \math m~variables bound by the 
\mbox{\math m~\math\gamma-quantifiers} in whose scope the 
\mbox{\math\delta-quantifier \math{Q x.}}
occurs and which {\em do actually occur in the scope \nlbmath H \hskip.1em
of the \mbox{\math\delta-quantifier}}. \hskip.4em 
Moreover, \hskip.2em
the \math\delta\ in \nolinebreak
``\math{\inpit{Q y\stopq H,x_1,\ldots,x_m}^\delta}\nolinebreak
\hskip.15em\nolinebreak'' 
denotes the application of the function \nlbmath\delta\ 
to \nlbmaths{(Q y\stopq H,x_1,\ldots,x_m)}. \hskip.6em 
The result of this application 
must be a function symbol of arity \nlbmath m \hskip.1em
that is not part of the original signature, \hskip.2em
such that, \hskip.2em
in \nolinebreak case of 
\bigmaths{\inpit{B,x_1,\ldots,x_m}^\delta=\inpit{C,x'_1,\ldots,x'_{m'}}^\delta}, 
either \inpit{C,x'_1,\ldots,x'_{m'}} \nolinebreak is 
\inpit{B,x_1,\ldots,x_m}, \hskip.3em 
or \math C
results from \nlbmath B \hskip.1em
by renaming of bound variables
and by a bijective renaming of the free variables \nlbmath{x_1,\ldots,x_{m}} 
via \bigmaths
{\{x_1\tight\mapsto x'_1\comma\ldots\comma x_m\tight\mapsto x'_{m'}\}}. 
\getittotheright\qed\end{definition}

\yestop\begin{remark}[\skolemizedForm s]\label
{remark skolemized forms}\mbox{}\par\noindent
\herbrandindex\herbrand\ has no name for the \outerskolemizedform\
 and he does not use the \deltaplusplusskolemizedform,
 which is the current standard in two-valued \firstorder\ logic
because it \nolinebreak is closely related to the liberalized \math\delta-rules.
\par\noindent
We cannot get along without the \deltaplusplusskolemizedform\ here
because there is no reasonable version of 
the generalized rule of \mbox{\math{\deltaplusplus\!}-quantification}
that is compatible with the \outerskolemizedform.
\end{remark}

\yestop\yestop\yestop\begin{definition}[\SententialTautology]\mbox{}\\
A \firstorder\ formula is a {\em\sententialtautology}\/ if it
    is quantifier-free and truth-functionally valid, provided its
    atomic subformulas are read as atomic sentential variables.%
\qed\end{definition}
  
\vfill\pagebreak

\subsection{\protect\PropertyC\ and \protect\PropertyCstar}

\yestop\halftop\begin{definition}[\PropertyC, \PropertyCstar]\label
{definition properties C and C star}\mbox{}
\par\noindent Let \math A be a rectified \firstorder\ formula. \ 
Let \math n be a positive natural number.
\par\noindent \math A \nolinebreak
{\em has \propertyC\ of order \nlbmath n}\/ \udiff\ the 
expansion \nlbmath{F^{\,\termsofdepth n F}} is a
\sententialtautology,
where \math F is
the \outerskolemizedform\ of \nlbmaths A.
\par\noindent\math A \nolinebreak
{\em has \propertyCstar\ of order \nlbmath n}\/ \udiff\ the 
expansion \nlbmath{F^{\,\termsofdepth n F}} is a
\sententialtautology,
where \math F is
the \deltaplusplusskolemizedform\ of \nlbmaths A.
\getittotheright\qed\end{definition}

\yestop\noindent
Note that \math{F^{\,\termsofdepth 1 F}} is defined
\ifandonlyif\ \math F does not contain any quantifier,
because we have \bigmaths{\termsofdepth 1 F\tightequal\emptyset},
and because an empty conjunction or disjunction for the expansion of
a quantifier is not admitted in \defiref{definition sub-expansion}. \hskip.3em
Thus, we get the following corollary, 
which ---~to avoid the usage of an undefined term~---
may also serve as a cleaner definition for \propertyC\ and \propertyCstar\
of order \nlbmaths 1.
\begin{corollary}\label{corollary sentential tautology}%
\par\noindent The following three are logically equivalent 
for a rectified formula \nlbmaths A:
\begin{enumerate}\noitem\item
\math A has \propertyC\ of order \nlbmaths 1.\item
\math A has \propertyCstar\ of order \nlbmaths 1.\item
The 
(no matter whether \deltaminus\nolinebreak\hskip-.16em\nolinebreak-
 or \math{\deltaplusplus\!}-)\/ 
\skolemizedform\ of \nlbmath A\\is a \sententialtautology.
\getittotheright\qed\end{enumerate}\end{corollary}

\yestop\yestop\noindent
Because all quantifiers in a \skolemizedform\ \nlbmath F 
are \math\gamma-quantifiers,
existential quantifiers occur only at positive positions 
(and are expanded with disjunctions) 
and universal quantifiers occur only at negative positions 
(and are expanded with conjunctions). \
Thus, 
we \nolinebreak get the following corollaries:
\begin{corollary}\label{corollary sub-expansion}%
\\If any sub-expansion of a \skolemizedform\ \nlbmath F 
\wrt\ \nlbmath{\mathcal T}
is a \sententialtautology, \hskip.1em
\\then \math{F^{\,\mathcal T}\!} is a \sententialtautology\ \aswell.
\getittotheright\qed\end{corollary}
\begin{corollary}\label{corollary sub-expansion two}\sloppy%
\\Let \math A be a rectified \firstorder\ formula. \ 
Let \math n be a positive natural \mbox{number.} \
If any sub-expansion of the \outerskolemizedformwithoutform\
(or \nolinebreak else: \math{\deltaplusplus\!}-\skolem ized) \hskip.1em
form \nlbmath F of \nlbmath A
\mbox{\wrt\ \termsofdepth n F} \mbox{is a \sententialtautology,} \hskip.2em
then \math A has \propertyC\ 
(or \nolinebreak else: \propertyCstar)
\mbox{of order \nlbmath n.}
\getittotheright\qed\end{corollary}
\vfill\pagebreak

\halftop\halftop\noindent
By using the idea of 
\cororef{corollary sub-expansion} for the forward direction and
by replacing all terms whose top function symbols are new
with the same old constant or free variable 
for the backward direction, 
we \nolinebreak get:

\begin{corollary}\label{corollary extend signature}%
\par\noindent Let \math A be a rectified \firstorder\ formula. \ 
Let \math n be a positive natural \mbox{number.}
\\Let\/ \math F be the \outerskolemizedformwithoutform\
(or \nolinebreak else: \math{\deltaplusplus\!}-\skolem ized) 
form of \nlbmaths A.
\\Let \math{\mathcal T'\!} be formed just as\/ \termsofdepth n F,
but over an extended set of 
function 
symbols.
\par\noindent Now the following two are 
logically equivalent:\begin{enumerate}\noitem\item
\ \math A has \propertyC\ 
(or \nolinebreak else: \propertyCstar) of order \nlbmaths n.\item
\ \math{F^{\,\mathcal T'}\!} is a \sententialtautology.
\getittotheright\qed\end{enumerate}\end{corollary}

\yestop\yestop\noindent
As the replacement of 
\deltaminus\nolinebreak\hskip-.16em\nolinebreak-\skolem\ terms 
with \deltaplusplus-\skolem\ terms,
dropping the arguments at the missing argument positions,
cannot increase the range of values of a sentential expression under 
all valuations 
(of its sentential variables named by atomic formulas), \hskip.1em
we \nolinebreak get:

\begin{corollary}\label{from C to C star}\mbox{}
Let\/ \math A be a rectified formula. \hskip.3em
Let\/ \math n be a positive natural number. \hskip.4em
\\If\/ \math A \hskip.1em has \propertyC\ of order\/ \nlbmaths n, \hskip.2em 
then\/ \math A \hskip.1em has \propertyCstar\ of order\/ \nlbmath n \hskip.1em
\aswell. \getittotheright\qed\end{corollary}

\halftop\halftop\halftop\noindent 
Note that the converse of \cororef{from C to C star} does not hold
in general:
\begin{example}\label{example Property C one}\sloppy\mbox{}\\
Let us take \nlbmath A to be the rectified formula \bigmaths\formulaexampleC. 
\\Its \deltaplusplusskolemizedform\ is
\bigmaths
{\exists a\stopq\inparentheses
{
{\neg
{\app p{\forallvari b{}}}}\oder{\app p a}}}. 
Now, expansion \wrt\ \nlbmath{\{\forallvari b{}\}} 
results in a \sententialtautology. \hskip.3em
Thus, \hskip.1em
\math A \nolinebreak has \propertyCstar\ of order \nlbmaths 2. \hskip.7em
\\Its \outerskolemizedform\ is \bigmaths\formulaexampleCouterskolemizedform.
Here expansion \wrt\ \nlbmath{\{\bullet\}} does not result in 
a \sententialtautology. \hskip.3em
Thus, \math A \nolinebreak does not have \propertyC\ of order \nlbmaths 2,
\hskip.2em
but only of order \nlbmath 3.\getittotheright\qed\end{example}
\vfill\cleardoublepage
\section{From \protect\propertyC\ to a Linear Derivation}\label
{section example}
In this section,
instead of proving our theorems,
we will exemplify them by a generalizable example.

In this example, 
we add some sugar to our formula syntax. \
%
Let us write
\\[+.9ex]\LINEnomath{``\math{A\implies B}\hskip.1em''}
\\[-.5ex]instead of
\\[-.5ex]\LINEnomath{``\hskip.09em\math{\neg A\oder B}\hskip.1em'',}
\\[+.9ex]and let ``\tightimplies'' have lower operator precedence than 
``\tightund'' and ``\tightoder''.

The following formula
says that if we have 
an upper bound of every two elements as well as transitivity, 
then we also have an upper bound of every three elements.
\mathcommand\lastlineofherbrandformula
{\exists\boundvari z{}\stopq\inpit{
\forallvari u{}\tightprec\boundvari z{}\und
\forallvari v{}\tightprec\boundvari z{}\und
\forallvari w{}\tightprec\boundvari z{}}}
\par\halftop\halftop\noindent\LINEmaths{\noparenthesesoplist{
\forall\boundvari x{}\stopq
\forall\boundvari y{}\stopq
\exists\boundvari m{}\stopq
\inpit{\boundvari x{}\tightprec\boundvari m{}\und 
\boundvari y{}\tightprec\boundvari m{}}
\oplistund\mediumheadroom
\forall\boundvari c{}\stopq
\forall\boundvari b{}\stopq
\forall\boundvari a{}\stopq
\inpit{
\boundvari a{}\tightprec\boundvari b{}
\und
\boundvari b{}\tightprec\boundvari c{}
\implies
\boundvari a{}\tightprec\boundvari c{}}
\oplistimplies\mediumheadroom
\forall\boundvari u{}\stopq
\forall\boundvari v{}\stopq
\forall\boundvari w{}\stopq
\exists\boundvari z{}\stopq
\inpit{
\boundvari u{}\tightprec\boundvari z{}\und
\boundvari v{}\tightprec\boundvari z{}\und
\boundvari w{}\tightprec\boundvari z{}}}}{}{\Large\math{(A)}}
\par\halftop\halftop\noindent
The (\deltaminus\nolinebreak\hskip-.16em\nolinebreak- \aswellas\ 
\math{\deltaplusplus\!}-) 
\skolemizedform\ of \nlbmath A is
\par\halftop\halftop\noindent\LINEmaths{\noparenthesesoplist{
\forall\boundvari x{}\stopq
\forall\boundvari y{}\stopq
\inpit{
\boundvari x{}\tightprec\app{\forallvari m{}}{\boundvari x{},\boundvari y{}}\und 
\boundvari y{}\tightprec\app{\forallvari m{}}{\boundvari x{},\boundvari y{}}}
\oplistund\mediumheadroom
\forall\boundvari c{}\stopq
\forall\boundvari b{}\stopq
\forall\boundvari a{}\stopq
\inpit{
\boundvari a{}\tightprec\boundvari b{}
\und
\boundvari b{}\tightprec\boundvari c{}
\implies
\boundvari a{}\tightprec\boundvari c{}}
\oplistimplies\mediumheadroom
\lastlineofherbrandformula}}{}{\Large\math{(F)}}\par\noindent
\subsection{Informal Proof}
Let us look for an informal proof of the last line of \nlbmath F,
assuming the first two lines to be given as lemmas. \hskip.2em
From semantical considerations, 
it is obvious 
that a solution for \nlbmath{\boundvari z{}} is given by the substitution
\\[-1.7ex]\LINEmaths{\begin{array}{l l l}\sigma_3
 &\ :=\ 
 &\{\boundvari z{}\mapsto\termzwei\}.
\\\end{array}}{}
\par\noindent We can prove this as follows:
First we apply the first line 
twice instantiated by the two substitutions
\\[-1.7ex]\LINEmaths{\begin{array}{l l l}\sigma_{1,1}
 &\ :=\ 
 &\{\boundvari x{}\mapsto\boxv
\comma
\boundvari y{}\mapsto\boxw\},
\\\sigma_{1,2}
 &\ :=\
 &\{\boundvari x{}\mapsto\boxu
  \comma
  \boundvari y{}\mapsto\boxeins
  \}
\\\end{array}}{}
\\\noindent to obtain all of
\par\indent\bigmaths{\forallvari v{}\tightprec\termeins
  \comma
  \forallvari w{}\tightprec\termeins
  \comma}{}
\\\indent\bigmaths{\forallvari u{}\tightprec\termzwei\comma
  \termeins\tightprec\termzwei}.
\par\noindent Now we can apply the second line twice
instantiated by the two substitutions
\par\noindent\LINEmaths{\left.\begin{array}{l l l}
   \sigma_{2,1}
 &\ :=\
 &\sigma_2\cup\{\boundvari a{}\mapsto\forallvari v{}\},
\\\sigma_{2,2}
 &\ :=\ 
 &\sigma_2\cup\{\boundvari a{}\mapsto\forallvari w{}\}
  ~~
\\\end{array}\right\}
\mbox{~~for~~~}
\sigma_2\ :=\ 
\{\boundvari c{}\mapsto\termzwei\comma\boundvari b{}\mapsto\termeins\}
}{}
\par\noindent to obtain also
\par\indent\bigmaths{\forallvari v{}\tightprec\termzwei\comma
  \forallvari w{}\tightprec\termzwei},
\par\noindent which completes our informal proof the last line of 
formula \nlbmaths F.

As the maximum height of the instantiations in this informal proof
is \nlbmaths 3, \hskip.2em
it is obvious that the formula \nlbmath A has \propertyC\ of order \nlbmaths 4.%
\pagebreak
\subsection
{Human-Oriented Formal Proof Construction in \heijenoort's Version of \herbrand's {\it Modus Ponens}\/\hskip.07em-Free Calculus}\label
{section formal proof}%
\yestop\yestop\noindent
Let us use this informal proof to construct a formal 
proof of the formula \nlbmath A
in \heijenoortindex\heijenoort's version of 
\herbrandindex\herbrand's {\it modus ponens}\/\hskip.07em-free calculus,
which we \nolinebreak presented in \sectref{section calculus}. \
More precisely,
our goal is to 
reduce the formula \nlbmath A
to a sub-expansion of \nlbmath F \wrt\ \nlbmaths{\termsofdepth n F}{}
(for \maths{n\tightequal 4})
that is a \sententialtautology; \hskip.2em
our \nolinebreak means will be the renaming of bound variables and
the generalized rules of 
\mbox{\math\gamma-simplification} and 
\mbox{\deltaminus- and} 
\mbox{\math\gamma-quantification}.

In any of the following reductive backward steps, \hskip.1em
the parts of a formula that will be removed 
will be shown in {\color{red}red} and those parts that are crucially involved
(but not removed) are shown in {\color{orange}orange}. \hskip.3em
The parts changed \wrt\ the {\em previous}\/ formula will be shown
in {\color{green}green} if they are crucial for triggering the previous step, 
and in {\color{blue}blue} otherwise.

\yestop\yestop
\subsubsection{Phase\,1 of the Reduction: Rename Bound \math\delta-Variables}
First we use an amazing trick invented by \herbrandindex\herbrand:
We use \skolem\ terms as names for 
bound \math\delta-variables and for the free variables
that result from backward application of the 
\ruleofuniversaloidquantification.

By this trick, \herbrand\ escapes the problems of giving 
\skolem\ functions a semantics and all 
involved problems, 
such as the need for a principle of choice.
These problems are present in \cite{loewenheim-1915}
and were overcome in \makeaciteoftwo{skolem-1920}{skolem-1923b}
only by using {\em\skolemnormalform}\/ instead of \skolemizedform.
\skolemizedform,
however, 
is more intuitive and more
useful in automated theorem proving than \skolemnormalform,
and \skolem\ returned to \skolemizedform\ in~\cite{skolem-1928}.

According to \herbrand's trick, 
every term that has a \skolem\ function as top symbol
is seen as a variable, such that 
\app{\forallvari m{}}{\forallvari v{},\forallvari w{}}
is {\em not}\/ a subterm of 
\app
{\forallvari m{}}
{\forallvari u{},{\app{\forallvari m{}}{\forallvari v{},\forallvari w{}}}},
because they are both just variables.
This interpretation does not affect the formula \nlbmath A,
simply because \math A \nolinebreak 
does not contain any \skolem\ function symbols.

Thus, we start by renaming each bound \math\delta-variable
in \nlbmath A
to its \skolem\ term in the \outerskolemizedform\ \nlbmaths F{} of \nlbmaths A,
and obtain:
\par\halftop\halftop\noindent\LINEmaths{\noparenthesesoplist{
{\color{orange}\forall\boundvari x{}\stopq
\forall\boundvari y{}\stopq
\exists\app{\forallvari m{}}{\boundvari x{},\boundvari y{}}\stopq
\inpit{
\boundvari x{}\tightprec\app
{\forallvari m{}}{\boundvari x{},\boundvari y{}}\und 
\boundvari y{}\tightprec\app
{\forallvari m{}}{\boundvari x{},\boundvari y{}}}}
\oplistund\headroom
\forall\boundvari c{}\stopq
\forall\boundvari b{}\stopq
{\color{orange}
\forall\boundvari a{}\stopq
\inpit{
\boundvari a{}\tightprec\boundvari b{}
\und
\boundvari b{}\tightprec\boundvari c{}
\implies
\boundvari a{}\tightprec\boundvari c{}}}
\oplistimplies\headroom
\forall\forallvari u{}\stopq
\forall\forallvari v{}\stopq
\forall\forallvari w{}\stopq
\lastlineofherbrandformula}}{}{\Large\math{(S)}}
\par\halftop\halftop\noindent
Here, \hskip.05em
in the first line of \nlbmaths S, \hskip.1em 
we must not read the occurrences of the  
variables \boundvari x{} and \nlbmath{\boundvari y{}} as 
arguments of \nlbmath{\forallvari m{}} as variables bound by the outer
universal quantifier, \hskip.1em 
but as a part of the name \nolinebreak``\math{\app
{\forallvari m{}}{\boundvari x{},\boundvari y{}}}''
of the variable bound by the existential quantifier.

\vfill\pagebreak

\halftop
\subsubsection
{Phase\,2 of the Reduction: Generalized Rule of \math\gamma-Simplification}
We know from our informal proof that we 
need the first and the second line twice.
Thus,
we now double them at \math{\forall x.} \hskip.1em
and \nlbmaths{\forall a.}, \hskip.2em
respectively,
by two backward applications of the generalized rule
of \math\gamma-simplification, and obtain:
\par\halftop\noindent\LINEmaths{\noparenthesesoplist{
{\color{blue}\forall\boundvari x{}\stopq
\forall\boundvari y{}\stopq
\exists\app{\forallvari m{}}{\boundvari x{},\boundvari y{}}\stopq
\inpit{
\boundvari x{}\tightprec\app
{\forallvari m{}}{\boundvari x{},\boundvari y{}}\und 
\boundvari y{}\tightprec\app
{\forallvari m{}}{\boundvari x{},\boundvari y{}}}}
\oplistund
{\color{green}\forall\boundvari x 1\stopq
\forall\boundvari y 1\stopq
\exists\app{\forallvari m{}}{\boundvari x 1,\boundvari y 1}\stopq
\inpit{
\boundvari x 1\tightprec\app
{\forallvari m{}}{\boundvari x 1,\boundvari y 1}\und 
\boundvari y 1\tightprec\app
{\forallvari m{}}{\boundvari x 1,\boundvari y 1}}}
\footroom\oplistund
\forall\boundvari c{}\stopq
\forall\boundvari b{}\stopq\inparenthesesoplist{
{\color{blue}\forall\boundvari a{}\stopq
\inpit{
\boundvari a{}\tightprec\boundvari b{}
\und
\boundvari b{}\tightprec\boundvari c{}
\implies
\boundvari a{}\tightprec\boundvari c{}}}
\oplistund
{\color{green}\forall\boundvari a 1\stopq
\inpit{
\boundvari a 1\tightprec\boundvari b{}
\und
\boundvari b{}\tightprec\boundvari c{}
\implies
\boundvari a 1\tightprec\boundvari c{}}}}
\oplistimplies\headroom
{\color{red}\forall\forallvari u{}\stopq
\forall\forallvari v{}\stopq\forall\forallvari w{}\stopq}
\exists\boundvari z{}\stopq\inpit{
{\color{orange}\forallvari u{}}\tightprec\boundvari z{}\und
{\color{orange}\forallvari v{}}\tightprec\boundvari z{}\und
{\color{orange}\forallvari w{}}\tightprec\boundvari z{}}
}}{}{\Large\math{(R)}}
\par\halftop\halftop\noindent Note that 
---~to keep our formula rectified~---
we have added the index \nlbmath 1
to the copies and also renamed the bound \math\delta-variable
\app{\forallvari m{}}{\boundvari x{},\boundvari y{}} consistently.%
\halftop
\subsubsection
{Phase\,3 of the Reduction: Generalized Rules of Quantification}
\halftop\noindent
Now we have to apply the generalized rules of 
\mbox{\deltaminus\hskip-.16em- and} 
\mbox{\math\gamma-quantification}
backward until 
no quantifiers remain. 
According to the side-condition of the 
\ruleofuniversaloidquantification, 
we \nolinebreak have to guarantee that every \mbox{\math\delta-quantifier} is 
removed \mbox{before} its variable has been used for the instantiation 
of \math\gamma-variables. 
According to the side-condition of the generalized rule of
\mbox{\math\gamma-quantification}, 
we \nolinebreak have to guarantee that every term substituted for a
\mbox{\math\gamma-variable} does not contain any 
variables bound in the scope of its \mbox{\math\gamma-quantifier}.
To satisfy these side-conditions, 
we \nolinebreak have to remove the accessible quantifiers
in the order given by the following procedure.

\begin{procedure}{\bf(Restricting the Order of the 
Applications of Generalized 
\\\getittotheright{Rules of \math\gamma- and 
 \deltaminus\nolinebreak\hskip-.16em\nolinebreak-Quantification 
 for Reduction to a Sub-Expansion)}}%
\label{procedure}%
\\\noindent Do the following 
for every natural number \nlbmath i, \hskip.2em
stepping (by \nlbmath 1) 
\mbox{from \math 1 to \nlbmath n}:\notop\begin{quote}\sloppy
Keep removing 
\begin{itemize}\notop\item
the accessible \math\gamma-quantifiers
whose \mbox{\math\gamma-variable} is to be
instantiated with a term \nlbmath t \mbox{with \nlbmaths{\CARD t\prec i}{}} \ 
(where variables whose names are \skolem\ terms are considered to have the
 height of the \skolem\ terms) \hskip.2em
\aswellas\noitem\item
the accessible \math\delta-quantifiers,\notop\end{itemize}
until no such quantifiers remain.
\getittotheright\qed\end{quote}\end{procedure}

\begin{remark}[Historical Correctness]\\
For restricting the order of the applications, 
\herbrandindex\herbrand\ considered only the case of a
reduction to an expansion, 
not to a proper sub-expansion. \hskip.2em
The order he described in the first part of 
the proof of his \fundamentaltheorem\ 
\cite[\litsectref{5.5.1}]{herbrand-PhD}
is \nolinebreak a \nolinebreak complete one, which 
---~for the case of a full expansion~---
is one of the orders our Procedure\,\ref{procedure} admits. \hskip.3em
If we applied \herbrand's description directly to sub-expansions, however,
there would be no guarantee that the side-conditions 
of the applications of the 
\ruleofuniversaloidquantification\
are satisfied. \hskip.2em
Nevertheless, Procedure\,\ref{procedure} just captures
the weakest condition to satisfy all side-conditions for the reason
why \herbrand's procedure satisfies them; \cfnlb\ \sectref
{section Interlude: Why Side-Conditions are Always Satisfied}.
\getittotheright\qed\end{remark}

\pagebreak\par\indent
Let us see what our Procedure\,\ref{procedure} means in case of our example.
\par\yestop\noindent
For \math{i\tightequal 1}, \hskip.2em
there cannot be any 
\math\gamma-quantifiers satisfying the condition, \hskip.1em
because no term can have height \nlbmaths 0. \hskip.3em
Thus,
we apply the \ruleofuniversaloidquantification\ 
thrice backward to the last line, and obtain:
\par\halftop\halftop\noindent\LINEmaths{\noparenthesesoplist{
{\color{red}\forall\boundvari x{}\stopq\forall\boundvari y{}\stopq}
\exists\app{\forallvari m{}}{\boundvari x{},\boundvari y{}}\stopq
\inpit{
\boundvari x{}\tightprec\app
{\forallvari m{}}{\boundvari x{},\boundvari y{}}\und 
\boundvari y{}\tightprec\app
{\forallvari m{}}{\boundvari x{},\boundvari y{}}}
\oplistund
{\color{red}\forall\boundvari x 1\stopq}\forall\boundvari y 1\stopq
\exists\app{\forallvari m{}}{\boundvari x 1,\boundvari y 1}\stopq
\inpit{
\boundvari x 1\tightprec\app
{\forallvari m{}}{\boundvari x 1,\boundvari y 1}\und 
\boundvari y 1\tightprec\app
{\forallvari m{}}{\boundvari x 1,\boundvari y 1}}
\footroom\oplistund
\forall\boundvari c{}\stopq
\forall\boundvari b{}\stopq\inparenthesesoplist{
\forall\boundvari a{}\stopq
\inpit{
\boundvari a{}\tightprec\boundvari b{}
\und
\boundvari b{}\tightprec\boundvari c{}
\implies
\boundvari a{}\tightprec\boundvari c{}}
\oplistund
\forall\boundvari a 1\stopq
\inpit{
\boundvari a 1\tightprec\boundvari b{}
\und
\boundvari b{}\tightprec\boundvari c{}
\implies
\boundvari a 1\tightprec\boundvari c{}}}
\oplistimplies\headroom
\exists\boundvari z{}\stopq\inpit{
{\color{blue}\forallvari u{}}\tightprec\boundvari z{}\und
{\color{blue}\forallvari v{}}\tightprec\boundvari z{}\und
{\color{blue}\forallvari w{}}\tightprec\boundvari z{}}}
}{}{\Large\math{(Q)}}\begin{sloppypar}
\par\halftop\halftop\noindent
Note that the introduced free variables all
have height \nlbmaths i. \hskip.3em

Now,
because no accessible \mbox{\math\delta-quantifiers} remain,
we increment \math i by \nlbmaths 1.\end{sloppypar}
\par\yestop\noindent For \math{i\tightequal 2}, \hskip.2em
the only accessible \math\gamma-quantifiers 
that are to be \mbox{instantiated} with terms of 
height smaller than \nlbmath 2 \hskip.1em
are \math x and \math y in the first line
(according to the substitution \nlbmaths{\sigma_{1,1}}), \hskip.2em
\aswellas\ \math{x_1} of the second line 
(according to the substitution \nlbmaths{\sigma_{1,2}}).

Therefore, 
we apply \nlbmath{\sigma_{1,1}} to all occurrences of \nlbmath x and \nlbmath y
in the first line and remove their \mbox{\math\gamma-quantifiers}. \
Formally this means that
we rename the bound \math\delta-variable
\app{\forallvari m{}}{\boundvari x{},\boundvari y{}} to 
\app{\forallvari m{}}{\forallvari v{},\forallvari w{}}, \hskip.2em
and then we apply the generalized rule of \math\gamma-quantification
twice (backward).
Doing the same to \boundvari x 1 according to the instantiation of 
the original \math x in \nlbmaths{\sigma_{1,2}}, \hskip.2em
we then obtain the following formula:
\par\halftop\halftop\noindent\LINEmaths{\noparenthesesoplist{
{\color{red}\exists\app{\forallvari m{}}{\forallvari v{},\forallvari w{}}\stopq}
\inpit{{\color{green}\forallvari v{}}\tightprec
{\color{orange}\app{\forallvari m{}}{\forallvari v{},\forallvari w{}}}\und 
{\color{green}\forallvari w{}}\tightprec
{\color{orange}\app{\forallvari m{}}{\forallvari v{},\forallvari w{}}}}
\oplistund
\forall\boundvari y 1\stopq
\exists\app{\forallvari m{}}{\forallvari u{},\boundvari y 1}\stopq
\inpit{
{\color{green}\forallvari u{}}\tightprec\app
{\forallvari m{}}{\forallvari u{},\boundvari y 1}\und 
\boundvari y 1\tightprec\app
{\forallvari m{}}{\forallvari u{},\boundvari y 1}}
\footroom\oplistund
\forall\boundvari c{}\stopq
\forall\boundvari b{}\stopq\inparenthesesoplist{
\forall\boundvari a{}\stopq
\inpit{
\boundvari a{}\tightprec\boundvari b{}
\und
\boundvari b{}\tightprec\boundvari c{}
\implies
\boundvari a{}\tightprec\boundvari c{}}
\oplistund
\forall\boundvari a 1\stopq
\inpit{
\boundvari a 1\tightprec\boundvari b{}
\und
\boundvari b{}\tightprec\boundvari c{}
\implies
\boundvari a 1\tightprec\boundvari c{}}}
\oplistimplies\headroom
\lastlineofherbrandformula}}{}{\Large\math{(P)}}
\par\halftop\halftop\noindent Now,
the only accessible \mbox{\math\delta-quantifier}
is ``\app{\exists\forallvari m{}}{\forallvari v{},\forallvari w{}}.''
in the first line. 
Let us remove it by a backward application
of the \ruleofuniversaloidquantification, 
and obtain:
\par\halftop\halftop\noindent\LINEmaths{\noparenthesesoplist{
\forallvari v{}\tightprec
{\color{blue}\app{\forallvari m{}}{\forallvari v{},\forallvari w{}}}\und 
\forallvari w{}\tightprec
{\color{blue}\app{\forallvari m{}}{\forallvari v{},\forallvari w{}}}
\oplistund
{\color{red}\forall\boundvari y 1\stopq}
\exists\app{\forallvari m{}}{\forallvari u{},\boundvari y 1}\stopq
\inpit{
\forallvari u{}\tightprec\app
{\forallvari m{}}{\forallvari u{},\boundvari y 1}\und 
\boundvari y 1\tightprec\app
{\forallvari m{}}{\forallvari u{},\boundvari y 1}}
\footroom\oplistund
\forall\boundvari c{}\stopq
\forall\boundvari b{}\stopq\inparenthesesoplist{
\forall\boundvari a{}\stopq
\inpit{
\boundvari a{}\tightprec\boundvari b{}
\und
\boundvari b{}\tightprec\boundvari c{}
\implies
\boundvari a{}\tightprec\boundvari c{}}
\oplistund
\forall\boundvari a 1\stopq
\inpit{
\boundvari a 1\tightprec\boundvari b{}
\und
\boundvari b{}\tightprec\boundvari c{}
\implies
\boundvari a 1\tightprec\boundvari c{}}}
\oplistimplies\headroom
\lastlineofherbrandformula}}{}{\Large\math{(O)}}
\par\halftop\halftop\noindent
Note that the variable of the removed \math\delta-quantifier has 
height \nlbmath i again.

Now,
because all accessible \mbox{\math\gamma-quantifiers}
are to be replaced with 
terms of height not smaller than \nlbmaths 2, \hskip.2em
\mbox{we increment \math i} by \nlbmaths 1.
\vfill\pagebreak\par\indent
For \math{i\tightequal 3}, \hskip.1em
according to Procedure\,\ref{procedure}, \hskip.1em
we now have to apply the
generalized rule of \mbox{\math\gamma-quantification} 
backward to the 
second line according to the substitution \nlbmath{\sigma_{1,2}}
(for the renamed \boundvari y 1 instead of the
 original \nlbmath{\boundvari y{}})
(after the involved renaming of the following
\math\delta-variable),
because every term in the range of \math{\sigma_{1,2}} has a height 
smaller than \nlbmaths 3. \ 
\mbox{We obtain}: 
\par\halftop\noindent\LINEmaths{\noparenthesesoplist{
{\forallvari v{}\tightprec\app
{\forallvari m{}}{\forallvari v{},\forallvari w{}}\und 
\forallvari w{}\tightprec\app
{\forallvari m{}}{\forallvari v{},\forallvari w{}}}
\oplistund{\color{red}\exists\termzwei\stopq}
\oplistnl
\mbox{~\mbox{}~\mbox{}~~~~~~}\inpit{
\forallvari u{}\tightprec
{\color{orange}\app{\forallvari m{}}{\forallvari u{},
{\app{\forallvari m{}}{\forallvari v{},\forallvari w{}}}}}\und 
{\color{green}
 \app{\forallvari m{}}{\forallvari v{},\forallvari w{}}}\tightprec
{\color{orange}\app
{\forallvari m{}}{\forallvari u{},
{\app{\forallvari m{}}{\forallvari v{},\forallvari w{}}}}}}
\footroom\oplistund
\forall\boundvari c{}\stopq
\forall\boundvari b{}\stopq\inparenthesesoplist{
\forall\boundvari a{}\stopq
\inpit{
\boundvari a{}\tightprec\boundvari b{}
\und
\boundvari b{}\tightprec\boundvari c{}
\implies
\boundvari a{}\tightprec\boundvari c{}}
\oplistund
\forall\boundvari a 1\stopq
\inpit{
\boundvari a 1\tightprec\boundvari b{}
\und
\boundvari b{}\tightprec\boundvari c{}
\implies
\boundvari a 1\tightprec\boundvari c{}}}
\oplistimplies\headroom
\lastlineofherbrandformula}}{}{\Large\math{(L)}}
\par\halftop\noindent After a subsequent application of the 
\ruleofuniversaloidquantification, 
we \nolinebreak obtain:
\par\halftop\noindent\LINEmaths{\noparenthesesoplist{
{\forallvari v{}\tightprec\app
{\forallvari m{}}{\forallvari v{},\forallvari w{}}\und 
\forallvari w{}\tightprec\app
{\forallvari m{}}{\forallvari v{},\forallvari w{}}}
\oplistund
\forallvari u{}\tightprec
{\color{blue}\app
{\forallvari m{}}{\forallvari u{},{\app{\forallvari m{}}{\forallvari v{},\forallvari w{}}}}}\und 
{\app{\forallvari m{}}{\forallvari v{},\forallvari w{}}}\tightprec
{\color{blue}\app{\forallvari m{}}{\forallvari u{},
{\app{\forallvari m{}}{\forallvari v{},\forallvari w{}}}}}
\footroom\oplistund
{\color{red}\forall\boundvari c{}\stopq
\forall\boundvari b{}\stopq}\inparenthesesoplist{
{\color{red}\forall\boundvari a{}\stopq}
\inpit{
\boundvari a{}\tightprec\boundvari b{}
\und
\boundvari b{}\tightprec\boundvari c{}
\implies
\boundvari a{}\tightprec\boundvari c{}}
\oplistund
{\color{red}\forall\boundvari a 1\stopq}
\inpit{
\boundvari a 1\tightprec\boundvari b{}
\und
\boundvari b{}\tightprec\boundvari c{}
\implies
\boundvari a 1\tightprec\boundvari c{}}}
\oplistimplies\headroom
{\color{red}\exists\boundvari z{}\stopq}\inpit{
\forallvari u{}\tightprec\boundvari z{}\und
\forallvari v{}\tightprec\boundvari z{}\und
\forallvari w{}\tightprec\boundvari z{}}}}{}{\Large\math{(K)}}
\par\halftop\noindent Note that the variable of the 
removed \math\delta-quantifier has 
height \nlbmath i again.

Moreover, 
note that the free variable 
\app{\forallvari m{}}{\forallvari v{},\forallvari w{}}
occurs twice in the first line but only once in second line,
because the remaining two occurrences are just part of the name 
of the free variable \app
{\forallvari m{}}
{\forallvari u{},{\app{\forallvari m{}}{\forallvari v{},\forallvari w{}}}}. \
Be reminded that 
every \skolem\ function symbol starts a name of a variable.%

Now, as all accessible quantifiers in the formula \nlbmath K are 
\mbox{\math\gamma-quantifiers}
whose variables are to be replaced with 
terms of height not smaller than \nlbmaths 3,
\mbox{we increment \math i} by \nlbmath 1 again.
\subsubsection{Interlude: Why Side-Conditions are Always Satisfied}\label
{section Interlude: Why Side-Conditions are Always Satisfied}%
\subsubsubsection{\protect\Ruleofuniversaloidquantification}\mbox{}
\par\noindent
Note that the side-condition of the 
\ruleofuniversaloidquantification\
(\ie\nolinebreak\ that its bound \mbox{\math\delta-variable} 
 does not occur free outside its scope) \hskip.1em
is observed in the last \nolinebreak step from the formula \nlbmath L \hskip.1em 
back to the formula \nlbmaths K. \hskip.2em

Let us show that 
\mbox{---~under} 
the reasonable conditions explained below on the terms introduced by 
backward applications of the generalized rule of 
\mbox{\math\gamma-quantification~---}
this is always the case for the order of removal of quantifiers 
given by Procedure\,\ref{procedure}:

If a bound \mbox{\math\delta-variable} is turned into a free
variable by the removal of its quantifier by reductive application of
the 
\ruleofuniversaloidquantification\ 
according to Procedure\,\ref{procedure}
during the removal of quantifiers for a given \nlbmaths i, \hskip.15em
then its height (seen as a term) must be \nlbmaths i. \hskip.4em

\pagebreak

Indeed, for \maths{i\tightequal 1},
the \skolem\ function symbol of the variable name 
whose \math\delta-quantifier is removed
must be a constant, \hskip.2em
simply because it cannot occur in the scope of \math\gamma-quantifiers. 
\hskip.3em
Moreover, 
for \maths{i\tightsucc 1},
its \skolem\ function symbol 
must have at least one argument
of height \nlbmaths{i\tight-1}, \hskip.2em
simply because otherwise 
its \mbox{\math\delta-quantifier} would have been removed 
 already for a smaller \nlbmaths i.

Here it is crucial that we took the
\outerskolemizedform\ instead of the \deltaplusplusskolemizedform: \hskip.3em
All \math\gamma-variables in whose scope the \math\delta-variable
occurs must be arguments of the \skolem\ term that
names the \math\delta-variable.

Now it is immediate that the \math\delta-variable cannot have been
introduced as a free variable by a previous backward
application of a generalized rule of \math\gamma-quantification, \hskip.2em
simply because all such free variables have a height strictly smaller 
than \nlbmaths i.

On the other hand, if the \math\delta-variable was
introduced as a free variable 
by a previous backward
application of a \ruleofuniversaloidquantification,
then 
---~because it must have the identical \skolem\ symbol~---
this previous application must have occurred in a parallel 
branch resulting from our initial backward applications of the
generalized rule of \mbox{\math\gamma-simplification.} \hskip.3em
If,
however,
we \nolinebreak never instantiate the \math\gamma-variables of
the top quantifiers of  
different copies introduced by backward applications of the 
generalized rule of \mbox{\math\gamma-simplification} with identical instances
(and there would be no benefit of this), \hskip.2em
then we \nolinebreak know 
that the \skolem\ terms that name the bound \math\delta-variables
have different arguments. 

Here \mbox{---~for} different instances of the \math\gamma-variables 
to result in different \mbox{\skolem\ terms~---}
it is again crucial to have the \deltaminus- instead
of the \deltaplusplusskolemizedform.

Therefore, 
the side-condition of the \ruleofuniversaloidquantification\
(\ie\nolinebreak\ that its bound \mbox{\math\delta-variable} 
 does not occur free outside its scope) \hskip.1em
is \nolinebreak always observed under reasonable conditions 
on the terms introduced by 
backward applications of the generalized rule of 
\mbox{\math\gamma-quantification.}

\subsubsubsection{Generalized rule of \math\gamma-quantification}\mbox{}
\begin{sloppypar}\par\noindent
In a similar way, we can see that also the side-condition of the generalized rule
of \mbox{\math\gamma-quantification}
---~namely that none of the free variables occurring in the
term \nlbmath t substituted for its bound \math\gamma-variable
is bound by quantifiers in the scope of the 
\mbox{\math\gamma-quantification~---}
is \nolinebreak always satisfied under reasonable conditions:\end{sloppypar}

If we do not admit a 
variable that does not start with a \skolem\ function symbol
in \math t unless it occurs free in \nlbmath A
(and there would be no benefit of this), \hskip.2em
then the \math\gamma-quantifications in the scope of the
original \mbox{\math\gamma-quantification} cannot bind such a variable,
simply because \math A is assumed to be rectified.

On the other hand, all \mbox{\math\delta-quantified} variables 
in the scope of the original \mbox{\math\gamma-quantification} 
contain \nlbmath t as a subterm of their names
and thus have a bigger height than all free variables
occurring in \nlbmaths t.

Here again it is  crucial to have the \deltaminus- instead
of the \deltaplusplusskolemizedform.

After these general considerations on the side-conditions of the 
generalized rules of quantification, \hskip.1em 
let us return to our example proof.

\vfill\pagebreak

\subsubsection{Continuation of Phase\,3}\label
{section Continuation of Phase 3}%
For \math{i\tightequal 4} \
(\ie\ for \math{i\tightequal n}), \
we remove the five \math\gamma-quantifiers 
remaining in the formula \nlbmath K 
by the generalized
rule of \math\gamma-quantification, 
using \math{\sigma_{2,1}} for the
third line,
\math{\sigma_{2,2}} \nolinebreak for the fourth line, and
\math{\sigma_3} for the last line,
and obtain the following formula:
\par\halftop\noindent\LINEmaths{\noparenthesesoplist{
{\forallvari v{}\tightprec\app
{\forallvari m{}}{\forallvari v{},\forallvari w{}}\und 
\forallvari w{}\tightprec\app
{\forallvari m{}}{\forallvari v{},\forallvari w{}}}
\oplistund
{
\forallvari u{}\tightprec\app
{\forallvari m{}}{\forallvari u{},{\app{\forallvari m{}}{\forallvari v{},\forallvari w{}}}}\und 
{\app{\forallvari m{}}{\forallvari v{},\forallvari w{}}}\tightprec\app
{\forallvari m{}}{\forallvari u{},{\app{\forallvari m{}}{\forallvari v{},\forallvari w{}}}}}
\oplistund 
\inparenthesesoplist{{\color{blue}\forallvari v{}}\tightprec
{\color{blue}\termeins}
\und{\color{blue}\termeins}\tightprec{\color{green}\termzwei}
\oplistimplies{\color{blue}\forallvari v{}}\tightprec{\color{green}\termzwei}}
\oplistund
\inparenthesesoplist{{\color{blue}\forallvari w{}}\tightprec
{\color{blue}\termeins}
\und{\color{blue}\termeins}\tightprec{\color{green}\termzwei}
\oplistimplies{\color{blue}\forallvari w{}}\tightprec{\color{green}\termzwei}}
\oplistimplies 
\inparenthesesoplist{
\forallvari u{}\tightprec{\color{green}\termzwei}\oplistund
\forallvari v{}\tightprec{\color{green}\termzwei}\oplistund
\forallvari w{}\tightprec{\color{green}\termzwei}}}}{}{\Large\math{(J)}}
\begin{sloppypar}\halftop\par\noindent
Now, 
considering all atomic formulas of \nlbmath J as names of sentential variables,
we \nolinebreak easily see that 
the formula \math J is a \sententialtautology. \
Moreover, \math J is a sub-expansion of the formula \nlbmath F \wrt\ 
\nlbmaths{\termsofdepth 4 F}, \hskip.2em
the champ fini of order \nlbmaths 4, \hskip.2em
formed over the function and free variable symbols of \nlbmath F,
\ie\ over \forallvari u{}, \forallvari v{}, \forallvari w{}, and
\nlbmaths{\app{\forallvari m{}}{\_,\_}}. \hskip.5em
Thus, \hskip.2em
by \cororefs{corollary sub-expansion}{corollary sub-expansion two}, \hskip.2em
the expansion \nlbmath{F^{\,\termsofdepth 4 F}} 
is a \sententialtautology\ as \nolinebreak well, \hskip.2em
and so \math A has \propertyC\ of order \nlbmaths 4.%
\end{sloppypar}%
%

\subsubsection{Stepwise {\em Deductive}\/ Construction of the Same Proof?}
Suppose we have found a sentential tautology that is a sub-expansion of 
the \outerskolemizedform\ \nlbmath F of \nlbmaths A. \hskip.3em
Can we now construct a proof of \nlbmath A deductively, \ie\ stepwise
from such a sentential tautology {\em forward}\/ to obtain formula \nlbmath A
in the end? \hskip.2em
Of course, the reductive construction shows that
there must be one such formula where this is possible, \hskip.2em
but even if we \nolinebreak found a formula where this is actually possible
(such as the formula \nlbmath J), \hskip.2em 
it still remains
much more difficult to find the steps deductively than reductively,
\hskip.1em for the following reasons:\begin{enumerate}\noitem\item
To find out, which of the variables named by \skolem\ terms
are to be turned into \math\gamma- and which in to \math\delta-quantifiers,
we have to check whether they were introduced by instantiation of \nlbmath F
(such as the occurrences of {\color{green}\termzwei\ in green} in \nlbmath J),
or \nolinebreak whether they resulted 
from a further instantiation of a \skolem\ term already
present in \nlbmath F 
(such as the occurrences of \termzwei\ in black in the \nth 2\,line of
\nlbmath J). \hskip.3em
The former have to be replaced with \math\gamma-quantifiers first,
so that the side-condition of the \ruleofuniversaloidquantification\
becomes satisfied for then turning the latter 
into \math\delta-quantifiers.\noitem\item
Before we can start with the introduction of \math\gamma-quantifiers for 
the instantiated terms with maximal depth 
(such as the occurrences of {\color{green}\termzwei\ in green} 
in the \nth 3 and \nth 4 lines of \nlbmath J, \hskip.1em
\ie\nolinebreak\ the quantifier
{\color{red}\math{\forall c.}} in formula \nlbmath K), \hskip.2em 
we \nolinebreak may first 
have to introduce the \math\gamma-quantifiers that occur
more innermost (such as the \mbox{\math\gamma-quantifiers} for the terms printed 
in {\color{blue}blue} in the formula \nlbmaths J, \hskip.1em
\ie\ the quantifiers 
{\color{red}\math{\forall a.}},
{\color{red}\math{\forall a_1.}}, 
{\color{red}\math{\forall b.}} in formula \nlbmath K).\noitem\end{enumerate}
\vfill\pagebreak
\subsection[Mechanical Proof Construction in \heijenoort's Version of\\\herbrand's {\it Modus Ponens}\/\hskip.07em-Free Calculus]{Mechanical Proof Construction in \heijenoort's Version of\\ \herbrand's {\it Modus Ponens}\/\hskip.07em-Free Calculus}
Knowing that \math A has \propertyC\ of order \nlbmaths 4, 
we can construct a proof of \nlbmath A in our correction of 
\heijenoortindex\heijenoort's version of 
\herbrandindex\herbrand's {\it modus ponens}\/\hskip.07em-free calculus
also {\em mechanically}.

Let \math N denote the cardinality of 
\nlbmaths{\termsofdepth 4 F}. \ 
Let \bigmaths{\termsofdepth 4 F=\{t_0,\ldots,t_{N-1}\}}. 
We \nolinebreak have \bigmaths{N=3+3^2+\inpit{3+3^2}^2=156}. \
Therefore, the number of leaf formulas 
found in the expansion \nlbmath{F^{\,\termsofdepth 4 F}}
of the form of the three original ones
of \nlbmaths F,
summed over the three original lines, is \maths{
N^2
+
N^3
+
N
=
3820908}.\par\indent
We can now reduce \math A mechanically to the huge\footnote{%
 \label{note disjunction versus expansion}%
 Note, however, that the {\em\herbrand\ disjunction}\/ used for \nlbmath A
 in~\makeaciteoftwo{herbrand-handbook}{SR--2009--01} \hskip.1em
 (instead of the expansion of the \outerskolemizedform) \hskip.1em
 is even more than a factor of \nlbmath{10^7} bigger.%
}
\sententialtautology\
\nlbmath{F^{\,\termsofdepth 4 F}}
as follows:
In phase\,1,
we rename every bound \mbox{\math\delta-variable} to 
its \skolem\ term of the \outerskolemizedform\ of \nlbmaths F. \hskip.4em
In phase\,2,
for every \math\gamma-quantifier, 
we apply the generalized rule of \mbox{\math\gamma-simplification}
\math{N\tight-1}~times
(backward)
and associate the terms \math{t_0,\ldots,t_{N-1}} to the resulting
\math N \math\gamma-quantifiers. \
Finally,
in phase\,3,
we apply Procedure\,\ref{procedure}
already described in \sectref{section formal proof}
for the reductive (\ie\ backward) 
application of the generalized rules
of \mbox{\deltaminus- and}
\mbox{\math\gamma-quantification}, 
including the renaming of the bound \math\delta-variables. \hskip.1em
When all quantifiers are removed, 
we \nolinebreak have obtained \nlbmaths{F^{\,\termsofdepth 4 F}\!}, \hskip.2em
and the formal proof is done.

Moreover, because the 
renaming of the bound \math\delta-variables is actually determined 
by the path on which they occur in the formula,
the final renaming of bound variables 
can be done right after the backward applications of the
generalized rule of \math\gamma-simplification.
Because our version of the generalized rule of 
\math\gamma-simplification admits variants of \maths H, \hskip.2em
we can therefore
restrict the 
renaming bound variables 
to 
the very beginning of a
reductive proof of
\nlbmaths A.

This result can be used to give a constructive proof of 
the following lemma.

\yestop
\begin{lemma}
[From \propertyC\ to a Linear Derivation]
\label{lemma from C to yields a la heijenoort}\par\noindent
Let\/ \math A be a rectified \firstorder\ formula. \
Let\/ \math F be the \outerskolemizedform\ of \nlbmath A.
\\Let\/ \math n be a positive natural number.
\par\noindent If\/ \math A has \propertyC\ of order \nlbmath n, 
then we can construct a derivation of \nlbmath A
of the following form, 
in which we read any term starting with a\/ \skolem\ function 
as an atomic variable:\par\noindent\begin{tabular}[b]{@{}l l}\majorheadroom
  {\bf Step\,1:} 
 &We \nolinebreak start with 
  the \sententialtautology\/ \nlbmath{F^{\,\termsofdepth n F}}
\\
 &(or else with 
some
sub-expansion of \nlbmaths F{}
   \wrt\ \nlbmaths{\termsofdepth n F}{}
that is a \sententialtautology).
\\\majorheadroom{\bf Step\,2:}
 &Then we may repeatedly apply the generalized rules of\/
  \math\gamma- and 
  \deltaminus\nolinebreak\hskip-.16em\nolinebreak-quanti\-fi\-cation.
\\\majorheadroom{\bf Step\,3:}
 &Then we may repeatedly apply the generalized rule of\/
  \math\gamma-simplification.
\\\majorheadroom{\bf Step\,4:}
 &Then we rename all bound \mbox{\math\delta-variables}
  to obtain \nlbmath A. 
\\\end{tabular}
\\[-2.6ex]\getittotheright\qed\end{lemma}\vfill\pagebreak
\subsection{Formal Proof in Our Free-Variable Calculus}\label
{section Formal Proof in Our Free-Variable Calculus}%
In this section we do the proof of \sectref{section formal proof}
again, but in our free-variable calculus. \ 
So we replace the
\ruleofuniversaloidquantification\ with our
generalized rule of \deltaplusplus-quantification.\footnotemark\ \
The advantage is not only the possible delay in 
the choice of witnessing terms,\footnotemark\
but also that Procedure\,\ref{procedure} of \sectref{section formal proof}
for restricting the order of 
applications of the generalized rules of quantification 
becomes superfluous and any order will do. \hskip.1em
Another difference is that now all \skolem\ functions are function
symbols of the first-order calculus and do not start the name of a
variable anymore. \
By application of the generalized rule of \math\gamma-simplification
just as before twice, we reduce \math A to the following formula:
\par\noindent\LINEmaths{\noparenthesesoplist{
\forall\boundvari x{}\stopq
\forall\boundvari y{}\stopq
\exists\boundvari m{}\stopq
\inpit{\boundvari x{}\tightprec\boundvari m{}\und 
\boundvari y{}\tightprec\boundvari m{}}
\oplistund
\forall\boundvari x 1\stopq
\forall\boundvari y 1\stopq
\exists\boundvari m 1\stopq
\inpit{\boundvari x 1\tightprec\boundvari m 1\und 
\boundvari y 1\tightprec\boundvari m 1}
\footroom\oplistund
\forall\boundvari c{}\stopq
\forall\boundvari b{}\stopq\inparenthesesoplist{
\forall\boundvari a{}\stopq
\inpit{
\boundvari a{}\tightprec\boundvari b{}
\und
\boundvari b{}\tightprec\boundvari c{}
\implies
\boundvari a{}\tightprec\boundvari c{}}
\oplistund
\forall\boundvari a 1\stopq
\inpit{
\boundvari a 1\tightprec\boundvari b{}
\und
\boundvari b{}\tightprec\boundvari c{}
\implies
\boundvari a 1\tightprec\boundvari c{}}}
\oplistimplies\headroom
\forall\boundvari u{}\stopq
\forall\boundvari v{}\stopq
\forall\boundvari w{}\stopq
\exists\boundvari z{}\stopq
\inpit{
\boundvari u{}\tightprec\boundvari z{}\und
\boundvari v{}\tightprec\boundvari z{}\und
\boundvari w{}\tightprec\boundvari z{}}
}}{}{\Large\math{(R')}}
\par\noindent Note that \math{R'} is 
similar to the formula \nlbmath R of \sectref{section formal proof}, 
but lacks the renaming of bound \mbox{\math\delta-variables}
to \skolem\ terms.
As the order of applications of generalized rules of quantification 
does not matter anymore,
let us remove all consecutively accessible quantifiers of the same kind 
and obtain:
\\[-.5ex]\noindent\LINEmaths{\noparenthesesoplist{
\exists\boundvari m{}\stopq
\inpit{\existsvari x 0\tightprec\boundvari m{}\und 
\existsvari y 0\tightprec\boundvari m{}}
\oplistund
\exists\boundvari m 1\stopq
\inpit{\existsvari x 1\tightprec\boundvari m 1\und 
\existsvari y 1\tightprec\boundvari m 1}
\footroom\oplistund
\inpit{
\existsvari a 0\tightprec\existsvari b{}
\und
\existsvari b{}\tightprec\existsvari c{}
\implies
\existsvari a 0\tightprec\existsvari c{}}
\oplistund
\inpit{
\existsvari a 1\tightprec\existsvari b{}
\und
\existsvari b{}\tightprec\existsvari c{}
\implies
\existsvari a 1\tightprec\existsvari c{}}
\oplistimplies\headroom
\lastlineofherbrandformula
}}{}{\Large\math{(Q')}}
\par\noindent If we denote\footnotemark\
\bigmaths{\inparenthesestight{\exists\boundvari m{}\stopq
\inpit{\existsvari x 0\tightprec\boundvari m{}\und 
\existsvari y 0\tightprec\boundvari m{}}\comma
\existsvari x 0\comma\existsvari y 0}^\delta}{}
with \nlbmaths{\forallvari m{}}, \hskip.2em 
then by application of the generalized rule of 
\math{\deltaplusplus\!}-quantification 
backward to the first line we obtain:
\par\noindent\LINEmaths{\noparenthesesoplist{
\existsvari x 0\tightprec
\app{\forallvari m{}}{\existsvari x 0,\existsvari y 0}\und 
\existsvari y 0\tightprec
\app{\forallvari m{}}{\existsvari x 0,\existsvari y 0}
\oplistund
\exists\boundvari m 1\stopq
\inpit{\existsvari x 1\tightprec\boundvari m 1\und 
\existsvari y 1\tightprec\boundvari m 1}
\footroom\oplistund
\inpit{
\existsvari a 0\tightprec\existsvari b{}
\und
\existsvari b{}\tightprec\existsvari c{}
\implies
\existsvari a 0\tightprec\existsvari c{}}
\oplistund
\inpit{
\existsvari a 1\tightprec\existsvari b{}
\und
\existsvari b{}\tightprec\existsvari c{}
\implies
\existsvari a 1\tightprec\existsvari c{}}
\oplistimplies\headroom
\lastlineofherbrandformula
}}{}{\Large\math{(L')}}
\par\noindent Because \bigmaths{\exists\boundvari m 1\stopq
\inpit{\existsvari x 1\tightprec\boundvari m 1\und 
\existsvari y 1\tightprec\boundvari m 1}}{} results
from \bigmaths{\exists\boundvari m{}\stopq
\inpit{\existsvari x 0\tightprec\boundvari m{}\und 
\existsvari y 0\tightprec\boundvari m{}}}{} by renaming of bound variables 
and by a bijective renaming of
free \math\gamma-variables via
\bigmaths{\{\existsvari x 0\tight\mapsto\existsvari x 1\comma
\existsvari y 0\tight\mapsto\existsvari y 1\}},
we can use the same function symbol \forallvari m{}
again for the second line. \hskip.3em
If we also remove the \math\gamma-quantifier
in the last line, \hskip.1em
we obtain:
\par\noindent\LINEmaths{\noparenthesesoplist{
\existsvari x 1\tightprec
\app{\forallvari m{}}{\existsvari x 1,\existsvari y 1}\und 
\existsvari y 1\tightprec
\app{\forallvari m{}}{\existsvari x 1,\existsvari y 1}
\oplistund
\existsvari x 2\tightprec
\app{\forallvari m{}}{\existsvari x 2,\existsvari y 2}\und 
\existsvari y 2\tightprec
\app{\forallvari m{}}{\existsvari x 2,\existsvari y 2}
\footroom\oplistund
\inpit{
\existsvari a 1\tightprec\existsvari b 1
\und
\existsvari b 1\tightprec\existsvari c 1
\implies
\existsvari a 1\tightprec\existsvari c 1}
\oplistund
\inpit{
\existsvari a 2\tightprec\existsvari b 2
\und
\existsvari b 2\tightprec\existsvari c 2
\implies
\existsvari a 2\tightprec\existsvari c 2}
\oplistimplies\headroom
\forallvari u{}\tightprec\existsvari z{}\und
\forallvari v{}\tightprec\existsvari z{}\und
\forallvari w{}\tightprec\existsvari z{}
}}{}{\Large\math{(J')}}
\par\noindent Finally, 
the \sententialtautology\ \nlbmath J 
of \nlbsectref{section Continuation of Phase 3}
results from 
\nlbmath{J'} by an appropriate instantiation of free \math\gamma-variables.%
\pagebreak\par
\addtocounter{footnote}{-2}%
\footnotetext{\Cfnlb\ \sectref{section delta-Quantification}.}%
\addtocounter{footnote}{1}%
\footnotetext{%
 \majorheadroom
 \Cfnlb\ \sectref{section gamma-Quantification}.}%
\addtocounter{footnote}{1}%
\footnotetext{%
 \majorheadroom
 For the function \nlbmaths\delta, 
 see the generalized rule of \deltaplusplus-quantification in
 \sectref{section delta-Quantification}, \hskip.2em
 and also \defiref{definition skolemized forms} in
 \nlbsectref{section Prerequisites for Property C}.%
}%
This example proof admits generalization to the following lemma.

\yestop
\begin{lemma}
[From \propertyCstar\ to a Linear Derivation]
\label{lemma from C star to yields a la heijenoort}\par\noindent
Let\/ \math A be a rectified \firstorder\ formula. \
Let\/ \math F be the \deltaplusplusskolemizedform\ of \nlbmaths A.
\\Let\/ \math n be a positive natural number.
\par\noindent If\/ \math A has \propertyCstar\ of order \nlbmath n, 
\\then we can construct a derivation of \nlbmath A
of the following form:
\par\noindent\begin{tabular}[b]{@{}l l}\majorheadroom
  {\bf Step\,1:}
 &We \nolinebreak start with a quantifier-free formula \nlbmath B 
  with free \math\gamma-variables, such that
\\
 &---~by replacing each free \math\gamma-variable \nlbmath{\existsvari x{}}
  occurring in \nlbmath B \hskip.1em
  with a
\\
 &~~~~term \nlbmath{t_{\existsvari x{}}} 
  (over the free symbols in \nlbmath F and possibly \nlbmath\bullet) 
  \hskip.2em
  with \bigmaths{\CARD{t_{\existsvari x{}}}\prec n}{}~---
\\
 &we obtain the \sententialtautology\ \nlbmath
  {F^{\,\termsofdepth n F}}
\\
 &(or else 
some
sub-expansion of \nlbmath F 
   \wrt\ \nlbmaths{\termsofdepth n F}{} that is a \sententialtautology).
\\\majorheadroom{\bf Step\,2:}
 &Then we may repeatedly apply the generalized rules of\/
\\
 &\mbox{}\hskip.1em\deltaplusplus\hskip-.2em-quantification
   and restricted \math\gamma-quantification.
\\\majorheadroom{\bf Step\,3:} 
 &Then we may repeatedly apply the generalized rule of\/ 
\\&\mbox{}\hskip.13em
  \math\gamma-simplification to obtain \nlbmaths A.
\\\end{tabular}%
\\[-2.6ex]\getittotheright\qed\end{lemma}
\vfill\pagebreak
\section
{\herbrandsfundamentaltheorem
}\label{section fundamental theorem}%
If we now have a look on what our inference rules do with a formula
that has \propertyC\ 
of order \nlbmath n, \hskip.2em
then we see that the results of applications of 
the generalized rules of 
simplification and
\deltaminus\nolinebreak\hskip-.16em\nolinebreak-quantification 
(and renaming of bound variables) 
are formulas that have
\propertyC\ 
of order \nlbmath n again.\footnote{
 Indeed, regarding the \ruleofuniversaloidquantification,
 the \skolem ization
 just removes the \mbox{\math\delta-quantifier} and replaces its bound variable
 with a constant that occurs at exactly the same places where the 
 variable occurred free before application of the 
 \ruleofuniversaloidquantification. \hskip.3em
 Moreover, regarding the generalized rule of simplification,
 the expansion of the \outerskolemizedform\ of
 \maths{A[H\circ H']}{} can be transformed into \math{S\oder S'}
 where \math S and \nlbmath{S'} are the
 expansions 
 of the \outerskolemizedform s
 of \math{A[H]} and \nlbmath{A[H']}, respectively,
 over \nolinebreak the \nolinebreak terms that include the 
 \skolem\ functions of \nlbmaths{A[H\circ H']}. \ 
 Clearly, all or none of \math S, \math{S'}, \mbox{\math{S\oder S'}}
 are \sententialtautologies, and the rest follows from 
 \cororef{corollary extend signature}.%
}
Moreover,
an application of the generalized rule of \mbox{\math\gamma-quantification}
where a term \nlbmath t is replaced with a bound \math\gamma-variable
in a formula that has \propertyC\ of order \nlbmath n \hskip.1em
results in 
a formula that has \propertyC\ of order \nlbmaths{n+\CARD t}.\footnote{%
 \majorheadroom
 Note that the replacement seemingly inverting the effect of the 
 application of the generalized rule of \math\gamma-quantification, 
 namely the removal of the new \math\gamma-quantifier and the replacement of 
 its bound variable with the term \nlbmath t would result only 
 in an order of \bigmaths{\max\{n\comma\CARD t+1\}}{}
 (because the free variables of \math t remain 
  free also after the substitution). \hskip.4em
 This is not a complete inversion, however, 
 because \math t may occur in the expansion also as 
 a newly added first argument of \skolem\ functions. \hskip.2em
 Indeed,
 the arity of all \skolem\ functions for the
 \mbox{\math\delta-quantifiers} in the 
 scope of the \mbox{\math\gamma-quantifier} may be increased; \hskip.2em
 this requires the order \nlbmath{n+\CARD t} in general. \hskip.3em
 Finally we have \bigmaths
 {\max\{n\comma\CARD t+1\comma n+\CARD t\}\ =\ n+\CARD t}.
 \par
 For examples of this, \hskip.1em
 consider the formulas of \hskip.1em\sectref{section formal proof}: \hskip.3em
 The formula \nlbmath L has \propertyC\ of order \nlbmath 2,\footnotemark\
 \hskip.2em
 but the formula \nlbmath O 
 \mbox{---~derived} from \nlbmath L by applications 
 the generalized rule of 
 \mbox{\math\gamma-quantification} 
 \mbox{with\addtocounter{footnote}{-1}\footnotemark\ \maths{\CARD t=1}{}~~---}
 \mbox{has \propertyC} of order \nlbmath 3. \hskip.4em
 Similarly,
 the formula \nlbmath P has \propertyC\ of order \nlbmath 3
 (it has to have the same order as formula \math O because it 
  results from it by application of the 
  \ruleofuniversaloidquantification), \hskip.2em
 but the formula \nlbmath Q 
 \mbox{---~derived} from \nlbmath P by applications the generalized rule of 
 \mbox{\math\gamma-quantification} 
 \mbox{with \maths{\CARD t=1}{}~~---}
 \mbox{has \propertyC} of order \nlbmaths 4. \hskip.4em
 In both examples, the order is not increased by inverse substitution
 of the term \nlbmath t,
 because then all formulas would have \propertyC\ invariantly of order 
 \bigmaths{\max\{n\comma \CARD t+1\}\ =\ \max{\{n\comma 2\}}\ =\ n}. 
 The order is increased because the depth of the \skolem\ term
 resulting from the previously introduced \math\delta-quantifier 
 is increased by a newly added first argument.%
}\footnotetext{%
 \majorheadroom
 In the peculiar notation of \sectref{section formal proof},
 where variables are denoted by terms with a \skolem\ function 
 as top symbol, this is a bit difficult to see: \hskip.2em
 To show that formula \nlbmath L has \propertyC\ of order \nlbmath 2,
 according to the substitution \nlbmath{\sigma_3},
 we have to replace the bound \mbox{\math\delta-variable} \nlbmath\termzwei\
 with the \skolem\ constant \nlbmaths{%
 \headroom
 \overline\termzwei}, \hskip.2em
 where the whole overlined term is to be interpreted
 as the name of the constant. \hskip.2em
 To show that formula \nlbmath O has \propertyC,
 we instantiate the bound \mbox{\math\delta-variable} 
 with the same term in principle,
 but now we have to denote this term in the form of \nlbmaths
 {\headroom\overline{\forallvari m{}(\forallvari u{},}\termeins\overline)},
 where \nlbmaths
 {\headroom\overline{\forallvari m{}(\forallvari u{},}\underline{~~}\overline)}{}
 is the \singulary\ \skolem\ function symbol introduced by the expansion,
 that was a nullary one before,
 now with the new argument \nlbmaths\termeins, \hskip.1em
 which is the name of a variable introduced by the generalized rule of
 \math\gamma-quantification.
 This term has height \nlbmaths 2; \hskip.2em
 so formula \nlbmath O has \propertyC\ of order \nlbmaths 3.%
} \hskip.4em
Using \cororef{corollary sentential tautology}, \hskip.15em
we \nolinebreak can summarize this as follows:
\begin{lemma}[From a Linear Derivation to \propertyC]
\label{lemma from yields to C a la heijenoort}\sloppy\\\noindent
If there is a derivation of the rectified \firstorder\ formula \nlbmath A 
from a \sententialtautology\
by applications of the 
generalized rules of simplification, and of\/
\math\gamma- and\/ 
\deltaminus\nolinebreak\hskip-.16em\nolinebreak-quantification 
\mbox{(and renaming of bound variables),}
\\then \math A has \propertyC\ of order\/ 
\nlbmaths{\displaystyle1+\sum_{i=1}^m\CARD{t_i}},
\par\noindent
where \math{t_1,\ldots,t_m} are the instances for the meta-variable \nlbmath t
of the generalized rule of\/ \mbox{\math\gamma-quantification} in its \math m
applications in the derivation of \nlbmaths A.
\getittotheright\qed\end{lemma}\vfill\pagebreak 

\begin{sloppypar}
\noindent
By the similarity of the expansion of the \deltaplusplusskolemizedform\
with the instantiation of the 
result of backward applications of the generalized rules of
\mbox{\math{\deltaplusplus\!}- and} 
restricted\/ \mbox{\math\gamma-quantification,} \hskip.1em
and by \lemmref{corollary extend signature}, \hskip.1em
we get:\end{sloppypar}

\begin{lemma}[From a Linear Derivation to \propertyCstar]
\label{lemma from yields to C star a la heijenoort}%
\\\noindent
Let \math n be a positive natural number. \hskip.3em
Let \math A be a rectified \firstorder\ formula. \hskip.3em
\\If there is a derivation of \nlbmath A
\begin{itemize}\noitem\item starting
from a quantifier-free formula
\nlbmath B with free \math\gamma-variables, such that
we \nolinebreak get a \sententialtautology\
by replacing each free \mbox{\math\gamma-variable} \nlbmath{\existsvari x{}}
occurring in \nlbmath B
with a term \nlbmath{t_{\existsvari x{}}} 
with \bigmaths
{\CARD{t_{\existsvari x{}}}\prec n},\noitem\item
proceeding
by applications of 
the generalized rules of simplification and of\/ 
\mbox{\math{\deltaplusplus\!}- and} restricted\/ \math\gamma-quantification 
(and renaming of bound variables),\noitem\end{itemize}
then \math A has \propertyCstar\ of order\/ \nlbmaths{n}.\qed\end{lemma}

\yestop\yestop\noindent
As a corollary of \lemmrefs
{lemma from C to yields a la heijenoort}
{lemma from yields to C a la heijenoort}, \hskip.2em
we can now state \herbrandsfundamentaltheorem\ in the form
presented by \heijenoortindex\heijenoort\ in \cite{heijenoort-herbrand}:

\yestop
\begin{theorem}[\herbrandsfundamentaltheorem\ \`a la \heijenoort]%
\label{theorem herbrand fundamental two}%
\par\noindent Let\/ \math A be a rectified \firstorder\ formula. 
\par\noindent The following two statements are logically equivalent. \hskip.3em
Moreover, 
we can construct a witness for each statement 
from a witness for the other one.\begin{enumerate}\item[1.]%
There is a positive natural number \nlbmath n such that
\math A has \propertyC\ of order \nlbmaths n.\item[2.]\sloppy 
There\,is\,a\,\sententialtautology\ \nlbmaths B, 
and 
\\there\,is\,a\,derivation\,of \math A from \nlbmath B
that\,consists\,in\,applications\,of the generalized rules of
\\simplification, \hskip.2em \mbox{\math\gamma-quantification,} \hskip.2em
 and\/ 
\deltaminus\nolinebreak\hspace*{-.16em}\nolinebreak-quantification
\\(and in the renaming of bound variables).
\getittotheright\qed\end{enumerate}\end{theorem}

\yestop\yestop\noindent As a corollary of \lemmrefs
{lemma from C star to yields a la heijenoort}
{lemma from yields to C star a la heijenoort}, \hskip.2em
we can now state \herbrandsfundamentaltheorem\ in an up-to-date form as follows.

\yestop
\begin{theorem}[\herbrandsfundamentaltheorem\ \`a la \wirth]%
\label{theorem herbrand fundamental three}%
\par\noindent Let\/ \math A be a rectified \firstorder\ formula. \hskip.3em
Let \math n be a positive natural number.
\par\noindent The following two statements are 
logically equivalent. Moreover, we can always construct witnesses for the second
statement if the first one holds.\begin{enumerate}\item[1.]%
\math A has \propertyCstar\ of order \nlbmaths n.\item[2.]%
There is a quantifier-free formula \math B 
with free \math\gamma-variables, \hskip.2em and
\\there is a substitution \math\sigma\ from these \math\gamma-variables
to terms \nlbmath t with \bigmaths{\CARD t\prec n}{}\\\getittotheright{such that
\math{B\sigma} \nolinebreak is a \sententialtautology, and} 
\\there is derivation of \math A from \nlbmath B
that consists in applications of the 
generalized rules of\linebreak simplification, \hskip.2em
\mbox{\math{\deltaplusplus\!}\hspace*{-.02em}-quantification}, \hskip.2em 
and restricted \math\gamma-quantification.
\getittotheright\qed\end{enumerate}\end{theorem}
\vfill\pagebreak
\section{Generalized Rules of Quantification in the Literature}\label
{section The Generalized Rules of Quantification in the Literature}
\begin{sloppypar}
Regarding \herbrandindex\herbrand's {\it modus ponens}-free calculus,
\heijenoortindex\heijenoort's invention is restricted to the {\em generalized}\/
rules of \math\gamma- and \mbox{\deltaminus\hskip-.1em-quantification.}
Let us briefly discuss their chronology.
\begin{itemize}%
\item\sloppy
In 1960, 
\schuettename\ \schuettelifetime\ 
published the generalized rules of 
\mbox{\math\gamma- and} \mbox{\deltaminus\hskip-.1em-quantification}
under the names ``S\,3'' (with a built-in application of the generalized rule
of \mbox{\math\gamma-simplification}) \hskip.2em
and ``S\,2'' 
in his first monograph on proof theory (in German)
\cite[\p\,78]{schuette60:_beweis}. \hskip.5em
\schuette's versions, however, come with an additional \mbox{restriction}
regarding the positions where quantifiers may be introduced. \
\schuette's top-down definition of 
positive and negative positions (``{\germanfont Positiv\-teile}'' 
and ``{\germanfont Negativ\-teile}'' \mbox{\cite[\p 11]{schuette60:_beweis}})
covers only \mbox{\math\alpha-operators} 
(according to \smullyan's classification).\footnote{%
 \Eg, 
 in the formula \nlbmath{A\tightund B}, the top position is
 positive, but the positions of the sub-formulas \nlbmath A and \nlbmath B
 are neither positive nor negative because a positive \nlbmath\tightund\ is a
 \mbox{\math\beta-operator}.} \
\schuette's rules restrict the introduction of 
quantifiers to positions that are either positive or negative
in this sense.%
\item
In 1968,
the analogs of the generalized rules of 
\mbox{\math\gamma- and \deltaminus\hskip-.1em-quantification}
without restrictions
for a tableau calculus are introduced 
in 
\heijenoortindex\heijenoort's unpublished \mbox{paper \cite{heijenoort-tree-herbrand}},
under the names 
 ``(V\math')''
 and 
 ``(IV\math')''.%
\item
In 1970, in \bernays' second\footnote{%
In the first edition \cite{grundlagen-first-edition-volume-two}, 
however,
these rules did not occur yet.%
} 
edition \cite[\p 166]{grundlagen-second-edition-volume-two} of \Vol\,II 
of 
the fundamental work on proof theory,
{\em Foundations of Mathematics} 
\makeaciteoftwo
{grundlagen-german-english-edition-volume-one-one} 
{grundlagen-german-english-edition-volume-one-two} 
by \hilbertname\ \hilbertlifetime\ and \bernaysname\ \bernayslifetime,
the generalized rules of 
\mbox{\math\gamma- and \deltaminus\hskip-.1em-quantification}
are introduced 
under the names 
\mbox{``\inpit{\mu^\ast}'' and ``\inpit{\nu^\ast}}\closequotecomma
but again in a restricted form. \
Let positive positions now be defined recursively over 
\math\alpha- and \math\beta-operators,
but without the mutually recursive definition of 
negative positions.\footnotemark\ \
\footnotetext{%
 \Eg, 
 in the formula \nlbmath{A\und \neg\neg B}, the top position and the position
 of the sub-formulas \nlbmath A and \nlbmath{\neg\neg B} are positive, 
 but the position of the sub-formulas \math{\neg B} \nolinebreak
 and \nlbmath B are neither positive nor negative.}%
\bernays' rules restrict the introduction of 
quantifiers to positive positions in this sense.
\item
In 1975, \heijenoortindex\heijenoort\ introduced the generalized rules of
\mbox{\math\gamma- and} \mbox{\deltaminus\hskip-.1em-quantification} 
(without restrictions and for a \hilbert\ calculus) \hskip.1em
in \nolinebreak the unpublished paper \cite{heijenoort-herbrand}, \hskip.1em
which we \nolinebreak quoted extensively in \nlbsectref{section calculus}.\item 
In 1977, 
in his second monograph on proof theory~\cite[\p\,20]{schuette-1977},
\hskip.3em
\schuette\ published versions of
the generalized rules of quantification
whose names and restrictions are in
principle the same as the before-mentioned ones of \cite{schuette60:_beweis},
but over a different set of sentential operators:
falsehood and implication instead of conjunction, disjunction, and negation.%
\item
In the 1980s, \heijenoortindex\heijenoort\ refers to the generalized rules of
\mbox{\math\gamma- and} \mbox{\deltaminus\hskip-.1em-quantification} in 
\makeaciteoftwo{heijenoort-oeuvre-herbrand}{heijenoort-work-herbrand},
but does not define them. \hskip.3em
See our \noteref{note on heijenoort-work-herbrand}.
\item In 2009, 
to the best of our knowledge,
the first correct definition of the generalized rules of 
\mbox{\math\gamma- and} \mbox{\deltaminus\hskip-.1em-quantification} and
\mbox{(\math\gamma-) simplification}
without restrictions
was published 
in our handbook article on \herbrand\ 
\makeaciteoftwo{herbrand-handbook}{SR--2009--01}.%
\end{itemize}%
\section{\herbrand's Original Rules}\label
{section herbrand's Original Rules}
To see that \heijenoortindex\heijenoort's inference rules are not just a slight variation
of \herbrandindex\herbrand's original rules, 
let us have a look at the latter ones.

\herbrandindex\herbrand\ had only the non-generalized versions of the rules of 
\mbox{\math\gamma- and} \deltaminus\hskip-.1em-quanti\-fi\-ca\-tion,
which 
(just as the non-generalized version of the rule of simplification)
result from our formalization of the generalized rules 
in \sectref{section calculus} 
by restricting \nlbmath{A[\ldots]} to the empty context
\hskip.1em
(\ie\nolinebreak\ \math{A[Q x\stopq H]}, \eg, is just \nlbmath{Q x\stopq H}). \
\herbrandindex\herbrand\ introduced them in \litsectref{2.2} of his \PhDthesis\
\cite{herbrand-PhD};
he named the rules of 
\mbox{\math\gamma- and \deltaminus\hskip-.16em-quantification}
``second'' and ``first rule of generalization'' 
\cite[\p\,74\f]{herbrand-logical-writings},
respectively 
(``{\frenchfont deuxi\`eme}''
 and 
 ``{\frenchfont premi\`ere r\`egle de g\'en\'eralisation}''
 \cite[\p\,68\f]{herbrand-ecrits-logiques}).

At the same places, 
we also find the (non-generalized) rule of simplification
(``{\frenchfont r\`egle de simplification}\/''), 
{\it modus ponens} (``{\frenchfont r\`egle d'implication}\/''),
and the ``rules of passage'' (``{\frenchfont r\`egles de passage}\/''). \
While \herbrandindex\herbrand\ does not need 
{\it modus ponens}\/ for completeness,
the deep inference rules of passage are needed 
(in anti-prenex direction)
because the shallow rules of quantification
--- contrary to the generalized \mbox{ones~---}
cannot introduce quantifiers at 
non-top positions.\end{sloppypar}
\begin{quote}{\bf Rules of Passage:} \
The following six \nolinebreak logical
  equivalences may be used for rewriting from left to right 
({\em prenex direction}\/) and from right to left 
({\em anti-prenex direction}\/), 
resulting in twelve \nolinebreak deep 
inference rules:
\par\noindent\LINEmath{\begin{array}{l c r@{~~~~}c@{~~~~}l}
       (1)
      & 
      &\neg\forall x\stopq A
      &\equivalent
      &\exists x\stopq \neg A
     \\(2)
      & 
      &\neg\exists x\stopq A
      &\equivalent
      &\forall x\stopq \neg A
     \\(3)
      & 
      &\inpit{\forall x\stopq A}\nottight\oder B
      &\equivalent
      &\forall x\stopq \inpit{A\hskip.06em\tightoder B}
     \\(4)
      & 
      &B\nottight\oder\forall x\stopq A
      &\equivalent
      &\forall x\stopq \inpit{B\hskip.06em\tightoder A}
     \\(5)
      & 
      &\inpit{\exists x\stopq A}\nottight\oder B
      &\equivalent
      &\exists x\stopq \inpit{A\hskip.06em\tightoder B}
     \\(6)
      & 
      &B\nottight\oder\exists x\stopq A
      &\equivalent
      &\exists x\stopq \inpit{B\hskip.06em\tightoder A}
      \\\end{array}}%
\par\noindent
  Here, \math B is a formula in which the variable \nlbmath x does not occur
  free.%
\end{quote} 
  As explained in \sectref{section calculus}, 
  if \nlbmath x \nolinebreak occurs \freely\ 
  in \nlbmath B, an implicit renaming of
  the bound occurrences of \nlbmath x in \nlbmath A
  is admitted to enable rewriting in
  prenex direction.

Note that \herbrand\ did not need rules of passage for conjunction
(besides the rules of passage for negation (1,\,2) and for disjunction 
(3,\,4,\,5,\,6)),
\hskip.2em
because he considered conjunction 
\bigmaths{A\hskip.09em\tightund\hskip.05em B}{} a
meta-level notion defined as \bigmaths
{\neg\inpit{\neg A\hskip.14em\tightoder\hskip.10em\neg B}}.

Finally, in \litsectref{5.6.A} of his \PhDThesis\ \cite{herbrand-PhD},
\herbrandindex\herbrand\ also introduces a ``generalized rule of simplification''
\cite[\p\,175]{herbrand-logical-writings} 
(``{\frenchfont r\`egle de simplification \mbox{g\'en\'eralis\'ee}}''
\cite[\p\,143]{herbrand-ecrits-logiques}). \
It is a deep inference rule for rewriting \bigmathnlb{H\oder H}{} to 
\nlbmaths H. \
Assuming that \math\neg\ and \tightoder\ are the only sentential operators
in our formula language, 
\herbrandindex\herbrand's rule 
it is tantamount to \heijenoortindex\heijenoort's
generalized rule of simplification
\mbox{presented} 
in \nlbsectref{section generalized rule of simplification}, \hskip.2em
because \herbrandindex\herbrand\ implicitly equates the variants \nlbmath H
and \nlbmath{H'}.
\vfill\pagebreak
\section{\herbrandsfalselemma\ and its Corrections}\label
{section herbrand's ``False Lemma'' and its Corrections}

For a given positive natural number \nlbmath n, \ 
{\em\herbrandsfalselemma}\/ \hskip.2em 
says that
\propertyC\ of order \nlbmath n \hskip.2em
is invariant under the application of the rules of passage.

\herbrandsfalselemma\ 
is wrong because the rules of passage may change the 
\outerskolemizedform. \ 
This happens whenever 
a \mbox{\math\gamma-quantifier} binding a variable \nlbmath x is moved
over a binary operator 
whose unchanged operand \nlbmath B 
contains a \mbox
{\math\delta-quantifier}.\footnote{\label{footnote same meta}%
 Here we use the same meta-variables
 as in our description of the rules of passage in 
 \nlbsectref{section herbrand's Original Rules}
 and assume that \math x \nolinebreak does not occur \freely\ in \nlbmath B.}

\begin{counterexample}\mbox{}\\
Let us again consider 
the rectified formula \nlbmath A of
\examref{example Property C one} of \sectref{section property C}, \hskip.3em
\ie\ the formula \bigmaths\formulaexampleC, which has the \outerskolemizedform\
\bigmaths\formulaexampleCouterskolemizedform, and thus 
\propertyC\ of order \nlbmath 3 
(by expansion \wrt\ \nlbmath{\{\bullet,\app{\forallvari b{}}\bullet\}}), 
\hskip.2em
but not \propertyC\ of order \nlbmaths 2. \hskip.6em
Let us move the \mbox{\math\gamma-quantifier} \nolinebreak
``\math{\exists a.}'' inward
by applying the last equivalence of the rules of passage in
anti-prenex direction. \hskip.2em
Then we obtain
\bigmaths{\inpit{\neg\exists b\stopq
{\app p b}}
\oder\exists a\stopq
{\app p a}}. \
The \mbox{(\deltaminus\hskip-.1em- \aswellas\ \math{\deltaplusplus\!}-)}
\skolemizedform\ 
of this formula is 
\bigmaths{\neg
{\app p{\forallvari b{}}}
\oder\exists a\stopq
{\app p a}}. \
Now, \mbox{expansion} \wrt\ \nlbmath{\{\forallvari b{}\}} 
results in a \sententialtautology. \
Thus, \hskip.2em
\math A \nolinebreak has \propertyC\ of order \nlbmaths 2.
\getittotheright\qed\end{counterexample}

\begin{sloppypar}
\par\halftop\noindent This, 
however, 
is not really 
a counterexample for \herbrandsfalselemma\
because \herbrandindex\herbrand\ treated our fresh constant 
\nolinebreak ``\nlbmath\bullet\nolinebreak\hskip.02em\nolinebreak''
from \defiref{definition champs finis}
as a variable and 
defined the height of a \skolem\ constant to be \nlbmath 1 
(just as we did in \defiref{definition champs finis}), \hskip.2em
but the height of a variable to be \nlbmath 0 (instead of \nlbmath 1), \hskip.2em
such that 
\ \mbox{\math{\CARD{{\forallvari b{}}{}}=1=
\CARD{\app{\forallvari b{}}{\bullet}}}}. \ 
As free variables and \skolem\ constants play exactly the same \role\
in the given context,
\herbrandindex\herbrand's definition of 
height is a bit counterintuitive and was possibly introduced
to avoid this counterexample. 
To \nolinebreak find a counterexample for
\herbrandsfalselemma\ also for \herbrandindex\herbrand's definition of height,
we just have to find a way to replace
\nolinebreak ``\nlbmath\bullet\nolinebreak\hskip.02em\nolinebreak''
with a \skolem\ constant, as shown in the following example.
\end{sloppypar}
\begin{counterexample}\label{example true counterexample}\mbox{}\\
Let us take \math B to be the 
rectified formula\\\LINEmaths{
\inparentheses{
\inpit{\neg\exists\boundvari b{}\stopq
{\app p{\boundvari b{}}}}
\oder
\exists\boundvari a{}\stopq
{\app q{\boundvari a{}}}}
\nottight{\oder}
\inparentheses{\inpit{\exists\boundvari x{}\stopq
{\app p{\boundvari x{}}}}
\und\neg\exists\boundvari y{}\stopq
{\app q{\boundvari y{}}}}}.
\\The \mbox{(\deltaminus\hskip-.1em- \aswellas\ \math{\deltaplusplus\!}-)}
  \skolemizedform\ of \nlbmath B is
\\\LINEmaths{
\inparentheses{
{\neg
{\app p{\forallvari b{}}}}
\oder
\exists\boundvari a{}\stopq
{\app q{\boundvari a{}}}}
\nottight{\oder}
\inparentheses{\inpit{\exists\boundvari x{}\stopq
{\app p{\boundvari x{}}}}
\und\neg
{\app q{\forallvari y{}}}}}.
\\Now, expansion \wrt\ \nlbmath{\{\forallvari b{},\forallvari y{}\}} 
results in a \sententialtautology. \
Thus, \math B \nolinebreak has \propertyC\ of order \nlbmaths 2. \hskip.4em
Let us move the \mbox{\math\gamma-quantifier} \nolinebreak
``\math{\exists a.}'' in formula \nlbmath B outward
by applying the last equivalence of the rules of passage in
prenex direction. \hskip.2em
Then we obtain 
\\\LINEmaths{\inparentheses{\exists\boundvari a{}\stopq
\inparentheses{
\inpit{\neg\exists\boundvari b{}\stopq
{\app p{\boundvari b{}}}}
\oder{\app q{\boundvari a{}}}}}
\nottight{\oder}
\inparentheses{\inpit{\exists\boundvari x{}\stopq
{\app p{\boundvari x{}}}}
\und\neg\exists\boundvari y{}\stopq
{\app q{\boundvari y{}}}}}.
\\Let us call this formula \nlbmaths{B'}. \hskip.5em
The \outerskolemizedformwithoutform\
(but not the \math{\deltaplusplus\!}-\skolem ized)\footnotemark\ \hskip.2em
form of \nlbmath{B'} is
\\[-1.7ex]\LINEmaths{\inparentheses{\exists\boundvari a{}\stopq
\inparentheses{
{\neg
{\app p{\app{\forallvari b{}}a}}}
\oder{\app q{\boundvari a{}}}}}
\nottight{\oder}
\inparentheses{\inpit{\exists\boundvari x{}\stopq
{\app p{\boundvari x{}}}}
\und\neg
{\app q{\forallvari y{}}}}}.
\\Now the smallest expansion that results in a sentential tautology is 
the one \wrt\ \nlbmaths{\{\forallvari y{},
\app{\forallvari b{}}{\forallvari y{}}\}}. \
Thus, the formula \nlbmath{B'} has \propertyC\ of order \nlbmaths 3,
but not of order \nlbmaths 2, no matter in which of the two ways we 
define the height of variables.
\getittotheright\qed\end{counterexample}
\pagebreak\par
\footnotetext{\label{note discussion inner}%
 From this example, one might get the idea that the flaw in 
 \herbrandsfalselemma\ would be a peculiarity of the 
 \outerskolemizedform.
 Indeed,
 for the \deltaplusplusskolemizedform\ 
 (\cfnlb\ \defiref{definition skolemized forms}), 
 moving \math\gamma-quantifiers
 with the rules of passage
 cannot change the number of arguments of the \skolem\ functions. 
 This does not help, however,
 to overcome the flaw in \herbrandsfalselemmacomma
 because, for the \deltaplusplusskolemizedform,
 moving \mbox{\math\delta-quantifiers} may change the number
 of arguments of the \skolem\ functions if the rule of passage is 
 applied within the scope of a \mbox{\math\gamma-quantifier} whose bound
 variable occurs in \nlbmath B but not in 
 \nlbmath A (\cfnlb\ \noteref{footnote same meta}). \ 
 The \deltaplusplusskolemizedform\ of 
 \bigmaths{
 {\exists\boundvari y 1\stopq
 \forall\boundvari z 1\stopq
 {\app p{\boundvari y 1,\boundvari z 1}}}
 \oder
 \exists\boundvari y 2\stopq
 \forall\boundvari z 2\stopq
 {\app q{\boundvari y 2,\boundvari z 2}}}{}
 is 
 \bigmathnlb{\exists\boundvari y 1\stopq
 {\app p{\boundvari y 1,\app{\forallvari z 1}{\boundvari y 1}}}
 \oder
 \exists\boundvari y 2\stopq
 {\app q{\boundvari y 2,\app{\forallvari z 2}{\boundvari y 2}}}}, 
 but the \deltaplusplusskolemizedform\ of any 
 prenex form
 has a {\em binary}\/ \skolem\ function, unless we use
 \henkin\ quantifiers as found in \hintikka's 
 \firstorder\ logic, \cfnlb\ \cite{hintikkaprinciples}, \cite{SR--2011--01}.}%
\begin{sloppypar}
The basic function of \herbrandsfalselemma\ in 
the proof of \herbrandsfundamentaltheorem\ 
is to establish the logical equivalence of 
\propertyC\ of a formula \nlbmath A
with \propertyC\ of the 
{prenex} and {anti-prenex forms} of \nlbmaths A.\end{sloppypar}

Let us see how this flaw of \herbrandsfalselemma\ has been
corrected.
\subsection{\bernays' Correction}
In\,1939, 
\bernays\ remarked that \herbrandindex\herbrand's proof is hard to
follow\footnote{%
 \majorheadroom
 In both the 1939 edition
 \cite[\litnoteref 1, \p 158]{grundlagen-first-edition-volume-two}
 and in the 1970 edition 
 \cite[\litnoteref 1, \p 161]{grundlagen-second-edition-volume-two},
 we read: \
 ``{\germanfontfootnote Die \herbrand sche Bewei\esi f\ue hrung ist
    schwer zu verfolgen}''} 
and 
---~for the first time~--- 
published a sound proof of a version of \herbrandsfundamentaltheorem.
This version is
restricted to prenex form, but more efficient in the number of terms
that have to be considered in a sub-expansion than \herbrandindex\herbrand's
quite global limitation to an expansion \wrt\ {\em all terms \nlbmath t with \
  \math{\CARD t\prec n}}, \ \ related to \propertyC\ of order 
\nlbmath n.\footnote{\label{footnote bernays}%
 \majorheadroom
 \Cfnlb\ \litsectref{3.3} of the 1939 edition 
 \cite{grundlagen-first-edition-volume-two}\@. \ In the 1970 edition
 \cite{grundlagen-second-edition-volume-two}, \hskip.2em
 \bernays\ also indicates how to remove the restriction to prenex formulas.}

\subsection{\goedel's and \dreben's Correction}
\begin{sloppypar}
According to a conversation with \heijenoortindex\heijenoort\ in
autumn\,1963,\footnote{%
 \majorheadroom
 \Cfnlb\ \cite[\p\,8, \litnoteref j]{herbrand-ecrits-logiques}\@.}  
\goedelname\ \goedellifetime\ noticed the 
flaw
in the proof
of \herbrandsfalselemma\ in 1943 and wrote a private note,
but did not publish \nolinebreak it. \
While \goedel's documented \mbox{attempts} to construct a counterexample to
\herbrandsfalselemma\ failed, \hskip.1em
he had actually worked out a {\em\repair}\/ of 
\herbrandsfalselemmacomma
which is sufficient for the proof of \herbrandsfundamentaltheorem.\footnote{%
 \majorheadroom
 \Cfnlb\ 
 \cite{goldfarb-herbrand-goedel}.%
}

In\,1962, when \goedelsrepair\ was still unknown,
a young student, 
\andrewsname\ \andrewslifetime, \hskip.1em
had the audacity to tell
his advisor \churchname\ \churchlifetime\ that there seemed to be a gap in 
the proof of \herbrandsfalselemmafullstop
\church\ sent \andrews\ to \drebenname\ \drebenlifetime, who finally came up
with a counterexample. 
And then \andrews\ constructed a simpler counterexample
(similar to the one we presented in \coexref{example true counterexample})
and joint work found a \repair\ similar to \goedel's.\footnote{%
 \majorheadroom
 \Cfnlb\ 
 \cite{andrews-herbrand-award},
 \cite{false-lemmas-in-herbrand},
 \cite{supplement-to-herbrand}.%
}\end{sloppypar}

\pagebreak

Roughly speaking, the \repaired\ lemma says that --- to keep \propertyC\ of
\nlbmath A invariant under (a single application of) a 
rule of passage --- we may have to step from 
order \nlbmath n \hskip.2em to order
\\[-.9ex]\noindent\LINEmaths{n\cdot\inpit{N^r+1}^n}.
\par\noindent Here \math{r} \nolinebreak is the
number of \mbox{\math\gamma-quantifiers} in whose scope the 
rule of passage is applied and \math N \nolinebreak is \nolinebreak 
the cardinality of \nlbmath{\termsofdepth n F} for the 
\outerskolemizedform\ \nlbmath F of \nlbmath A.\footnote{%
 \Cfnlb\ \cite[\p\,393]{supplement-to-herbrand}.}  

\goedelsanddrebensrepair\ is not particularly elegant 
because --- iterated several times until a 
prenex form is reached \nolinebreak--- it can lead to pretty high orders. \ 
Thus, although this 
\repair\ serves
well for a \mbox{finitistic proof} of \herbrandsfundamentaltheorem,
it results in a
complexity that is completely
unacceptable in practice (\eg\ intractable in automated reasoning), \hskip.2em
and this already for small non-prenex formulas.%

\subsection{\heijenoort's Correction}\label
{section Heijenoort's Correction}

\begin{sloppypar}
The lack of generalized versions of the rules of quantification force
\herbrandindex\herbrand's formal derivations 
in his {\it modus ponens}\/-free calculus
to take a detour over 
prenex form, which 
was standard at \herbrandindex\herbrand's time.
\mbox{For example},
\loewenheimindex\loewenheim\ and \skolem\ had always reduced their problems to 
prenex forms
of various kinds. 
The reduction of a proof task to prenex form
has several disadvantages, however, such as 
serious negative effects on proof complexity.\footnote{%
 \majorheadroom
 \Cfnlb\ \eg\ \cite{baazdelta}, \cite{baazleitschcolllog}.%
}

The 
surprising --- but, as a matter of fact, correct ---
thesis of \anellisname\ \anellislifetime\ in~\cite{anellis-loewenheim}
is that, building on the \loewenheimskolemtheorem,
it \nolinebreak was \herbrandindex\herbrand's work in elaborating 
\hilbert's concept
of ``being a proof'' that gave rise to the development
of the variety of \firstorder\ calculi in the 1930s,
such as the ones of the
\hilbert\ school
\cite{grundlagen-german-english-edition-volume-one-two}, \hskip.2em
and such as natural deduction and sequent calculi in~\cite{gentzen}. \
So \nolinebreak instead of fiddling around with \hilbertsepsilon\ 
and other existing calculi,
\herbrandindex\herbrand\ showed that the design of tailor-made calculi
suiting the given situation and its requirements
can be most fertile.\footnote{%
 \majorheadroom 
 This is a finding that still today is not sufficiently respected 
 in many areas, especially in those where logic is applied as a tool.
 \Cfnlb\ \eg\ \cite{wirth-hilbert-seki} for a discussion of this
 topic in the representation and computation of the semantics
 of natural language dialogs.%
}

As a consequence of this attitude of \herbrandindex\herbrand,
if he had become aware of the flaw in his 
\herbrandsfalselemmawithoutherbrandcomma
he would probably have avoided the whole detour over 
prenex forms.
And as there is not much of a choice for this avoidance,
he would probably have proceeded 
in the way of \heijenoortindex\citet{heijenoort-herbrand},
which we have presented \sectfromtoref
{section calculus}
{section The Generalized Rules of Quantification in the Literature}. \
Indeed, not only \heijenoortindex\heijenoort's generalized rules of quantification
are {\em deep}\/ inference rules,
already \herbrandindex\herbrand's rules of passage 
and of generalized simplification were
such rules which manipulate formulas at an arbitrary depth.

\pagebreak

In \herbrandindex\herbrand's proof of the equivalence 
of \propertyC\
and derivability in \herbrandindex\herbrand's {\em modus ponens}\/-free calculus
(stated in \herbrandsfundamentaltheorem),
there is only one application of 
\herbrandsfalselemmacomma
namely for dealing with applications of the rules of passage.
As the rules of passage are not part of \heijenoortindex\heijenoort's version of
\herbrandindex\herbrand's {\em modus ponens}\/-free calculus,
there is no need for \herbrandsfalselemma\ in 
{\em``\/\heijenoortsrepair''}\/ anymore.

Contrary to \bernaysrepair,
\heijenoortsrepair\ avoids 
the detour over the Extended 
First \nlbmath\varepsilon-Theorem of 
the proof of 
\bernays\ mentioned before; \cfnlb\ \noteref{footnote bernays}. \ 

Contrary to the \goedelsanddrebensrepair,
the corrections of \heijenoortsrepairindex\heijenoort\ and 
\bernaysrepairindex\bernays\ 
have the minor disadvantage\footnote
{Note that a corrected version of \herbrandsfalselemma\ 
 (in the sense of \goedelsanddrebensrepair) 
 is still needed if
 we want to prove that a derivation in \herbrandindex\herbrand's calculus 
 that includes 
 {\em modus ponens}\/ implies \propertyC\ in a finitistic sense. \
 As the example on top of page\,\,201 
 in \cite{herbrand-logical-writings} shows,
 however,
 an intractable increase of the order of \propertyC\ 
 cannot be avoided in general for 
 inference steps by \mbox{\em modus ponens}.}
that they 
do not provide a \repaired\
version of \herbrandsfalselemmafullstop

We see a unique quality of \heijenoortindex\heijenoort's correction in the following:
It \repairverb s exactly what went wrong in the execution \herbrandindex\herbrand's
proof plan. \
As we have shown in 
\sectrefs{section example}{section fundamental theorem}, \hskip.2em
\herbrandindex\herbrand's original proof of his \fundamentaltheorem\
can easily be adapted to a proof of 
its version \`a la \heijenoortindex\heijenoort\
(\ie\ our \theoref{theorem herbrand fundamental two})
(omitting the detour over prenex form). \ 
The idea to such a \repair\ of \herbrandindex\herbrand's proof
seems to have been published first in 
 \cite[\litnoteref{77}, \p\,555]{heijenoort-source-book} and
 \cite[\litnoteref{60}, \p171]{herbrand-logical-writings},
where, however, no generalized versions of the rules of 
quantification are mentioned.

Moreover,
the version of \herbrandsfundamentaltheorem\ \`a la \heijenoortindex\heijenoort\
resulting from \heijenoortsrepair\ 
is more elegant and concise than \herbrandindex\herbrand's \mbox{original one:}
\hskip.3em
It is not only more easily memorized by humans. \hskip.2em
Also the intractable and un\-intuitive\footnote{%
 \majorheadroom
 Unintuitive \eg\ in the sense of \cite{tait-2006}.%
}
rise in complexity
introduced by \goedelsanddrebensrepair\ is avoided. \hskip.2em
And because the relation to \propertyC\ is more straightforward
after \heijenoortsrepair,
the \fundamentaltheorem\ 
\`a 
la \heijenoortindex\heijenoort\ gives human beings
a better chance to construct a proof of a manageable size from 
a sub-expansion of the \skolemizedform\ 
as we \nolinebreak have done in \sectref{section example}.

Finally,
our new free-variable version of \herbrandsfundamentaltheorem\ 
(\ie\ \nolinebreak\theoref{theorem herbrand fundamental three})
adapts \heijenoortindex\heijenoort's correction to 
our novel free-variable calculus,
and 
makes the relation to modern reductive free-variable \mbox{calculi}\footnote{%
 \majorheadroom
 Such as free-variable tableau, sequent, and matrix calculi \cite
 {fitting}, \makeaciteoftwo{wirthcardinal}{SR--2011--01}.%
}
for computer-assisted theorem proving with strong automation support\footnote{%
 \majorheadroom
 \Cfnlb\ \cite{sergediss}, \cite{pds}, \cite{dodi-diss}.%
}
most obvious and clear.
\pagebreak\end{sloppypar}

\begin{figure}[h]
\begin{center}%
\includegraphics
[width=130mm,height=55mm,bb=0 0 459 238]%
{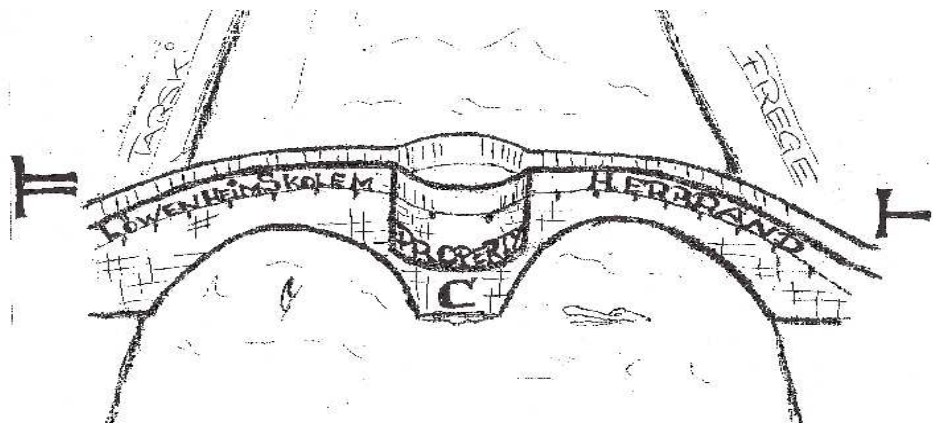}
\end{center}
\caption{\label{figure bridge}%
\mbox{}
The bridge of the \protect\loewenheimskolemtheorem\ and \protect\herbrandsfundamentaltheorem, based on the sentential \protect\propertyC\ standing firm in the river 
that divides the banks of valid and derivable formulas in the land of first-order predicate logic.
}\yestop\mbox{}
\end{figure}

\section
{Why \herbrandsfundamentaltheorem\ is Central,
\\Though Not Deep}%
\label{section Why is Central}

\begin{sloppypar}
Right at the beginning of \sectref{section calculus},
we have quoted from~\cite{heijenoort-herbrand},
where \heijenoortindex\heijenoort\ again quotes a statement of 
\citet{bernays-herbrand}, 
saying that \herbrandsfundamentaltheorem\
would be the central theorem of predicate logic. 
\heijenoortindex\heijenoort\ has quoted this statement also 
in~\cite[\p 1]{herbrand-ecrits-logiques} and 
in~\cite{heijenoort-spanish-lecture-course}.

Referring to the place of appearance of \herbrandsfundamentaltheorem\
in \cite[\litchapref 5]{herbrand-PhD}, \hskip.2em
\heijenoortindex\heijenoort\ also ends his article~\cite{heijenoort-logic-calculus-language}
as follows:\begin{quote}``Let me say simply, 
in conclusion,
that {\em\Begriffsschrift}~\cite{begriffsschrift},
\loewenheim's paper~\cite{loewenheim-1915},
and Chapter\,\,5 of \herbrand's thesis~\cite{herbrand-PhD}
are the three cornerstones of modern logic.''\end{quote}
We did not find any place, however,
where \heijenoortindex\heijenoort\ called \herbrandsfundamentaltheorem\
``deep'' 
(or
``one of the deepest results of logic''
as reported by \feferman\ in our quotation at the end of 
\sectref{section Introduction}).\footnote{%
 We also did not find a French publication of \heijenoortindex\heijenoort\
 where he calls \herbrandsfundamentaltheorem\ 
 ``profond''
 (``deep''). \
 This does not mean, of course, that such a publication does not exist,
 although we have carefully searched through all English and French texts of 
 \heijenoortindex\heijenoort\ listed in our references
 (some of them even by machine).}
Among the statements of \heijenoortindex\heijenoort\ 
we found on \herbrandsfundamentaltheorem,
closest to ``deep''
(but not really close) 
comes the following sentence in~\cite[\p\,247]{heijenoort-modern-logic}:
\begin{quote}``In his doctoral 
dissertation~\cite{herbrand-PhD}, \herbrandname\ \herbrandlifetime\
presented a theorem that revealed a profound feature of first-order logic.''
\end{quote}
\pagebreak

\noindent
On the one hand, 
it is obvious that \herbrandsfundamentaltheorem\
reveals a profound feature of first-order logic
and that it is one of the corner\-stones of modern logic.

On the other hand, however,
we cannot really justify why it should be called
``deep\closequotecomma
and agree with \feferman's statement in our quotation at the end of 
\sectref{section Introduction} that such an assessment would be
``debatable\closequotefullstopextraspace
In the vernacular of working mathematicians,
a ``deep theorem'' 
or ``deep result''
is a theorem whose proof is very hard to find,
especially
because it requires ingenuity,\footnote{
 A ``deep theorem\closequotecomma
 however, does not necessarily require computational power,
 such as in the case of the very hard proof of
 the \fourcolortheorem\
 \cite{four-color-theorem}.}
such as the proof of 
\fermatslasttheorem\
(\cfnlb\ \cite{wiles-fermat-correction}, \cite{wiles-fermat}). \hskip.4em
Although the proof of \herbrandsfundamentaltheorem\ in\,1929 
was more ``avantgarde''
than most other proofs in the history of modern logic 
at their respective times,
and although the theorem itself is ingenious,
the proof of \herbrandsfundamentaltheorem\ is {\em not}\/ hard to find. \
As we have seen in 
\sectref{section example}, \hskip.2em
the proof of \herbrandsfundamentaltheorem\
basically verifies the equivalence of 
\PropertyC\ with
\herbrandindex\herbrand's
{\it modus ponens}\/-free calculus,
which was up to his own design.\footnote{%
 \majorheadroom
 In terms of \figuref{figure bridge}, 
 we could say that 
 the place where \herbrandindex\herbrand's part of the bridge 
 entered the bank of 
 derivable formulas in
 the land of first-order predicate logic
 was up to \herbrandindex\herbrand's choice.} \
\herbrandindex\herbrand\ was a creative genius in designing all the novel
features he needed for that proof at his time,
but we do not think that all this classifies that proof as being 
deep in the sense we described here.

Compared to the great achievements of \herbrandindex\herbrand\ in his \PhDthesis,
the gap in his original proof is negligible,
at least \wrt\ \heijenoortindex\heijenoort's correction as discussed 
\mbox{in \sectref{section Heijenoort's Correction}}. \
Moreover,
there can be no doubt that he would have been able to patch the
gap if he had only become aware of it. \hskip.1em
Therefore, 
we think that this gap does not justify to say that 
\herbrandindex\herbrand\ failed to prove his \fundamentaltheorem,
nor that the theorem is a deep one.

All in all, 
\heijenoortindex\heijenoort\ was definitely right in calling 
\herbrandindex\herbrand\ a {\frenchfont\em g\'enie cr\'ea\-teur}
\mbox{\cite[\p 1]{herbrand-ecrits-logiques}}, \hskip.1em
whereas it remains an open question here whether he really called 
\herbrandsfundamentaltheorem\ ``deep\closequotefullstop

To find out, \hskip.1em
on the other hand, \hskip.1em
why \herbrandsfundamentaltheorem\ should be called 
``central\closequotecomma
let us now have a look at a translation of 
\bernays' briefly quoted statement in its context:
\begin{quote}
``In its proof-theoretic form,  
\herbrandindex\herbrand's Theorem can be seen
as the central theorem of predicate logic.
It expresses the relation of predicate logic to propositional logic
in a concise and felicitous form.''\footnote{%
 \majorheadroom
 The German original text is:
 ``Der \herbrand'sche Satz in seiner bewei\esi theoretischen Fassung 
 kann als das zentrale Theorem der Pr\ae dikaten\-logik angesehen werden.
 In ihm wird die Beziehung der Pr\ae dikatenlogik zur Au\esi sagen\-logik 
 auf eine pr\ae gnante Form gebracht.''%
}\end{quote}
So it is quite obvious that \bernays\ did not want to say that 
\herbrandsfundamentaltheorem\ 
is situated in the land of first-order predicate logic
like Paris in France;\footnote{%
 \majorheadroom
 You can hardly avoid Paris when traveling in France.%
} \hskip.3em
he \nolinebreak just wanted to say that it establishes an 
(in the limit complete) 
reduction of the problem of derivability in one calculus to the problem of
derivability 
in a simpler one 
(a \nolinebreak{\em\herbrand\ reduction} as it is called today), \hskip.1em
and that this reduction is ``concise and felicitous\closequotefullstopnospace
\footnote{%
 \majorheadroom
 The original German word ``{\germanfontfootnote pr\ae gnant}'' 
 (which we translate as ``concise and felicitous'') 
 has the same etymological root as
 the English word ``pregnant\closequotecomma
 but it has mostly lost that direct meaning.%
}

Let us elaborate a bit on the questions
why \herbrandsfundamentaltheorem\ is central
and why the form of its expression is concise and felicitous.

As an explicit notion, 
\herbrand's \propertyC\
was first formulated in \herbrand's thesis~\cite{herbrand-PhD}. 
Implicitly it was used already before, namely
by 
\loewenheimname\ \loewenheimlifetime\ in~\cite{loewenheim-1915} 
and  by 
\skolemname\ \skolemlifetime\ in~\cite{skolem-1928}. \hskip.3em
It \nolinebreak is the main property of
\herbrandindex\herbrand's work and may well be called the central property of
\firstorder\ logic, for reasons to be explained in the following.

Contrary to the 
the unjustified\footnote
{\Cf\ \cite[\litsectref{3.12}, \esp\,\litnoteref{98}]{herbrand-handbook} 
 or \cite[\litsectref{14}, \esp\,\litnoteref{114}]{SR--2009--01}} 
criticism of \skolem, \herbrandindex\herbrand, and
\heijenoortindex\heijenoort, 
there are no essential gaps in the proof of 
\loewenheim\ 
in~\cite{loewenheim-1915} 
of the later so-called \loewenheimskolemtheorem,
which says
that \propertyC\ is equivalent to model-theoretic
first-order validity. \
In his \PhDthesis,
\herbrandindex\herbrand\ also showed the equivalence of 
his calculi with those of the \hilbert\ school
\cite{grundlagen-german-english-edition-volume-one-two} \hskip.1em
and the \PM~\cite{PM}. \
Therefore,
as a consequence of the \loewenheimskolemtheorem,
the completeness of all these calculi is an immediate
corollary of \herbrandsfundamentaltheorem. \
But \herbrandindex\herbrand\ did not trust the left arc of the bridge in our illustration
of \figuref{figure bridge}. \
\herbrandindex\herbrand's standpoint was so radically finitistic that 
in the area of logic he
did not accept model theory or set theory at \nolinebreak all. \
And so \goedel\ proved the completeness of first-order logic first 
when he
submitted his thesis~\cite{goedel-completeness} in 1929, 
in the same year as \herbrandindex\herbrand, and the
theorem is now called {\em \goedel's Completeness Theorem}\/ 
in all textbooks on logic.\end{sloppypar}

So,
although \propertyC\ is purely sentential and does not really 
belong to the land of first-order predicate logic with its
semantical and syntactical banks (\cfnlb\ \figuref{figure bridge}), \hskip.25em
it can be seen as the central property of 
first-order predicate logic 
and \herbrandsfundamentaltheorem\ connects this 
property with derivability in the standard calculi. \hskip.3em
Moreover, the way this connection is established,
is concise and felicitous because the formal proof in 
\herbrandindex\herbrand's {\it modus ponens}\/-free calculus
is linear, has the ``sub-formula'' property, 
and does not take any detours because it does not have rules
that eliminate intermediate results, 
such as {\em modus ponens}\/ and \gentzen's cut rule do. \hskip.3em
Furthermore, the proof construction of \herbrandsfundamentaltheorem\
out of \propertyC\
gets even more concise and felicitous in \heijenoortsrepair;
\cfnlb\ \sectref{section Heijenoort's Correction}.

\begin{sloppypar}
Moreover,
\PropertyC\ has strongly influenced the early history of 
automated deduction in the \nth 2\,half of the \nth{20}\,century
\makeaciteoftwo{herbrand-handbook}{SR--2009--01}. \
The different treatment of 
\mbox{\math\delta-quantifiers} and \math\gamma-quantifiers
in \propertyC,
namely by \skolemization\ and expansion, respectively,
as found in~\cite{herbrand-PhD}, \cite{skolem-1928}, \hskip.1em
rendered the reduction to 
sentential logic by hand for small examples (and later with a computer)
practically executable for the first time.\footnote{%
 \majorheadroom
 For instance, the elimination of both \math\gamma- and
 \math\delta-quantifiers with the help of
 \hilbert's \mbox{\nlbmath\varepsilon-operator} suffers from an exponential
 complexity in formula size. \ 
 As \nolinebreak a result, already small formulas grow so large 
 that the mere size makes them inaccessible to human inspection; \hskip.3em
 and this is still the case
 for the term-sharing representation of 
 \mbox{\math\varepsilon-terms} of~\cite[\litexamref 8]{wirth-jal},
 which is further elaborated in~\cite[\litexamref{3.7}]{SR--2011--01}.%
} \
This different treatment of the two kinds of quantification
is inherited from the 
\peirce--\schroeder\ 
tradition\footnote{\label{note peirce schroeder tradition}%
 \majorheadroom
 For the heritage of \peirce\ see
 \cite{brady}, \cite{peirce-1885}; \hskip.2em
 for that of 
 \schroeder\ see \cite{schroeder-handbook}, \cite{schroeder-vorlesungen-III}.%
}
which came on \herbrandindex\herbrand\ via 
\loewenheimindex\loewenheim\ and \skolem. \ 
\russell\ and \hilbert\ had already merged
that tradition with the one of \frege,
sometimes emphasizing their \frege\ heritage over the
one of \peirce\ and \schroeder.\footnote
{While this emphasis on \frege\ 
 will be understood by everybody who ever had the 
 fascinating experience of reading \frege,
 it put some unjustified bias to the historiography of modern logic,
 still present in the selection of \heijenoortindex\heijenoort's 
 famous source book \cite{heijenoort-source-book}; \
 \cfnlb\ \eg\ 
 \cite[\litchapref 3]{anellis-heijenoort-long}.%
}
It \nolinebreak was \herbrandindex\herbrand\ who completed the bridge 
between these two traditions with his \fundamentaltheorem,
as depicted in 
\figuref{figure bridge}.

So one could well say with some justification that, 
as a consequence of \propertyC\
being the central property of first-order predicate logic,
the \loewenheimskolemtheorem\ and \herbrandsfundamentaltheorem\
are the central theorems of first-order predicate logic.
\end{sloppypar}

Moreover,
\herbrandindex\herbrand's {\it modus ponens}\/-free calculus
--- at least after \heijenoortsrepair\ and our adaption to
    free-variable calculi ---
is most close to modern approaches in automated theorem proving
that aim at a synergetic
combination of the semantical strength of human mathematicians
with the syntactical and computational strength of computing machines,
as we have indicated already at the end of \sectref
{section Heijenoort's Correction}. \
The manner in which automated theorem-proving systems based on 
modern sequent, tableau, and matrix calculi
organize proof search\footnote{%
 \majorheadroom
 \Cfnlb\ \eg\ \cite{sergecore}, \cite{wallen}, \cite{wirthcardinal}.}
does not follow the \hilbert\ school and their 
\math\varepsilon-elimination theorems, but \gentzen's and \herbrandindex\herbrand's calculi. \ 
Moreover, regarding their \skolemization, their deep inference,\footnote{%
 \majorheadroom
 Note that although the deep inference rules of 
 {\em generalized}\/ quantification are an 
 extension of \herbrandindex\herbrand's calculi by \heijenoortindex\heijenoort,
 the deep inference rules of passage and of 
 generalized  simplification are \herbrandindex\herbrand's original contributions.}
and their focus on \mbox{\math\gamma-quantifiers} and their multiplicity,
these modern 
proof-search calculi
are even more in \herbrandindex\herbrand's tradition than in 
\gentzen's.
\vfill\cleardoublepage
\section{Conclusion and Future Work}\label
{section Conclusion}

\begin{sloppypar}
We have most clearly described and discussed 
\heijenoortindex\heijenoort's version of 
\herbrandindex\herbrand's {\it modus ponens}\/-free calculus
and its relation to \herbrandsfundamentaltheorem\
in \sectfromtoref{section calculus}{section fundamental theorem}, \hskip.2em
where our presentation had a didactical focus because
(as stated already in \sectref{section Motivation}) \hskip.1em
we hope that our improvements of this calculus will become part of the standard
education of logicians, 
just as well as the famous construction of Cut-free
proofs according to \gentzensHauptsatz~\cite{gentzen}.

Already from our brief discussion in \sectref
{section The Generalized Rules of Quantification in the Literature}
on the publication history of the {\em generalized}\/ rules of quantification
alone, it becomes clear that 
\heijenoortindex\heijenoort's unpublished handout~\cite{heijenoort-herbrand} would have
deserved publication in an international journal, provided that \heijenoortindex\heijenoort\
had worked out the proof properly, involving our correction 
of \heijenoortindex\heijenoort's generalized rule of simplification.\footnote
{That this proof would have been within
\heijenoortindex\heijenoort's reach is out of question; \hskip.3em
although he had not worked it out, 
as we can see from
the bug in his generalized rule of simplification.}

After clarifying \herbrand's original rules in 
\sectref{section herbrand's Original Rules}, \hskip.2em
it became possible to discuss \herbrandsfalselemma\
in \sectref{section herbrand's ``False Lemma'' and its Corrections}: \hskip.3em
After briefly describing \bernays', \goedel's, and \dreben's corrections,
we explained why we think that ``\heijenoortindex\heijenoort's correction''
\mbox{(consisting} in the material we presented in \nolinebreak
 \sectfromtoref{section calculus}{section fundamental theorem}) \hskip.1em
is actually the one that fits 
\herbrandindex\herbrand's style in his work on logic 
best and offers the most tractable solution, 
which also has a close relation to automated and human-oriented theorem proving.

Finally,
in \sectref{section Why is Central},
we have checked in which context
\bernays\ called \herbrandsfundamentaltheorem\ ``central''
and 
--- to defend \heijenoortindex\heijenoort's and \bernays'
    assessment of \herbrandsfundamentaltheorem\ 
    against \feferman's critique quoted at the end 
    \mbox{of \sectref{section Introduction} ---}
elaborated on the question why its form of expression is 
``concise and felicitous\closequotecomma
though not deep.

\halftop\halftop\noindent
The strong dependencies\begin{enumerate}\notop\item
between the \outerskolemizedform\ 
and the proof of \herbrandsfundamentaltheorem\
(actually \lemmref{lemma from C to yields a la heijenoort})
described in \sectref{section formal proof}, \hskip.1em
and\noitem\item between the \deltaplusplusskolemizedform\ and the 
liberalized \math\delta-rule in the form of our
generalized rule of \deltaplusplus-quantification
(which forced us to introduce \propertyCstar\ 
 in \defiref{definition properties C and C star}),\notop\end{enumerate}
became clear only in this article and
may deserve further consideration.

\halftop\halftop\noindent
A careful bilingual edition of \herbrandindex\herbrand's complete works
on the basis of the previous editorial achievements 
\makeaciteoftwo{herbrand-1936}{herbrand-logical-writings}, 
which was on the agenda until 
\heijenoortindex\heijenoort\ died~\cite[\p\,383]{heijenoort-work}, \hskip.2em
is still in high demand.

Last but not least,
we would like to suggest our 
handbook article \makeaciteoftwo{herbrand-handbook}{SR--2009--01} 
for further reading on \herbrandindex\herbrand's
work in logic, 
also because
a comparison with the material we have further elaborated in this
article will make some additional
aspects of \herbrandsfundamentaltheorem\ clear.\footnotemark
\vfill\pagebreak\par\footnotetext{%
For example, as indicated already in our quotation from 
\heijenoortindex\heijenoort's~\cite{heijenoort-herbrand} at the beginning of 
\sectref{section Introduction},
besides the expansions used in this article
(nowadays called ``\herbrand\ expansions''), \hskip.2em
there is the alternative of \herbrand\ disjunctions
used in \makeaciteoftwo{herbrand-handbook}{SR--2009--01}. \
Considering the intractable size 
of a \herbrand\ disjunction compared to a \herbrand\ expansion
(\cfnlb\ \noteref{note disjunction versus expansion}), \hskip.1em
we have come to the conclusion that expansions should be preferred
in general, not only because they are \herbrand's own choice,
but also because they are more intuitive and more tractable,
especially in combination with the ``sub-expansions'' 
we \nolinebreak have introduced in 
\defiref{definition sub-expansion}.
 This means, for instance, that presentations of \herbrandsfundamentaltheorem\
 via \herbrand\ disjunctions, 
 such as in \cite[\PP{247}{249}]{heijenoort-modern-logic},
 should be avoided because they are none of the following:
 historically adequate, tractable, intuitive.}

\mbox{}


\noindent{\bf Acknowledgments:}
I would like to thank 
\autexiername, 
\fefermanname, \goldfarbname,
and especially \anellisname\
for their kind help via \eMAIL\ with some of my
open questions.
\end{sloppypar}
\vfill\cleardoublepage
\halftop
\addcontentsline{toc}{section}{References}
\small
\nocite{grundlagen-first-edition-volume-one,grundlagen-second-edition-volume-one,grundlagen-first-edition-volume-two,grundlagen-second-edition-volume-two,grundlagen-german-english-edition-volume-one-one}
\bibliography{herbrandbib}

\begin{thebibliography}{}

\bibitem[\protect\citeauthoryear{Abeles}{1994}]{abeles-herbrand-unification}
Francine Abeles.
\newblock {\herbrand's Fundamental Theorem} and the beginning of logic
  programming.
\newblock {\em Modern Logic}, 4:63--73, 1994.

\bibitem[\protect\citeauthoryear{Andrews{\protectedandrewsindex}}{2003}]{andre%
ws-herbrand-award}
Peter~B. Andrews{\protectedandrewsindex}.
\newblock {\herbrand\ Award} acceptance speech.
\newblock {\em J. Automated Reasoning}, 31:169--187, 2003.

\bibitem[\protect\citeauthoryear{Anellis{\protectedanellisindex}}{1991}]{anell%
is-loewenheim}
Irving~H. Anellis{\protectedanellisindex}.
\newblock The {\loewenheimskolemtheorem}, theories of quantification, and proof
  theory.
\newblock 1991.
\newblock In \cite[\PP{71}{83}]{drucker}.

\bibitem[\protect\citeauthoryear{Anellis{\protectedanellisindex}}{1992}]{anell%
is-heijenoort-long}
Irving~H.{\protectedheijenoortindex} Anellis{\protectedanellisindex}.
\newblock {\em Logic and Its History in the Work and Writings of\/
  {\heijenoortname}}.
\newblock Modern Logic Publ., Ames (IA), 1992.

\bibitem[\protect\citeauthoryear{Autexier{\protectedautexierindex}
  \bgroup\&al.\egroup }{2006}]{pds}
Serge Autexier{\protectedautexierindex}, {\mbox{Ch}}ristoph
  Benz{\-}m{\ue}ller{\protectedbenzmuellerindex}, Dominik Dietrich, Andreas
  Meier, and Claus-Peter Wirth{\protectedwirthindex}.
\newblock A generic modular data structure for proof attempts alternating on
  ideas and granularity.
\newblock 2006.
\newblock In \cite[\PP{126}{142}]{forthMKMfive}, \www\url{/p/pds}.

\bibitem[\protect\citeauthoryear{Autexier{\protectedautexierindex}}{2003}]{ser%
gediss}
Serge Autexier{\protectedautexierindex}.
\newblock {\em Hierarchical Contextual Reasoning}.
\newblock PhD thesis, \addressuniSBshort, 2003.

\bibitem[\protect\citeauthoryear{Autexier{\protectedautexierindex}}{2005}]{ser%
gecore}
Serge Autexier{\protectedautexierindex}.
\newblock The {\sc core} calculus.
\newblock 2005.
\newblock In \cite[\PP{84}{98}]{twentiethCADEfive}.

\bibitem[\protect\citeauthoryear{Baaz \bgroup\&\ \egroup
  Ferm{\"u}ller}{1995}]{baazdelta}
Matthias Baaz and {\mbox{Ch}}ristian~G. Ferm{\"u}ller.
\newblock Non-elementary speedups between different versions of tableaux.
\newblock 1995.
\newblock In \cite[\PP{217}{230}]{fourthTABLEAUninetyfive}.

\bibitem[\protect\citeauthoryear{Baaz \bgroup\&\ \egroup
  Leitsch}{1995}]{baazleitschcolllog}
Matthias Baaz and Alexander Leitsch.
\newblock Methods of functional extension.
\newblock {\em Collegium Logicum --- Annals of the {\goedelname\ Society}},
  1:87--122, 1995.

\bibitem[\protect\citeauthoryear{Baumgartner \bgroup\&al.\egroup
  }{1995}]{fourthTABLEAUninetyfive}
Peter Baumgartner, Reiner H{\"a}hnle, and Joachim Posegga, editors.
\newblock {\em {\thefourthTABLEAUninetyfive}}, number 918 in Lecture Notes in
  Artificial Intelligence. {\springerverlag}, 1995.

\bibitem[\protect\citeauthoryear{Beckert \bgroup\&al.\egroup
  }{1993}]{deltaplusplus}
Bernhard Beckert, Reiner H{\"a}hnle, and Peter~H. Schmitt.
\newblock The even more liberalized {\math\delta}-rule in free-variable
  semantic tableaus.
\newblock 1993.
\newblock In \cite[\PP{108}{119}]{goedelcolloquium1993}.

\bibitem[\protect\citeauthoryear{Berka \bgroup\&\ \egroup
  Kreiser}{1973}]{logiktexte}
Karel Berka and Lothar Kreiser, editors.
\newblock {\em {Logik-Texte -- Kommentierte Au\esi wahl zur Geschichte der
  modernen Logik}}.
\newblock \akademieverlag, 1973.
\newblock \nth 2\,\rev\,\edn\ (\nth 1\,\edn\,1971; \nth
  4\,\rev\,\rev\,\edn\,1986).

\bibitem[\protect\citeauthoryear{Bernays\protect\bernaysindex}{1957}]{bernays-%
herbrand}
Paul Bernays\protect\bernaysindex.
\newblock {\Ue ber den Zusammenhang des \herbrand schen Satzes mit den neueren
  Ergebnissen von \schuette{\protectedschuetteindex} und {\sc Stenius}}.
\newblock In {\em Proceedings of the International Congress of Mathematicians
  1954}, Groningen and Amsterdam, 1957. Noordhoff and \northholland.

\bibitem[\protect\citeauthoryear{Brady}{2000}]{brady}
Geraldine Brady.
\newblock {\em From {\protectedpeirceindex\peirce} to {\skolem}: A Neglected
  Chapter in the History of Logic}.
\newblock {\northholland}, 2000.

\bibitem[\protect\citeauthoryear{Cohen \bgroup\&\ \egroup
  Wartofsky}{1967}]{boston-studies-3}
Robert~S. Cohen and Marx~W. Wartofsky, editors.
\newblock {\em {\Proc\ of the Boston Colloquium for the Philosophy of Science,
  1964--1966}: In Memory of {\namefont Norwood Russell Hanson}}.
\newblock Number~3 in Boston Studies in the Philosophy of Science.
  \DReidelpublishing, 1967.

\bibitem[\protect\citeauthoryear{Dietrich}{2011}]{dodi-diss}
Dominik Dietrich.
\newblock {\em Assertion Level Proof Planning with Compiled Strategies}.
\newblock Optimus Verlag, Alexander Mostafa, \Goettingen, 2011.
\newblock \PhDthesis, \Dept\ Informatics, \addressuniSBshort.

\bibitem[\protect\citeauthoryear{Dreben{\protecteddrebenindex} \bgroup\&\
  \egroup Denton}{1963}]{supplement-to-herbrand}
Burton Dreben{\protecteddrebenindex} and John Denton.
\newblock A supplement to {\herbrand}.
\newblock {\em {\jslname}}, 31:393--398, 1963.

\bibitem[\protect\citeauthoryear{Dreben{\protecteddrebenindex}
  \bgroup\&al.\egroup }{1963}]{false-lemmas-in-herbrand}
Burton Dreben{\protecteddrebenindex}, Peter~B. Andrews{\protectedandrewsindex},
  and St{{\aa}}l Aanderaa.
\newblock False lemmas in {\herbrand}.
\newblock {\em {\bullamsname}}, 69:699--706, 1963.

\bibitem[\protect\citeauthoryear{Drucker}{1991}]{drucker}
Thomas Drucker, editor.
\newblock {\em Perspectives on the History of Mathematical Logic}.
\newblock {\birkhaeuser}, 1991.

\bibitem[\protect\citeauthoryear{Feferman}{1993a}]{heijenoort-life}
Anita~Burdman Feferman.
\newblock {\em Politics, Logic and Love --- The Life {of\/ \heijenoortname}}.
\newblock A K Peters, Wellesley (MA), 1993.

\bibitem[\protect\citeauthoryear{Feferma{\protectedfefermanindex}n}{1993b}]{he%
ijenoort-work}
Sol(omon) Feferma{\protectedfefermanindex}n.
\newblock {\heijenoortname's} scholarly work.
\newblock 1993.
\newblock In \cite[\PP{371}{390}]{heijenoort-life}.

\bibitem[\protect\citeauthoryear{Fitting}{1990}]{Fitting90}
Melvin Fitting.
\newblock {\em First-order logic and automated theorem proving}.
\newblock {\springerverlag}, 1990.
\newblock \nth 1\,\edn\ (\nth 2\,\rev\,\edn\ is \cite{fitting}).

\bibitem[\protect\citeauthoryear{Fitting}{1996}]{fitting}
Melvin Fitting.
\newblock {\em First-order logic and automated theorem proving}.
\newblock {\springerverlag}, 1996.
\newblock \nth 2\,\rev\,\edn\ (\nth 1\,\edn\ is \cite{Fitting90}).

\bibitem[\protect\citeauthoryear{Frege\protectedfregeindex}{1879}]{begriffssch%
rift}
Gottlob Frege\protectedfregeindex.
\newblock {\em {\Begriffsschrift, eine der arithmetischen nachgebildete
  Formelsprache de\es\ reinen Denken\es}}.
\newblock {Verlag von L. Nebert, \Halle}, 1879.
\newblock Corrected facsimile in \cite{begriffsschrift-und-andere}. \ Reprint
  of \PP{III}{VIII} and \PP{1}{54} in \cite[\PP{48}{106}]{logiktexte}. \
  {English} translation in \cite[\PP{1}{82}]{heijenoort-source-book}.

\bibitem[\protect\citeauthoryear{Frege\protectedfregeindex}{1964}]{begriffssch%
rift-und-andere}
Gottlob Frege\protectedfregeindex.
\newblock {\em {\Begriffsschrift\ und andere Auf\-s\ae tze}}.
\newblock \wissenschaftlichebuchgesellschaftdarmstadt, 1964.
\newblock {Zweite Au\fli age, mit \husserlname s und \scholzname' Anmerkungen,
  herau\esi gegeben von {\namefont Ignacio Angelelli}}.

\bibitem[\protect\citeauthoryear{Gabbay{\protectedgabbayindex} \bgroup\&\
  \egroup Woods}{2004\ff}]{handbook-of-the-history-of-logic}
Dov Gabbay{\protectedgabbayindex} and John Woods, editors.
\newblock {\em Handbook of the History of Logic}.
\newblock \northholland, 2004\ff.

\bibitem[\protect\citeauthoryear{Gentzen}{1935}]{gentzen}
Gerhard Gentzen.
\newblock {Untersuchungen \ue ber das logische Schlie\sz en}.
\newblock {\em Mathematische Zeitschrift}, 39:176--210,405--431, 1935.
\newblock Also in \cite[\PP{192}{253}]{logiktexte}. {English} translation in
  \cite{gentzen-collected}.

\bibitem[\protect\citeauthoryear{Gentzen}{1969}]{gentzen-collected}
Gerhard Gentzen.
\newblock {\em The Collected Papers of \gentzenname}.
\newblock {\northholland}, 1969.
\newblock \Ed\ by \szaboname.

\bibitem[\protect\citeauthoryear{G{\oe
  de{\protectedgoedelindex}l}}{1930}]{goedel-completeness}
Kurt G{\oe de{\protectedgoedelindex}l}.
\newblock {Die Vollst\ae ndigkeit der Axiome des logischen Funktionen\-kalk\ue
  ls}.
\newblock {\em {\monatsheftempname}}, 37:349--360, 1930.
\newblock With {English} translation also in \cite[\Vol\,I,
  \PP{102}{123}]{goedelcollected}.

\bibitem[\protect\citeauthoryear{G{\oe
  de{\protectedgoedelindex}l}}{1986ff.}]{goedelcollected}
Kurt G{\oe de{\protectedgoedelindex}l}.
\newblock {\em Collected Works}.
\newblock {\oxfordunipress}, 1986ff.
\newblock \Ed\ by \fefermanname, \dawsonname\protect\dawsonindex,
  \goldfarbname, \heijenoortname, \kleenename, \parsonsname, \siegname,
  \etalabbrev.

\bibitem[\protect\citeauthoryear{Goldfarb{\protectedgoldfarbindex}}{1970}]{rev%
iew-herbrand-ecrits-logiques}
Warren Goldfarb{\protectedgoldfarbindex}.
\newblock Review of \cite{herbrand-ecrits-logiques}.
\newblock {\em The Philosophical Review}, 79:576--578, 1970.

\bibitem[\protect\citeauthoryear{Goldfarb{\protectedgoldfarbindex}}{1993}]{gol%
dfarb-herbrand-goedel}
Warren Goldfarb{\protectedgoldfarbindex}.
\newblock {\herbrand's} error and {\goedel's} correction.
\newblock {\em Modern Logic}, 3:103--118, 1993.

\bibitem[\protect\citeauthoryear{Gonthier}{2008}]{four-color-theorem}
Georges Gonthier.
\newblock Formal proof --- the {\fourcolortheorem}.
\newblock {\em {\noticesamsname}}, 55:1382--1393, 2008.

\bibitem[\protect\citeauthoryear{Gottlob \bgroup\&al.\egroup
  }{1993}]{goedelcolloquium1993}
Georg Gottlob, Alexander Leitsch, and Daniele Mundici, editors.
\newblock {\em {Computational Logic and Proof Theory, \Proc\ \nth 3
  \goedelname\ Colloquium}}, number 713 in Lecture Notes in Computer Science.
  {\springerverlag}, 1993.

\bibitem[\protect\citeauthoryear{Heijenoort{\protectedheijenoortindex}}{1967}]%
{heijenoort-logic-calculus-language}
Jean~van Heijenoort{\protectedheijenoortindex}.
\newblock Logic as a calculus and logic as a language.
\newblock {\em Synthese}, 17:324--330, 1967.
\newblock Also in \cite[\PP{440}{446}]{boston-studies-3}. Also in
  \cite[\PP{11}{16}]{heijenoort-selected-essays}.

\bibitem[\protect\citeauthoryear{Heijenoort{\protectedheijenoortindex}}{1968}]%
{heijenoort-tree-herbrand}
Jean~van Heijenoort{\protectedheijenoortindex}.
\newblock On the relation between the falsifiability tree method and the
  {\herbrand} method in quantification theory.
\newblock Unpublished typescript, \Nov\,20, 1968, \PPcount{12};
  \heijenoortname\ Papers, 1946--1988, Archives of American Mathematics, Center
  for American History, \unitexasaustin, Box\,3.8/86-33/1. \ Copy in
  \anellisarchives, 1968.

\bibitem[\protect\citeauthoryear{Heijenoort{\protectedheijenoortindex}}{1971}]%
{heijenoort-source-book}
Jean~van Heijenoort{\protectedheijenoortindex}.
\newblock {\em From {\frege} to {\goedel}: A Source Book in Mathematical Logic,
  1879--1931}.
\newblock {\harvardunipress}, 1971.
\newblock {\nth 2\,\rev\ \edn\ (\nth 1\,\edn\,1967)}.

\bibitem[\protect\citeauthoryear{Heijenoort{\protectedheijenoortindex}}{1975}]%
{heijenoort-herbrand}
Jean~van Heijenoort{\protectedheijenoortindex}.
\newblock \herbrand.
\newblock Unpublished typescript, \May\,18, 1975, \PPcount{15};
  \heijenoortname\ Papers, 1946--1988, Archives of American Mathematics, Center
  for American History, \unitexasaustin, Box\,3.8/86-33/1. \ Copy in
  \anellisarchives, 1975.

\bibitem[\protect\citeauthoryear{Heijenoort{\protectedheijenoortindex}}{1976}]%
{heijenoort-spanish-lecture-course}
Jean~van Heijenoort{\protectedheijenoortindex}.
\newblock {\em El Desarrollo de la Teoria de la Cuantifcaci\'on}.
\newblock Instituto de Investigaciones Filos\'oficas, Universidad Nacional
  Aut\'onoma de M\'exico, 1976.

\bibitem[\protect\citeauthoryear{Heijenoort{\protectedheijenoortindex}}{1982}]%
{heijenoort-oeuvre-herbrand}
Jean~van Heijenoort{\protectedheijenoortindex}.
\newblock {L'\oefranz uvre} logique de {\herbrand} et son contexte historique.
\newblock 1982.
\newblock In \cite[\PP{57}{85}]{herbrand-symposium}. \Rev\ English translation
  is \cite{heijenoort-work-herbrand}.

\bibitem[\protect\citeauthoryear{Heijenoort{\protectedheijenoortindex}}{1986a}%
]{heijenoort-engels}
Jean~van Heijenoort{\protectedheijenoortindex}.
\newblock {\engelsname} and mathematics.
\newblock 1986.
\newblock In \cite[\PP{123}{151}]{heijenoort-selected-essays}. Previously
  unpublished manuscript written in~1948.

\bibitem[\protect\citeauthoryear{Heijenoort{\protectedheijenoortindex}}{1986b}%
]{heijenoort-work-herbrand}
Jean~van Heijenoort{\protectedheijenoortindex}.
\newblock {\herbrand's} work in logic and its historical context.
\newblock 1986.
\newblock In \cite[\PP{99}{121}]{heijenoort-selected-essays}. \Rev\ English
  translation of \cite{heijenoort-oeuvre-herbrand}.

\bibitem[\protect\citeauthoryear{Heijenoort{\protectedheijenoortindex}}{1986c}%
]{heijenoort-selected-essays}
Jean~van Heijenoort{\protectedheijenoortindex}.
\newblock {\em Selected Essays}.
\newblock Bibliopolis, Napoli, copyright 1985. Also published by Librairie
  Philosophique J. Vrin, \Paris, 1986, 1986.

\bibitem[\protect\citeauthoryear{Heijenoort{\protectedheijenoortindex}}{1992}]%
{heijenoort-modern-logic}
Jean~van Heijenoort{\protectedheijenoortindex}.
\newblock Historical development of modern logic.
\newblock {\em Modern Logic}, 2:242--255, 1992.
\newblock Written in 1974.

\bibitem[\protect\citeauthoryear{Herbran{\protectedherbrandindex}d}{1930}]{her%
brand-PhD}
Jacques Herbran{\protectedherbrandindex}d.
\newblock {\em {\herbrandPhDtitle}}.
\newblock PhD thesis, Universit\'e de Paris, 1930.
\newblock {Th\`eses} pr\'esent\'ees \`a la facult\'e des Sciences de Paris pour
  obtenir le grade de docteur\`es sciences math\'ematiques --- \frenchnth 1
  \nolinebreak th\`ese: Recherches sur la th\'eorie de la d\'emonstration ---
  \frenchnth 2 \nolinebreak th\`ese: Propositions donn\'ees par la facult\'e,
  Les \'equations de Fredholm --- Soutenues le 1930 devant la commission
  d'examen --- Pr\'esident: M. \vessiot, Examinateurs: MM. {\namefont Denjoy,
  Frechet} --- Vu et approuv\'e, \Paris, le 20\,Juin\,1929, Le doyen de la
  facult\'e des Sciences, {\namefont C. Maurain} --- Vu et permis d'imprimer,
  \Paris, le 20\,Juin\,1929, Le recteur de l'Academie de \Paris, {\namefont S.
  Charlety} --- No.\,d'ordre\,2121, S\'erie\,A, No.\,de\,S\'erie\,1252 ---
  Imprimerie J. Dziewulski, Varsovie --- \Univ\ de \Paris. \ Also in Prace
  \WarsawScientificSocietyIII\ Nauk Matematyczno-Fizychnych, Nr.\,33,
  \Warszawa. \ A contorted, newly typeset reprint is
  \cite[\PP{35}{153}]{herbrand-ecrits-logiques}. \ Annotated {E}nglish
  translation {\em Investigations in Proof Theory}\/ by \goldfarbname\
  (\litchapfromtoref 1 4) and \drebenname\ and \heijenoortname\ (\litchapref 5)
  with a brief introduction by \goldfarb\ and extended notes by \goldfarb\
  \mbox{(\litnotefromtoref A C, K--M, O)}, \dreben\ (\litnotefromtoref F I),
  \dreben\ and \goldfarb\ (\litnoterefss D J N), and \dreben, \huffname, and
  \hailperinname\ (\litnoteref E) in
  \cite[\PP{44}{202}]{herbrand-logical-writings}. \ {E}nglish translation of
  \litsectref 5 with a different introduction by \heijenoort\ and some
  additional extended notes by \dreben\ also in
  \cite[\PP{525}{581}]{heijenoort-source-book}. \ \bibcrit{\herbrand's
  \PhDthesis, his cardinal work, dated \herbrandPhDdating; submitted at the
  \Univ\ of \Paris; defended at the Sorbonne \herbrandPhDdefensedate; printed
  in \WarszawaEnglish, 1930.}

\bibitem[\protect\citeauthoryear{Herbran{\protectedherbrandindex}d}{1936}]{her%
brand-1936}
Jacques Herbran{\protectedherbrandindex}d.
\newblock {\em Le D\'eveloppement Moderne de la Th\'eorie des Corps
  Alg\'ebriques --- Corps de classes et lois de r\'eciprocit\'e}.
\newblock M\'emorial des Sciences Math\'ematiques, Fascicule\,LXXV.
  {\gauthiervillars}, 1936.
\newblock \Ed\ and with an appendix by \chevalleyname.

\bibitem[\protect\citeauthoryear{Herbran{\protectedherbrandindex}d}{1968}]{her%
brand-ecrits-logiques}
Jacques Herbran{\protectedherbrandindex}d.
\newblock {\em {\'E}crits Logiques}.
\newblock Presses Universitaires de France, Paris, 1968.
\newblock Cortorted \edn~of \herbrand's logical writings by \heijenoortname.
  Review in \cite{review-herbrand-ecrits-logiques}. \hskip.2em {E}nglish
  translation is \cite{herbrand-logical-writings}.

\bibitem[\protect\citeauthoryear{Herbran{\protectedherbrandindex}d}{1971}]{her%
brand-logical-writings}
Jacques Herbran{\protectedherbrandindex}d.
\newblock {\em Logical Writings}.
\newblock {\harvardunipress}, 1971.
\newblock \Ed\ by \goldfarbname. Translation of \cite{herbrand-ecrits-logiques}
  with additional annotations, brief introductions, and extended notes by
  \goldfarb, \drebenname, and \heijenoortname. \
  \bibcrit{\textonherbrandlogicalwritings}.

\bibitem[\protect\citeauthoryear{Hilbe{\protectedhilbertindex}rt \bgroup\&\
  \egroup
  Bernays\protect\bernaysindex}{1934}]{grundlagen-first-edition-volume-one}
David Hilbe{\protectedhilbertindex}rt and Paul Bernays\protect\bernaysindex.
\newblock {\em {Die Grundlagen der Mathematik --- Erster Band}}.
\newblock Number~XL in {Die Grundlehren der Mathematischen Wissenschaften in
  Einzeldarstellungen}. {\springerverlag}, 1934.
\newblock \nth 1\,\edn\ (\nth 2\,\edn\ is
  \cite{grundlagen-second-edition-volume-one}). English translation is
  \makeaciteoftwo{grundlagen-german-english-edition-volume-one-one}{grundlagen%
-german-english-edition-volume-one-two}.

\bibitem[\protect\citeauthoryear{Hilbe{\protectedhilbertindex}rt \bgroup\&\
  \egroup
  Bernays\protect\bernaysindex}{1939}]{grundlagen-first-edition-volume-two}
David Hilbe{\protectedhilbertindex}rt and Paul Bernays\protect\bernaysindex.
\newblock {\em {Die Grundlagen der Mathematik --- Zweiter Band}}.
\newblock Number~L in {Die Grundlehren der Mathematischen Wissenschaften in
  Einzeldarstellungen}. {\springerverlag}, 1939.
\newblock \nth 1\,\edn\ (\nth 2\,\edn\ is
  \cite{grundlagen-second-edition-volume-two}).

\bibitem[\protect\citeauthoryear{Hilbe{\protectedhilbertindex}rt \bgroup\&\
  \egroup
  Bernays\protect\bernaysindex}{1968}]{grundlagen-second-edition-volume-one}
David Hilbe{\protectedhilbertindex}rt and Paul Bernays\protect\bernaysindex.
\newblock {\em {Die Grundlagen der Mathematik~I}}.
\newblock Number~40 in {Die Grundlehren der Mathematischen Wissenschaften in
  Einzeldarstellungen}. {\springerverlag}, 1968.
\newblock \nth 2\,\rev\,\edn\,of \cite{grundlagen-first-edition-volume-one}.
  English translation is
  \makeaciteoftwo{grundlagen-german-english-edition-volume-one-one}{grundlagen%
-german-english-edition-volume-one-two}.

\bibitem[\protect\citeauthoryear{Hilbe{\protectedhilbertindex}rt \bgroup\&\
  \egroup
  Bernays\protect\bernaysindex}{1970}]{grundlagen-second-edition-volume-two}
David Hilbe{\protectedhilbertindex}rt and Paul Bernays\protect\bernaysindex.
\newblock {\em {Die Grundlagen der Mathematik~II}}.
\newblock Number~50 in {Die Grundlehren der Mathematischen Wissenschaften in
  Einzeldarstellungen}. {\springerverlag}, 1970.
\newblock \nth 2\,\rev\,\edn\,of \cite{grundlagen-first-edition-volume-two}.

\bibitem[\protect\citeauthoryear{Hilbe{\protectedhilbertindex}rt \bgroup\&\
  \egroup
  Bernays\protect\bernaysindex}{2013a}]{grundlagen-german-english-edition-volu%
me-one-one}
David Hilbe{\protectedhilbertindex}rt and Paul Bernays\protect\bernaysindex.
\newblock {\em {Grundlagen der Mathematik~I --- Foundations of Mathematics~I,
  Part\,A: \whatisinVolIPartA}}.
\newblock \url{http://wirth.bplaced.net/p/hilbertbernays}, 2013.
\newblock Thoroughly \rev\,\nth 2\,\edn\ (\nth 1\,\edn\ \collegepublications,
  2011). First English translation and bilingual facsimile \edn\ of the \nth 2
  German \edn\ \cite{grundlagen-second-edition-volume-one}, \incl\ the
  annotation and translation of all differences of the \nth 1 German \edn\
  \cite{grundlagen-first-edition-volume-one}. Translated and commented by
  \wirthname. \Ed\ by \wirthname, \siekmannname, \michaelgabbayname,
  \gabbayname. Advisory Board: \siegname\ (chair), \anellisname, \awodeyname,
  \baazname, \buchholzname, \buldtname, \kahlename, \mancosuname, \parsonsname,
  \peckhausname, \taitname, \tappname, \zachname.

\bibitem[\protect\citeauthoryear{Hilbe{\protectedhilbertindex}rt \bgroup\&\
  \egroup
  Bernays\protect\bernaysindex}{2013b}]{grundlagen-german-english-edition-volu%
me-one-two}
David Hilbe{\protectedhilbertindex}rt and Paul Bernays\protect\bernaysindex.
\newblock {\em {Grundlagen der Mathematik~I --- Foundations of Mathematics~I,
  Part\,B: \whatisinVolIPartB}}.
\newblock \url{http://wirth.bplaced.net/p/hilbertbernays}, 2013.
\newblock Thoroughly \rev\,\nth 2\,\edn. First English translation and
  bilingual facsimile \edn\ of the \nth 2 German \edn\
  \cite{grundlagen-second-edition-volume-one}, \incl\ the annotation and
  translation of all deleted texts of the \nth 1 German \edn\
  \cite{grundlagen-first-edition-volume-one}. Translated and commented by
  \wirthname. \Ed\ by \wirthname, \siekmannname, \michaelgabbayname,
  \gabbayname. Advisory Board: \siegname\ (chair), \anellisname, \awodeyname,
  \baazname, \buchholzname, \buldtname, \kahlename, \mancosuname, \parsonsname,
  \peckhausname, \taitname, \tappname, \zachname.

\bibitem[\protect\citeauthoryear{Hintikka}{1996}]{hintikkaprinciples}
K.~Jaakko~J. Hintikka.
\newblock {\em The Principles of Mathematics Revisited}.
\newblock {\cambridgeunipress}, 1996.

\bibitem[\protect\citeauthoryear{Kohlhase}{2006}]{forthMKMfive}
Michael Kohlhase, editor.
\newblock {\em {\theforthMKMfive}, Revised Selected Papers}, number 3863 in
  Lecture Notes in Artificial Intelligence. {\springerverlag}, 2006.

\bibitem[\protect\citeauthoryear{Kuhn}{1962}]{kuhn}
Thomas~S. Kuhn.
\newblock {\em The Structure of Scientific Revolutions}.
\newblock \chicagounipress, 1962.
\newblock \nth 1\,\edn.

\bibitem[\protect\citeauthoryear{L{\"o}wenheim{\protectedloewenheimindex}}{191%
5}]{loewenheim-1915}
Leopold L{\"o}wenheim{\protectedloewenheimindex}.
\newblock {\Ue ber M\oe glich\-keiten im Relativ\-kalk\ue l}.
\newblock {\em {\MathematischeAnnalen}}, 76:228--251, 1915.
\newblock {English} translation {\em On Possibilities in the Calculus of
  Relatives} by \bauermengelbergname\ with an introduction by \heijenoortname\
  in \cite[\PP{228}{251}]{heijenoort-source-book}.

\bibitem[\protect\citeauthoryear{Martelli \bgroup\&\ \egroup
  Montanari}{1982}]{martelli-montanari}
Alberto Martelli and Ugo Montanari.
\newblock An efficient unification algorithm.
\newblock {\em \ACM\ Transactions on Programming Languages and Systems},
  4:258--282, 1982.

\bibitem[\protect\citeauthoryear{Menzler-Trott}{2001}]{menzler-gentzen-german}
Eckart Menzler-Trott.
\newblock {\em {\gentzen's Problem -- Mathematische Logik im
  nationalsozialistischen Deutschland}}.
\newblock \birkhaeuser, 2001.
\newblock \Rev\ English translation is \cite{menzler-gentzen-english}.

\bibitem[\protect\citeauthoryear{Menzler-Trott}{2007}]{menzler-gentzen-english}
Eckart Menzler-Trott.
\newblock {\em Logic's Lost Genius --- The Life {of\/ \gentzenname}}.
\newblock \AMSname, 2007.
\newblock \Rev\ English translation of \cite{menzler-gentzen-german}.

\bibitem[\protect\citeauthoryear{Nieuwenhuis}{2005}]{twentiethCADEfive}
Robert Nieuwenhuis, editor.
\newblock {\em {\thetwentiethCADEfive}}, number 3632 in Lecture Notes in
  Artificial Intelligence. {\springerverlag}, 2005.

\bibitem[\protect\citeauthoryear{Paterson \bgroup\&\ \egroup
  Wegman}{1978}]{PatersonWegman78}
Michael~S. Paterson and Mark~N. Wegman.
\newblock Linear unification.
\newblock {\em \jcssname}, 16:158--167, 1978.

\bibitem[\protect\citeauthoryear{Peckhau{\protectedpeckhausindex}s}{2004}]{sch%
roeder-handbook}
Volker Peckhau{\protectedpeckhausindex}s.
\newblock {\schroeder's} logic.
\newblock 2004.
\newblock {\protectedschroederindex}In \cite[\Vol\,3: The Rise of Modern Logic:
  From \leibniz\ to \frege, \PP{557}{610}]{handbook-of-the-history-of-logic}.

\bibitem[\protect\citeauthoryear{Peirce\protectedpeirceindex}{1885}]{peirce-18%
85}
Charles~S. Peirce\protectedpeirceindex.
\newblock On the algebra of logic: A contribution to the philosophy of
  notation.
\newblock {\em American J. of Mathematics}, 7:180--202, 1885.
\newblock Also in \cite[\PP{162}{190}]{peirce-chrono-five}.

\bibitem[\protect\citeauthoryear{Peirce\protectedpeirceindex}{1993}]{peirce-ch%
rono-five}
Charles~S. Peirce\protectedpeirceindex.
\newblock {\em {Writings of\/ \peircename\ --- A Chronological Edition,
  \Vol\,5, 1884--1886}}.
\newblock \indianaunipress, 1993.
\newblock \Ed\ by {\namefont Christian J. W. Kloesel}.

\bibitem[\protect\citeauthoryear{Schr{\"o}der{\protectedschroederindex}}{1895}%
]{schroeder-vorlesungen-III}
Ernst Schr{\"o}der{\protectedschroederindex}.
\newblock {\em {Vorlesungen \ue ber die Algebra der Logik, \Vol\,3, Algebra der
  Logik und der Relative, Vorlesungen I-XII}}.
\newblock \teubnerverlag, \Leipzig, 1895.
\newblock {E}nglish translation of some parts in \cite{brady}.

\bibitem[\protect\citeauthoryear{Sch{\"u}tte{\protectedschuetteindex}}{1960}]{%
schuette60:_beweis}
Kurt Sch{\"u}tte{\protectedschuetteindex}.
\newblock {\em {Bewei\esi theorie}}.
\newblock Number 103 in Grundlehren der mathematischen Wissenschaften.
  {\springerverlag}, 1960.
\newblock Thoroughly revised {English} translation: \cite{schuette-1977}.

\bibitem[\protect\citeauthoryear{Sch{\"u}tte{\protectedschuetteindex}}{1977}]{%
schuette-1977}
Kurt Sch{\"u}tte{\protectedschuetteindex}.
\newblock {\em Proof theory}.
\newblock Number 225 in Grundlehren der mathematischen Wissenschaften.
  {\springerverlag}, 1977.
\newblock Translated from a thorough revision of \cite{schuette60:_beweis} by
  \crossleyname.

\bibitem[\protect\citeauthoryear{Skolem{\protectedskolemindex}}{1920}]{skolem-%
1920}
Thoralf Skolem{\protectedskolemindex}.
\newblock {Logisch-kombinatorische Untersuchungen \ue ber die Erf\ue llbarkeit
  und Bewei\esi barkeit mathematischer S\ae tze nebst einem Theorem \ue ber
  dichte Mengen}.
\newblock {\em {\skriftername, \skrifterpublisher}}, 1920/4:1--36, 1920.
\newblock Also in \cite[\PP{103}{136}]{skolem-1970}. \ {E}nglish translation of
  \litsectref 1 {\em Logico-Combinatorial Investigations in the Satisfiability
  or Provability of Mathematical Propositions: A simplified proof of a theorem
  by \loewenheimname\ and generalizations of the theorem} by
  \bauermengelbergname\ with an introduction by \heijenoortname\ in
  \cite[\PP{252}{263}]{heijenoort-source-book}. \
  \bibcrit{\loewenheimskolemtheorem\ via \skolemnormalform\ and choice of a
  sub-model}.

\bibitem[\protect\citeauthoryear{Skolem{\protectedskolemindex}}{1923}]{skolem-%
1923b}
Thoralf Skolem{\protectedskolemindex}.
\newblock {Einige Bemerkungen zur axiomatischen Begr\ue ndung der Mengenlehre}.
\newblock In {\em \Proc\ \nth 5 Scandinaviska Matematikerkongressen,
  {H}elsingfors, July\,4--7,\,1922}, pages 217--232, Helsingfors, 1923.
  Akademiska Bokhandeln.
\newblock Also in \cite[\PP{137}{152}]{skolem-1970}. \ {E}nglish translation
  {\em Some remarks on Axiomatized Set Theory} by \bauermengelbergname\ with an
  introduction by \heijenoortname\ in
  \cite[\PP{290}{301}]{heijenoort-source-book}. \ \bibcrit{The best written
  and, together with \cite{skolem-1920}, the most relevant publication on the
  \loewenheimskolemtheorem; although there is some minor gap in the proof
  according to \citet{wang-skolem}, cured in \cite{skolem-1929}. \ Includes
  \skolemsparadox\ and a proof of the \loewenheimskolemtheorem\ which does not
  require any weak forms of the \axiomofchoice, via \skolemnormalform\ and
  construction of a model without assuming the previous existence of another
  one}.

\bibitem[\protect\citeauthoryear{Skolem{\protectedskolemindex}}{1928}]{skolem-%
1928}
Thoralf Skolem{\protectedskolemindex}.
\newblock {\Ue ber die mathematische Logik (Nach einem Vortrag gehalten im
  Norwegischen Mathematischen Verein am {22.\,Oktober\,1928})}.
\newblock {\em Nordisk Matematisk Tidskrift}, 10:125--142, 1928.
\newblock Also in \cite[\PP{189}{206}]{skolem-1970}. \ {E}nglish translation
  {\em On Mathematical Logic}\/ by \bauermengelbergname\ and \foellesdalname\
  with an introduction by \drebenname\ and \heijenoortname\ in
  \cite[\PP{508}{524}]{heijenoort-source-book}. \ \bibcrit{First explicit
  occurrence of \skolemization\ and \skolem\ functions}.

\bibitem[\protect\citeauthoryear{Skolem{\protectedskolemindex}}{1929}]{skolem-%
1929}
Thoralf Skolem{\protectedskolemindex}.
\newblock {\Ue ber einige Grundlagenfragen der Mathematik}.
\newblock {\em {\skriftername, \skrifterpublisher}}, 1929/4:1--49, 1929.
\newblock Also in \cite[\PP{227}{273}]{skolem-1970}. \ \bibcrit{Detailed
  discussion of \skolemsparadox\ and the \loewenheimskolemtheorem}.

\bibitem[\protect\citeauthoryear{Skolem{\protectedskolemindex}}{1970}]{skolem-%
1970}
Thoralf Skolem{\protectedskolemindex}.
\newblock {\em Selected Works in Logic}.
\newblock Universitetsforlaget \Oslo, 1970.
\newblock \Ed\ by \fenstadname. \bibcrit{Without index, but with most funny
  spellings in the newly set titles}.

\bibitem[\protect\citeauthoryear{Smullyan}{1968}]{smullyan}
Raymond~M. Smullyan.
\newblock {\em First-Order Logic}.
\newblock {\springerverlag}, 1968.

\bibitem[\protect\citeauthoryear{Stern}{1982}]{herbrand-symposium}
Jacques Stern, editor.
\newblock {\em {\Proc\ of the \herbrand\ Symposium, Logic Colloquium'81,
  Marseilles, France, \Jul\,1981}}.
\newblock {\northholland}, 1982.

\bibitem[\protect\citeauthoryear{Tai{\protectedtaitindex}t}{2006}]{tait-2006}
William~W. Tai{\protectedtaitindex}t.
\newblock {\goedel's} correspondence on proof theory and constructive
  mathematics.
\newblock {\em {\philosphiamathname}}, 14:76--111, 2006.

\bibitem[\protect\citeauthoryear{Taylor \bgroup\&\ \egroup
  Wiles}{1995}]{wiles-fermat-correction}
Richard Taylor and Andrew Wiles.
\newblock Ring theoretic properties of certain {\hecke} algebras.
\newblock {\em Annals of Mathematics}, 141:553--572, 1995.
\newblock Received \Oct\,7, 1994. Appendix due to \faltingsname\ received
  \Jan\,26, 1995.

\bibitem[\protect\citeauthoryear{Wallen}{1990}]{wallen}
Lincoln~A. Wallen.
\newblock {\em Automated Proof Search in Non-Classical Logics --- efficient
  matrix proof methods for modal and intuitionistic logics}.
\newblock \mitpress, 1990.
\newblock Phd thesis.

\bibitem[\protect\citeauthoryear{Wang{\protectedwangindex}}{1970}]{wang-skolem}
Hao Wang{\protectedwangindex}.
\newblock A survey of {\skolem's} work in logic.
\newblock 1970.
\newblock {\protectedskolemindex}In \cite[\PP{17}{52}]{skolem-1970}.

\bibitem[\protect\citeauthoryear{Whitehead{\protect\whiteheadindex} \bgroup\&\
  \egroup Russe{\protectedrussellindex}ll}{1910--1913}]{PM}
Alfred~North Whitehead{\protect\whiteheadindex} and Bertrand
  Russe{\protectedrussellindex}ll.
\newblock {\em Principia Mathematica}.
\newblock {\cambridgeunipress}, 1910--1913.
\newblock {\nth 1\,\edn}.

\bibitem[\protect\citeauthoryear{Wiles}{1995}]{wiles-fermat}
Andrew Wiles.
\newblock Modular elliptic curves and {{\fermat's Last Theorem}}.
\newblock {\em Annals of Mathematics}, 141:443--551, 1995.
\newblock Received \Oct\,14, 1994.

\bibitem[\protect\citeauthoryear{Wirth{\protectedwirthindex}
  \bgroup\&al.\egroup }{2009}]{herbrand-handbook}
Claus-Peter Wirth{\protectedwirthindex}, J{\"o}rg Siek{\-}mann,
  {\mbox{Ch}}ristoph Benz{\-}m{\ue}ller{\protectedbenzmuellerindex}, and Serge
  Autexier{\protectedautexierindex}.
\newblock \herbrandname: Life, logic, and automated deduction.
\newblock 2009.
\newblock In \cite[\Vol\,5: Logic from \russell\ to \church,
  \PP{195}{254}]{handbook-of-the-history-of-logic}.

\bibitem[\protect\citeauthoryear{Wirth{\protectedwirthindex}
  \bgroup\&al.\egroup }{2014}]{SR--2009--01}
Claus-Peter Wirth{\protectedwirthindex}, J{\"o}rg Siek{\-}mann,
  {\mbox{Ch}}ristoph Benz{\-}m{\ue}ller{\protectedbenzmuellerindex}, and Serge
  Autexier{\protectedautexierindex}.
\newblock {\em Lectures on {\herbrandname} as a Logician}.
\newblock {SEKI-Report SR--2009--01 (ISSN 1437--4447)}. {SEKI Publications},
  {DFKI Bremen GmbH, Safe and Secure Cognitive Systems, Cartesium, Enrique
  Schmidt Str.\,5, D--28359 Bremen, Germany}, 2014.
\newblock \Rev\,\edn\ \May\,2014, \PPcount{ii+82},
  {\url{http://arxiv.org/abs/0902.4682}}.

\bibitem[\protect\citeauthoryear{Wirth{\protectedwirthindex}}{2004}]{wirthcard%
inal}
Claus-Peter Wirth{\protectedwirthindex}.
\newblock {\DescenteInfinie\ + Deduction}.
\newblock {\em {\ljigplname}}, 12:1--96, 2004.
\newblock \www\url{/p/d}.

\bibitem[\protect\citeauthoryear{Wirth{\protectedwirthindex}}{2006}]{nonpermut}
Claus-Peter Wirth{\protectedwirthindex}.
\newblock {\em {$\lim$$+$}, {$\delta^+$}, and Non-Permutability of
  {$\beta$}-Steps}.
\newblock {SEKI-Report SR--2005--01 (ISSN 1437--4447)}. {SEKI Publications},
  Saarland Univ., 2006.
\newblock \Rev\,\edn\ \Jul\,2006 (\nth 1\,\edn\,2005), \PPcount{ii+36},
  \url{http://arxiv.org/abs/0902.3635}. Thoroughly improved version is
  \cite{wirth-jsc-non-permut}.

\bibitem[\protect\citeauthoryear{Wirth{\protectedwirthindex}}{2008}]{wirth-jal}
Claus-Peter Wirth{\protectedwirthindex}.
\newblock \hilbert's epsilon as an operator of indefinite committed choice.
\newblock {\em J. Applied Logic}, 6:287--317, 2008.
\newblock \url{http://dx.doi.org/10.1016/j.jal.2007.07.009}.

\bibitem[\protect\citeauthoryear{Wirth{\protectedwirthindex}}{2012a}]{wirth-he%
ijenoort}
Claus-Peter Wirth{\protectedwirthindex}.
\newblock {\herbrandsfundamentaltheorem} in the eyes of {\heijenoortname}.
\newblock {\em {\logicauniversalisname}}, 6:485--520, 2012.
\newblock Received \Jan\,12, 2012. Published online \Jun\,22, 2012,
  \url{http://dx.doi.org/10.1007/s11787-012-0056-7}.

\bibitem[\protect\citeauthoryear{Wirth{\protectedwirthindex}}{2012b}]{wirth-js%
c-non-permut}
Claus-Peter Wirth{\protectedwirthindex}.
\newblock {\math{\lim\tight+}, \math{\delta^+}, and \NonPermutability\ of
  \math\beta-Steps}.
\newblock {\em {\jscname}}, 47:1109--1135, 2012.
\newblock Received \Jan\,18, 2011. Published online \Jul\,15, 2011,
  \url{http://dx.doi.org/10.1016/j.jsc.2011.12.035}. More funny version is
  \cite{nonpermut}.

\bibitem[\protect\citeauthoryear{Wirth{\protectedwirthindex}}{2012c}]{wirth-ma%
nifesto-ljigpl}
Claus-Peter Wirth{\protectedwirthindex}.
\newblock Human-oriented inductive theorem proving by descente infinie --- {a
  \nolinebreak manifesto}.
\newblock {\em {\ljigplname}}, 20:1046--1063, 2012.
\newblock Received \Jul\,11, 2011. Published online \Mar\,12, 2012,
  \url{http://dx.doi.org/10.1093/jigpal/jzr048}.

\bibitem[\protect\citeauthoryear{Wirth{\protectedwirthindex}}{2012d}]{wirth-hi%
lbert-seki}
Claus-Peter Wirth{\protectedwirthindex}.
\newblock {\em {\hilbert's} epsilon as an Operator of Indefinite Committed
  Choice}.
\newblock {SEKI-Report SR--2006--02 (ISSN 1437--4447)}. {SEKI Publications},
  Saarland Univ., 2012.
\newblock \Rev\,\edn\ \Jan\,2012, \PPcount{ii+73},
  {\url{http://arxiv.org/abs/0902.3749}}.

\bibitem[\protect\citeauthoryear{Wirth{\protectedwirthindex}}{2013}]{SR--2011-%
-01}
Claus-Peter Wirth{\protectedwirthindex}.
\newblock {\em {A Simplified and Improved Free-Variable Framework for
  \hilbert's epsilon as an Operator of Indefinite Committed Choice}}.
\newblock {SEKI Report SR--2011--01 (ISSN 1437--4447)}. {SEKI Publications},
  {DFKI Bremen GmbH, Safe and Secure Cognitive Systems, Cartesium, Enrique
  Schmidt Str.\,5, D--28359 Bremen, Germany}, 2013.
\newblock \Rev\,\edn\ \Jan\,2013 (\nth 1\,\edn\,2011), \PPcount{ii+65},
  {\url{http://arxiv.org/abs/1104.2444}}.

\end{thebibliography}
\vfill\cleardoublepage
\addcontentsline{toc}{section}{Index}
\index{Dreben's correction|see{G{\"o}del's and Dreben's correction}}
\printindex
\end{document}